\pdfoutput=1

\documentclass[10pt]{article}

\usepackage[includeheadfoot,top=2cm,bottom=2cm,left=2cm,right=2cm]{geometry}

\usepackage[utf8]{inputenc} 
\usepackage[T1]{fontenc}

\usepackage{hyperref}      
\usepackage{url}  
\hypersetup{
    colorlinks,
    linkcolor={blue!},
    citecolor={blue!},
    urlcolor={blue!}
}

\usepackage{amsfonts}   
\usepackage{amssymb}
\usepackage{amsthm}
\usepackage{mathtools}
\usepackage{upgreek}
\usepackage{multirow}

\usepackage[nocompress]{cite}
\usepackage[noabbrev,capitalise]{cleveref}
\Crefname{proposition}{Proposition}{Propositions}
\crefname{equation}{}{}
\Crefname{equation}{}{}

\Crefname{proposition}{Proposition}{Propositions}
\newtheorem{lemma}{Lemma}[section]
\newtheorem{proposition}[lemma]{Proposition}

\newtheorem{theorem}[lemma]{Theorem}
\newtheorem{remark}[lemma]{Remark}

\crefname{equation}{}{}
\Crefname{equation}{}{}
\numberwithin{equation}{section}

\usepackage{algorithm}
\usepackage{algorithmicx}
\usepackage{algpseudocode}
\usepackage{longtable}

\usepackage{graphicx}
\usepackage{caption}
\usepackage{subcaption}
\usepackage{float}
\usepackage{booktabs}
\usepackage{makecell}
\usepackage{float}

\usepackage{tikz}
\usetikzlibrary{datavisualization, plotmarks}
\pgfdvdeclarestylesheet{cust}{
    default style/.style={visualizer color=black},
    1/.style={visualizer color=red},
    2/.style={visualizer color=black},
    3/.style={visualizer color=red},
    4/.style={visualizer color=black}
}

\usepackage{appendix}
\usepackage{booktabs}       
\usepackage{nicefrac}       
\usepackage{microtype}      
\usepackage{xcolor}       
\usepackage{lipsum}  

\newcommand{\cond}{\textup{cond}}

\newcommand{\mapsfrom}{\mathrel{\reflectbox{$\mapsto$}}}

\newcommand{{{\bxi}}}{\mathbf{\upxi}}

\newcommand{\bb}{\mathbf{b}}
\newcommand{\bu}{\mathbf{u}}

\newcommand{\bv}{\mathbf{v}}
\newcommand{\by}{\mathbf{y}}
\newcommand{\bt}{\mathbf{t}}
\newcommand{\bm}{\mathbf{m}}

\newcommand{\bz}{\mathbf{z}}
\newcommand{\bx}{\mathbf{x}}

\newcommand{\bw}{\mathbf{w}}

\newcommand{\br}{\mathbf{r}}
\newcommand{\bnull}{\mathbf{0}}
\newcommand{\bTheta}{\mathbf{\Theta}}
\newcommand{\bSigma}{\mathbf{\Sigma}}

\newcommand{\bR}{\mathbf{R}}
\newcommand{\bA}{\mathbf{A}}
\newcommand{\bE}{\mathbf{E}}
\newcommand{\bN}{\mathbf{N}}
\newcommand{\bB}{\mathbf{B}}

\newcommand{\bX}{\mathbf{X}}
\newcommand{\bZ}{\mathbf{Z}}
\newcommand{\bT}{\mathbf{T}}
\newcommand{\bM}{\mathbf{M}}
\newcommand{\bK}{\mathbf{K}}

\newcommand{\bD}{\mathbf{D}}
\newcommand{\bG}{\mathbf{G}}
\newcommand{\bW}{\mathbf{W}}
\newcommand{\bOmega}{\mathbf{\Omega}}

\newcommand{\bS}{\mathbf{S}}

\newcommand{\bQ}{\mathbf{Q}}
\newcommand{\bC}{\mathbf{C}}
\newcommand{\bP}{\mathbf{P}}
\newcommand{\bI}{\mathbf{I}}

\newcommand{\bU}{\mathbf{U}}

\newcommand{\bV}{\mathbf{V}}
\newcommand{\bF}{\mathbf{F}}
\newcommand{\bPi}{\mathbf{\Pi}}

\newcommand{\bff}{\mathbf{f}}

\def\BibTeX{{\rm B\kern-.05em{\sc i\kern-.025em b}\kern-.08em
    T\kern-.1667em\lower.7ex\hbox{E}\kern-.125emX}}

\newcommand{\oleg}[1]{{}}
\newcommand{\blue}[1]{\textcolor{black}{#1}}
\newcommand{\footremember}[2]{%
  \footnote{#2}\newcounter{#1}\setcounter{#1}{\value{footnote}}%
}
\newcommand{\footrecall}[1]{%
  \footnotemark[\value{#1}]%
}

    \title{Preconditioning via Randomized Range Deflation (RandRAND)
    }
	\author{%
    Oleg Balabanov\footremember{ICSI}{International Computer Science Institute} \footremember{LBNL}{Lawrence Berkeley National Laboratory} \footremember{UCB}{Department of Statistics, University of California, Berkeley}~,
    Caleb Ju\footremember{cju}{Industrial \& Systems Engineering, Georgia Institute of Technology}~,
    Kaiwen He\footremember{KH}{Department of Computer Science, Purdue University}~,
    Aryaman Jeendgar\footrecall{ICSI} \footremember{TUM}{School of Computation, Information and Technology, Technical University of Munich}~,
    Michael W. Mahoney\footrecall{ICSI} \footrecall{LBNL} \footrecall{UCB}
}	
	\date{}
    
\begin{document}
	\maketitle
	
	\begin{abstract}
    
    We introduce RandRAND, a new class of randomized preconditioning methods for large-scale linear systems. 
    RandRAND deflates the spectrum via efficient orthogonal projections onto random subspaces, without computing eigenpairs or low-rank approximations. 
    This leads to advantages in computational cost and numerical stability. 
    We establish rigorous condition number bounds that depend only weakly on the problem size and that reduce to a small constant when the dimension of the deflated subspace is comparable to the effective spectral dimension.
    RandRAND can be employed without explicit operations with the deflation basis, enabling the effective use of fast randomized transforms. 
    In this setting, the costly explicit basis orthogonalization is bypassed by using fast randomized Q-less QR factorizations or iterative methods for computing orthogonal projections. These strategies balance the cost of constructing RandRAND preconditioners and applying them  within linear solvers, and can ensure robustness to rounding errors.
	\end{abstract}
	
	\begin{keywords}
		preconditioning, deflation, Krylov method, conjugate gradient, minres, randomized numerical linear algebra (RandNLA), sketching, kernel method, least-squares, ridge regression, portfolio optimization, inverse~PDE
	\end{keywords}

	\section{Introduction} 

Recent decades have witnessed the rise of randomization and preconditioning as central tools in computational science and machine learning. 
This line of work exploits the synergy between the two areas to enhance the efficiency and effectiveness of common computational methods for data analysis, machine learning, and beyond. 
    
In this paper, we introduce novel randomized linear solvers which are based on Randomized RANge Deflation (abbreviated ``RandRAND''). 
We focus on a numerical solution of large-scale systems of the~form:
\begin{equation}\label{eq:initsys}
        (\bA + \mu \bI)\bx = \bb,~~~\bA \in \mathbb{R}^{n \times n},
\end{equation}
where the matrix $\bA$ has a {relatively} small {number} of singular values significantly larger than the rest. The parameter $\mu \in \mathbb{R}$ describes a spectral shift.\footnote{It may be set to zero if no shift is applied or if the shift is unknown.} The efficient solution of such systems is crucial in many applications, including:
\begin{itemize}
\item Gaussian processes and kernel ridge regression.
These methods require solving regularized systems~\cref{eq:initsys} with psd $\bA$ and $\mu \geq 0$.
Real-world datasets often yield matrices with rapidly or moderately decaying spectra~\cite{frangella2023randomized,daskalakis2022good,rathore2024have,gp_intro_2}, making them well suited for RandRAND acceleration.
\item Second-order optimization. Newton steps in second-order optimization methods, including interior-point frameworks, often require solving large linear systems~\cite{zhao2022nysadmm,chu2024randomized,kkt_intro_1}.
RandRAND is well suited for this setting, as it applies to both psd and indefinite $\bA$ and only requires access to $\bA$ through matrix–vector products (matvecs).
\item Inverse problems with PDE constraints.
After reduction, these problems yield linear systems represented by inverse discretized PDE operators.
The systems are typically symmetric, often available in sparse factorized form, and exhibit strong spectral decay due to dominant physical modes~\cite{pde_cons_intro_1}. RandRAND’s flexible design enables robust preconditioning of these systems at low complexity and memory cost.
\item Large-scale eigenvalue computations.
Classical shift-and-invert Krylov and power methods require solving sequences of indefinite systems of the form~\cref{eq:initsys}.
Spectral analyses of operators with several extreme eigenvalues, such as neural tangent kernels and Hessians, for instance are crucial for studying memorization, convergence, and generalization in neural networks~\cite{pyhessian,spectrum_intro_1}.
\end{itemize}
In general, RandRAND can be applied to a wide range of scenarios (in both computational science and engineering applications as well as in artificial intelligence and machine learning applications), including both sparse and dense systems, which may be positive definite, indefinite, or non-symmetric.
{It naturally extends to complex-valued shifts and operators. Furthermore, while it is effective for solving a single system~\cref{eq:initsys}, it becomes even more advantageous when applied to sequences of shifted systems with multiple values of $\mu$, thanks to its ability to recycle computations.} 
    
{We also note that the ideas developed in this paper for the effective use of orthogonal projections and fast basis-less Cholesky QR extend beyond solving linear systems, and could, for instance, be applied to improve the numerical robustness of low-rank approximation algorithms. } 
     
\paragraph{Background.}
    Solving the system~\cref{eq:initsys} via a direct factorization typically requires $\mathcal{O}(n^3)$ flops, which in large-scale settings becomes prohibitively expensive. Furthermore, if $\bA$ is given as an implicit operator operated through matrix-vector multiplications, then forming $\bA$ explicitly would cost $T_\mathrm{mv} n^2$ flops and $\mathcal{O}(n^2)$ memory, where $T_\mathrm{mv} \ge n$ denotes the cost of a single matvec with $\bA$. These factors impose even more efficiency constraints on the factorization~methods.

    For these reasons, projection-based iterative methods, such as the celebrated conjugate gradient (CG) and minimal residual (MINRES) algorithms~\cite{paige1975solution,barrett1994templates,saad2003iterative}, are generally regarded as primary methods for solving~\cref{eq:initsys}. 
    The effectiveness of these methods is closely tied to the spectrum of the operator. 
    For instance, for spd systems, CG and MINRES reduce the error by at least a factor of $\mathrm{e}^2$ ($\approx 7.39$) every $\sqrt{\kappa}$ iterations, where $\kappa$ is the condition number defined as the ratio of the largest singular value to the smallest singular value of the operator.

     The convergence of an iterative solver can be accelerated via \textit{preconditioning}. 
     Denote the regularized operator $\bA + \mu \bI$ as $\bA_\mu$. 
     Preconditioning consists in constructing a preconditioner $\bP \in \mathbb{R}^{n \times n}$ such that $\bA_\mu \bP$ (or $\bP^{\frac{1}{2}}\bA_\mu \bP^{\frac{1}{2}}$) has more favorable spectral properties, e.g., it may have a smaller conditioner number. 
     As the result, solving the preconditioned system
     $$ \bA_\mu \bP \, \by = \bb, ~~~ \bx = \bP \, \by$$
     can drastically reduce the number of iterations, compared to solving the original system~\cref{eq:initsys}. 
     
     There are two basic approaches for constructing $\bP$: 
     first, by approximating the inverse of $\bA_\mu$, often using low-rank approximations or incomplete factorizations (see, e.g.,~\cite{frangella2023randomized,derezinski2024solving,benzi2003robust,benzi2002preconditioning} and the references therein); and 
     second, by using a deflation method, improving the spectrum of $\bA_\mu$ by deflation of a subspace from the range of $\bA_\mu$ (see, e.g.,~\cite{coulaud2013deflation,gutknecht2012spectral,gaul2013framework} and the references therein). 
     In this paper, we adopt the latter approach, enhancing its efficiency with Randomized Numerical Linear Algebra (RandNLA) techniques~\cite{mahoney2011randomized,halko2011finding,DM16_CACM,MD16_chapter,randnla_kdd24_TR}, while preserving its effectiveness with a high probability.
     In RandNLA, there has been work on preconditioning with randomized approximation methods (see, e.g.,~\cite{rokhlin2008fast, AMT10, MSM14_SISC, frangella2023randomized, derezinski2020improved}).     
     Notably, this is the first randomized preconditioning method that does not rely explicitly on operator approximation.

    \paragraph{Preliminaries.}   
    For a matrix $\bC \in \mathbb{R}^{n \times n}$, we denote the singular values of $\bC$ by $\sigma_1(\bC)$, $\sigma_2(\bC), \hdots, \sigma_n(\bC)$ and the eigenvalues by $\lambda_1(\bC)$, $\lambda_2(\bC), \hdots, \lambda_n(\bC)$, both ordered in decreasing magnitude. 
    Let $\lambda_{\mathrm{min}}(\bC)$ and $\sigma_{\mathrm{min}}(\bC)$ denote the smallest eigenvalue (if it is real) and singular value of $\bC$, respectively. 
    Define the $k$-stable rank of $\bC$  as $\mathrm{sr}_{k}(\bC):=  \sum_{j\geq k}\frac{\sigma_j(\bC)}{\sigma_k(\bC)}$, {which 
    quantifies the decay of the trailing $n-k+1$ singular values of $\bC$.}
    In addition, define the $k$-stable condition number of $\bC$ with respect to regularization parameter $\mu$ as $\cond_{k, \mu}(\bC):= \frac{1}{n-k+1}  \sum_{j\geq k}\frac{\sigma_j(\bC)}{\sigma_n(\bC)+|\mu|}$.  
    {If $\bC$ is psd and $\mu \geq 0$, this quantity compares the trailing eigenvalues of the unregularized matrix to the smallest eigenvalue of the regularized matrix $\bC+\mu\bI$, thereby quantifying the effectiveness of the deflation of the top $k$ eigenvalues from the spectrum of $\bC$ for preconditioning the regularized operator $\bC+\mu\bI$. } 
    Denote $\cond_{k, 0}(\bC)$ by $\cond_{k}(\bC)$. Furthermore, the eigenvalues and singular values $\lambda_i(\bA)$ and $\sigma_i(\bA)$ of the operator in~\cref{eq:initsys} are simply denoted by $\lambda_i$ and $\sigma_i$. 
    We use $c, \rho, \tau$ to denote user-specified or problem-dependent constants, and $C_i$ with $i=1, 2,\hdots$ to denote universal constants.
    {We call a preconditioner \emph{quasi-optimal} if it is constructed using an $l$-dimensional projection or a rank-$l$ approximation,
    and the resulting preconditioned operator $\bB$ satisfies
    $$
    \sigma_1(\bB) \leq c_2\, \sigma_k(\bA_\mu) \quad \text{and} \quad \sigma_{\min}(\bB) \geq c_3\, \sigma_{\min}(\bA_\mu),
    $$
    for $k \geq c_1 l$ for some coefficients $c_1, c_2, c_3$ bounded independently of the properties of $\bA_\mu$ and with only minor dependence on the dimension $n$.\footnote{In that $1/c_1, c_2, 1/c_3$ remain small, say $\leq 10$, for all reasonable problem dimensions}} Finally, in algorithms we use $\mathrm{qr}(\cdot)$ and $\mathrm{chol}(\cdot)$ to denote the economic QR factors and the Cholesky factor of a matrix, respectively.

    \paragraph{Roadmap.}
    The article is organized as follows. \Cref{our_approach,spd_solvers} introduce RandRAND preconditioning, including a numerical illustration, and a summary of its quasi-optimality and complexity bounds; and \Cref{related_work} discusses connections with prior work, including in particular a new interpretation of the Nystr{\"o}m approximation, which links some RandRAND preconditioners to Nystr{\"o}m methods. \Cref{range_approx} develops results on shifted randomized range approximation, extending classical results from~\cite{halko2011finding,boutsidis2013improved,tropp2011improved,cohen2015optimal,gittens2016revisiting,gittens2013revisiting} to the shifted subspace iteration and providing the foundation for performance guarantees of RandRAND. \Cref{sec:allprecond}, the main part of the paper, presents several variants of RandRAND and establishes rigorous condition number bounds. \Cref{orthogonalization} addresses implementation aspects, including efficient basis construction and numerical robustness. \Cref{adaptations} extends RandRAND to sweeps over parameter $\mu$ with computation recycling, and introduces variants that use (randomized) iterative methods to compute orthogonal projections, thereby avoiding or only approximately computing factorizations. %\Cref{adaptations} also discusses the adaptive selection of the random basis size. 
    \Cref{solvers} integrates RandRAND preconditioners into CG and MINRES linear solvers. \Cref{emp_eval} provides an extensive numerical study of RandRAND on synthetic systems and real-world applications, including kernel and random features regression, portfolio optimization, and PDE-constrained inverse problems. 
    \Cref{appendix_a} contains proofs and extensions of the theoretical results from~\cref{range_approx}; \Cref{appendix_b} contains proofs of the results from~\cref{sec:allprecond}; \Cref{appendix_c} offers additional details on the numerical experiments; and \Cref{appendix_d_0} provides additional characterization of certain RandRAND variants.
    
    \subsection{Our Approach} \label{our_approach}
    
    \paragraph{Framework.}
    A key ingredient of RandRAND is to use a random basis $\bV$ that, with high probability, accurately approximates the extreme eigenvectors or singular vectors of $\bA_\mu$ and allows efficient computation of orthogonal projections onto $\mathrm{range}(\bV)$. 
    The RandRAND procedure for spd systems is summarized in \cref{alg:framework}, where $\tau \geq 0$ is an efficiently computable coefficient.  
    For general systems, the framework in~\cref{alg:framework} can be extended by replacing $\bB$ and $\bA_\mu$ in~\cref{eq:option1} with the normal operators $\bB \bB^{\mathrm{T}}$ (or $\bB^{\mathrm{T}} \bB$)  and $\bA_\mu \bA_\mu^{\mathrm{T}}$ (or $\bA_\mu^{\mathrm{T}} \bA_\mu$), respectively.

    \begin{algorithm}[t] 
    \caption{RandRAND Framework} 
    \label{alg:framework}
    \begin{algorithmic}[1]
        \State Construct a random basis $\bV$ that approximates the extreme eigenvectors or singular vectors of $\bA_\mu$ and allows efficient computation of orthogonal projections onto $\mathrm{range}(\bV)$.
        \State Build a preconditioner for $\bA_\mu$ using orthogonal projections such that the resulting preconditioned operator $\bB$ satisfies
        \begin{equation}\label{eq:option1}
            \langle \bv, \bB \bv \rangle \approx \tau \|\bv\|^2, \quad \forall \bv \in \mathrm{range}(\bV),
            \qquad
            \langle \bv, \bB \bv \rangle \approx \langle \bv, \bA_\mu \bv \rangle, \quad \forall \bv \in \mathrm{range}(\bV)^\perp.
        \end{equation}
        \State Use the preconditioner within an iterative solver for~\cref{eq:initsys}.
    \end{algorithmic}
\end{algorithm}

\paragraph{Randomized Range Approximation.}
We select the basis $\bV$ as $\bA_\mu \bOmega$, where the test matrix $\bOmega \in \mathbb{R}^{n \times l}$, with $l \ll n$, is given by
\begin{equation} \label{eq:testmatrix}
    \bOmega := 
    \begin{cases} 
        \bA^q \bX^\mathrm{T}, & \text{if $\bA$ is symmetric},\\
        (\bA^\mathrm{T} \bA)^q \bX^\mathrm{T}, &  \text{if $\bA$ is non-symmetric,}
    \end{cases}
\end{equation}
where $\bX \in \mathbb{R}^{l \times n}$ is a random matrix from a family of oblivious subspace embeddings (OSEs), and $q$ is a small integer power, say $q = 0$ or $1$.
Prominent OSE operators include Gaussian matrices, randomized sparse operators, and operators based on fast Hadamard and trigonometric transforms~\cite{DMMS07_FastL2_NM10, halko2011finding, cohen2015optimal, tropp2011improved, ailon2009fast}. 
%{We show} that to capture the dominant $k$ singular vectors of $\bA_\mu$, it suffices to use the sketch dimension $l = \mathcal{O}(k)$ and a small power $q=0$ or $q=1$.
Note that $\bOmega$ does not depend on the shift $\mu$.
Moreover, the matrix $\bX$ is drawn independently from $\bA$ and therefore allows RandRAND to be implemented through matvecs $\bu \mapsto \bA\bu$, without direct access to the entries of $\bA$.

{Note also that while we primarily consider construction with OSEs, RandRAND is also compatible with data-aware embeddings~$\bX$, including sampling matrices~\cite{DMM08_CURtheory_JRNL, alaoui2015fast,gittens2016revisiting,gittens2013revisiting,woodruff2014sketching,chen2025randomly,diaz2023robust} and Krylov bases generated with a random starting block~\cite{tropp2023randomized,musco2015randomized}.}

\subsubsection{Preconditioners}

 Let $\bPi$ denote the orthogonal projector onto $\mathrm{range}(\bV)$. We introduce three RandRAND-based preconditioners. 
     \begin{subequations} \label{eq:RandRANDintro}
     \begin{itemize}
         \item R-RandRAND preconditioner for spd systems: 
         \begin{equation}
            \bP = \bA_\mu^{-1} \left (( \bI - \bPi) \bA_\mu (\bI - \bPi) + \tau \bPi \right )
          \end{equation}
          that achieves~\cref{eq:option1} through \textit{reducing} the action of $\bA_\mu$ in $\mathrm{range}({\bV})$. In particular, $\bP$  satisfies $\bA_\mu \bP \bv = \tau \bv$ for all $\bv \in \mathrm{range}({\bV})$. 
         \item C-RandRAND preconditioner  for spd systems:
         \begin{equation}
            \bP = (\bI - \bPi) + \tau \bPi \bA_\mu^{-1} \bPi
          \end{equation}
          that achieves~\cref{eq:option1} through \textit{correcting} the action of $\bA_\mu$ in $\mathrm{range}({\bV})$. In particular, $\bP$ satisfies $\bA_\mu \bP \bv = \tau \bw $, where $\bw = \bA_\mu (\bPi \bA_\mu^{-1}\bv)$, for all $\bv \in \mathrm{range}({\bV})$.      
         \item  G-RandRAND left preconditioner:
          \begin{equation}
            \bP = (\bI - \bPi) + \tau \left  ((\bA_\mu ^{-1}\bPi )^\mathrm{T} \bA_\mu^{-1} \bPi  \right )^{\frac{1}{2}}
          \end{equation}
           that \textit{generalizes} C-RandRAND to preconditioning of indefinite and non-symmetric systems. 
           \end{itemize}
     \end{subequations}
    Clearly, C-RandRAND and G-RandRAND preconditioners $\bP$ are spd, as is the operator $\bA_\mu \bP$ preconditioned by R-RandRAND. 
    
     {In this work} we show that the orthogonal projections in~\cref{eq:RandRANDintro} effectively deflate the range of $\bA_\mu\bOmega$ from the range of $\bA_\mu$. 
     Specifically, with an appropriate choice of the coefficient $\tau$, R-RandRAND provides a preconditioned matrix $\bB = \bA_\mu \bP$ which satisfies  
     \begin{equation} \label{eq:cndbound0}
        \sigma_1(\bB) \leq \|(\bI - \bPi) \bA_\mu \|~~~\text{and}~~~\sigma_n(\bB) \geq \sigma_n(\bA_\mu).
    \end{equation}  
    The preconditioners C-RandRAND and G-RandRAND satisfy a similar guarantee to~\cref{eq:cndbound0} up to constants.  
    From this it follows that RandRAND preconditioning greatly improves the condition number when $\bA_\mu \bOmega$ accurately captures the action of $\bA_\mu$. In turn, we prove that for $\bOmega$ in~\cref{eq:testmatrix}, a sketch size of $l \approx 2d$, where $d$ is the intended deflation dimension, and a power $q = 0$ or $1$, suffices to yield a quasi-optimal approximation basis, and achieve
     $$\sigma_1(\bB) \leq \|(\bI - \bPi) \bA_\mu \| \lesssim |\mu| + \sigma_d(\bA).$$
    We defer the formal discussions to~\cref{range_approx,sec:allprecond}.

    \paragraph{Implementation.}  

    The efficiency of RandRAND relies on two core operations: computing orthogonal projections and computing the action of $\bA_\mu^{-1}$ on vectors in $\mathrm{range}(\bA_\mu \bOmega)$.
The latter can be performed without explicit multiplication by $\bA_\mu^{-1}$ using the identity
\begin{equation} \label{eq:invAPi}
\bA_\mu^{-1} \bPi = \bOmega (\bA_\mu \bOmega)^\dagger,
\end{equation}
which follows directly from $\bPi = \bA_\mu \bOmega (\bA_\mu \bOmega)^\dagger$.
These operations can be incorporated in multiple ways, giving RandRAND the flexibility to adapt to specific problems and computational environments.
A natural approach is to orthogonalize the basis via QR factorization, $\bA_\mu \bOmega = \bQ \bR$, so that
$\bPi \bu = \bQ (\bQ^\mathrm{T}\bu)$ and $(\bA_\mu \bOmega)^\dagger \bu = \bR^{-\mathrm{T}} (\bQ^{\mathrm{T}}\bu)$.
However, explicit orthogonalization can be costly and may even dominate the overall runtime and memory consumption.
To mitigate this, RandRAND admits basis-less implementations that employ Q-less QR factorizations based on fast randomized transforms or sparse randomized operators.
In some settings, even factorization-free approaches are viable, where projections are computed iteratively through least-squares solves with $\bA_\mu \bOmega$ as the left-hand side. See~\cref{orthogonalization,ls_solves} for details.

   \subsubsection{Numerical illustration}
   
    As a numerical illustration, consider a linear system arising from random features regression applied to the shuttle dataset~\cite{vanschoren2014openml}. 
    We used $n = 10^4$ Fourier features associated with the rbf kernel, with the hyperparameter $\gamma = 1.33$. 
    The matrix $\bA$ is provided in factorized form, enabling efficient matvecs. 
    The regularization parameter was set to $\mu = 2.3 \cdot 10^{-13}$, resulting in $\mathrm{cond}(\bA) \approx 4 \cdot 10^{13}$. {See~\cref{sxn:basis-explicit-randrand} for more detailed description of the problem.}
    To precondition $\bA_\mu$, we constructed a R-RandRAND preconditioner using $\bOmega = \bX$, where $\bX$ is a Gaussian matrix with dimension $l = 500$. 

\begin{figure}[t] %[h!]
    \centering
    \begin{subfigure}[t]{0.235\textwidth}  
        \centering
        \includegraphics[width=\textwidth]{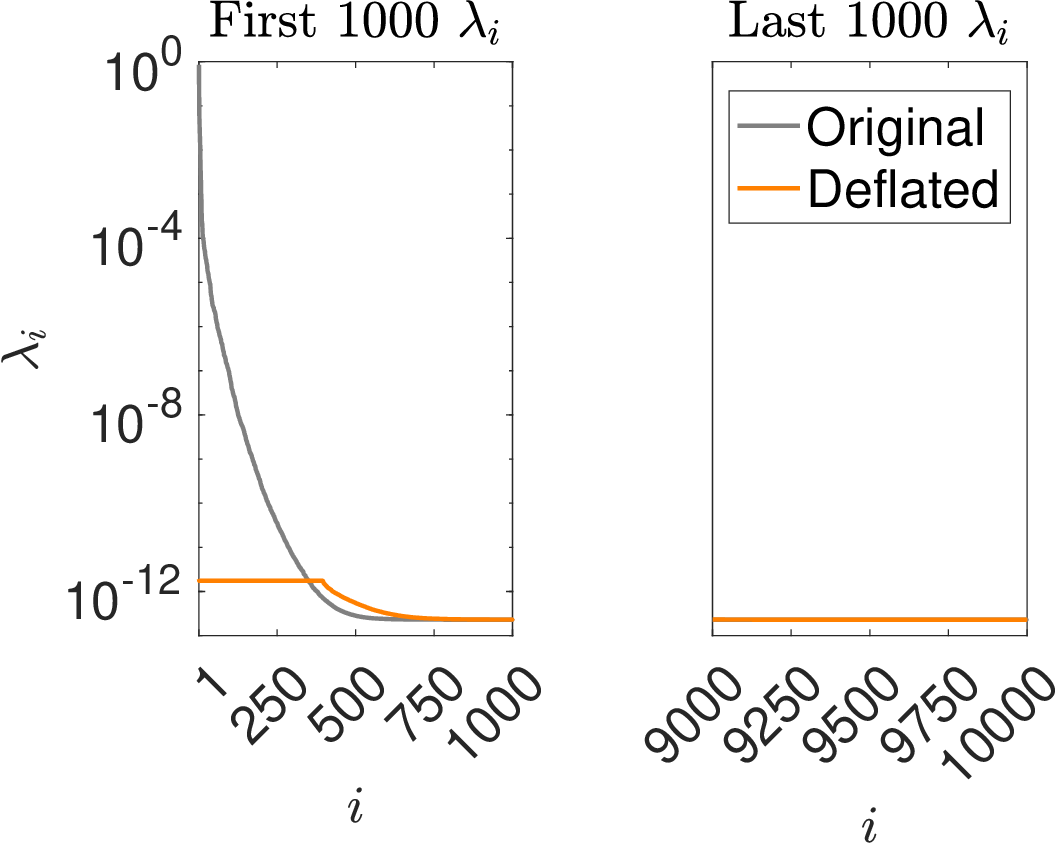} 
        \caption{\footnotesize{Eigenvalues, $\gamma = 1.33$.} }
        \label{fig:illustration1}
    \end{subfigure}
    \hspace{0.005\textwidth}
    \begin{subfigure}[t]{0.235\textwidth}
        \centering     \includegraphics[width=\textwidth]{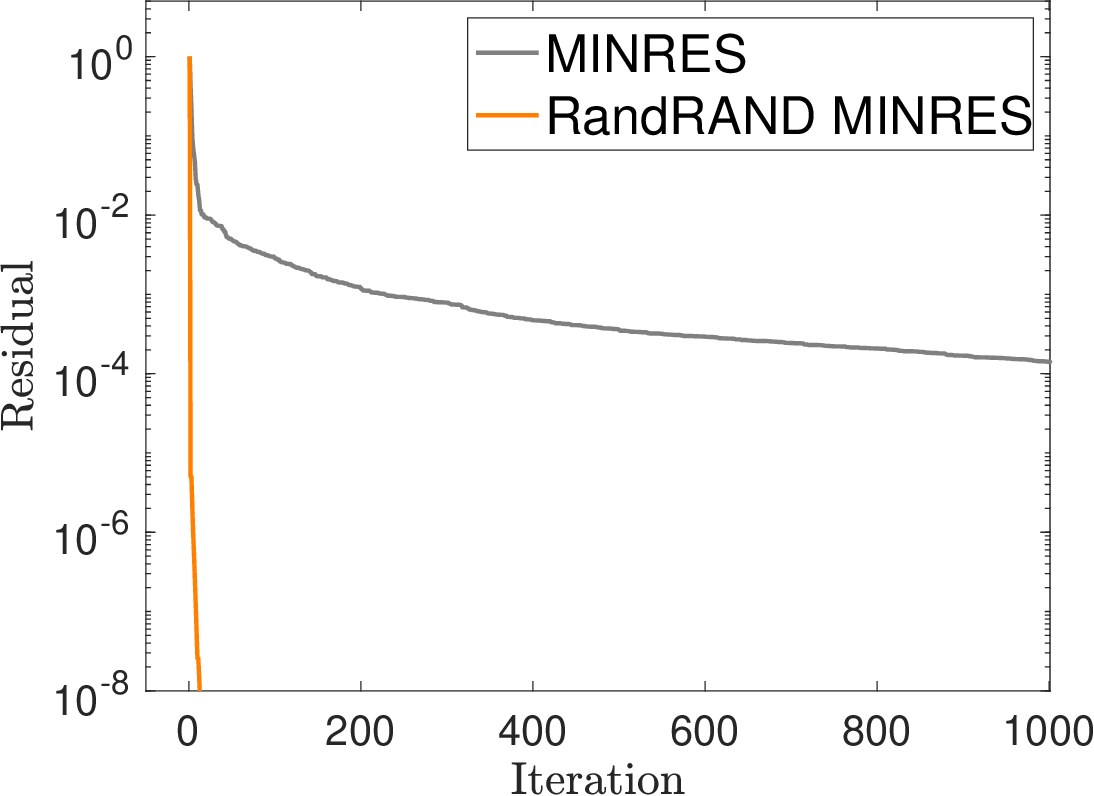} 
        \caption{MINRES error, $\gamma = 1.33$.}
        \label{fig:illustration2}
    \end{subfigure}
     \hspace{0.005\textwidth}       
    \begin{subfigure}[t]{0.235\textwidth}  
        \centering
        \includegraphics[width=\textwidth]{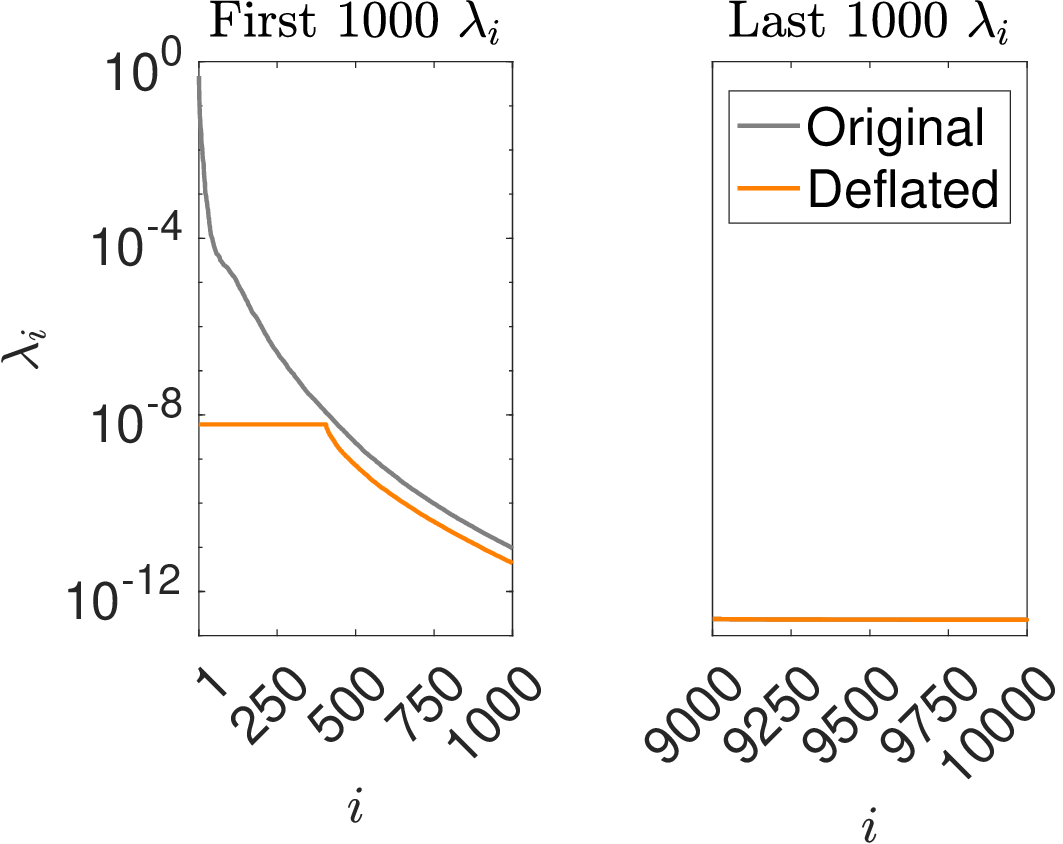} 
        \caption{Eigenvalues, $\gamma = 4.33$.}
        \label{fig:illustration3}
    \end{subfigure}
    \hspace{0.005\textwidth}
    \begin{subfigure}[t]{0.235\textwidth}
        \centering
        \includegraphics[width=\textwidth]{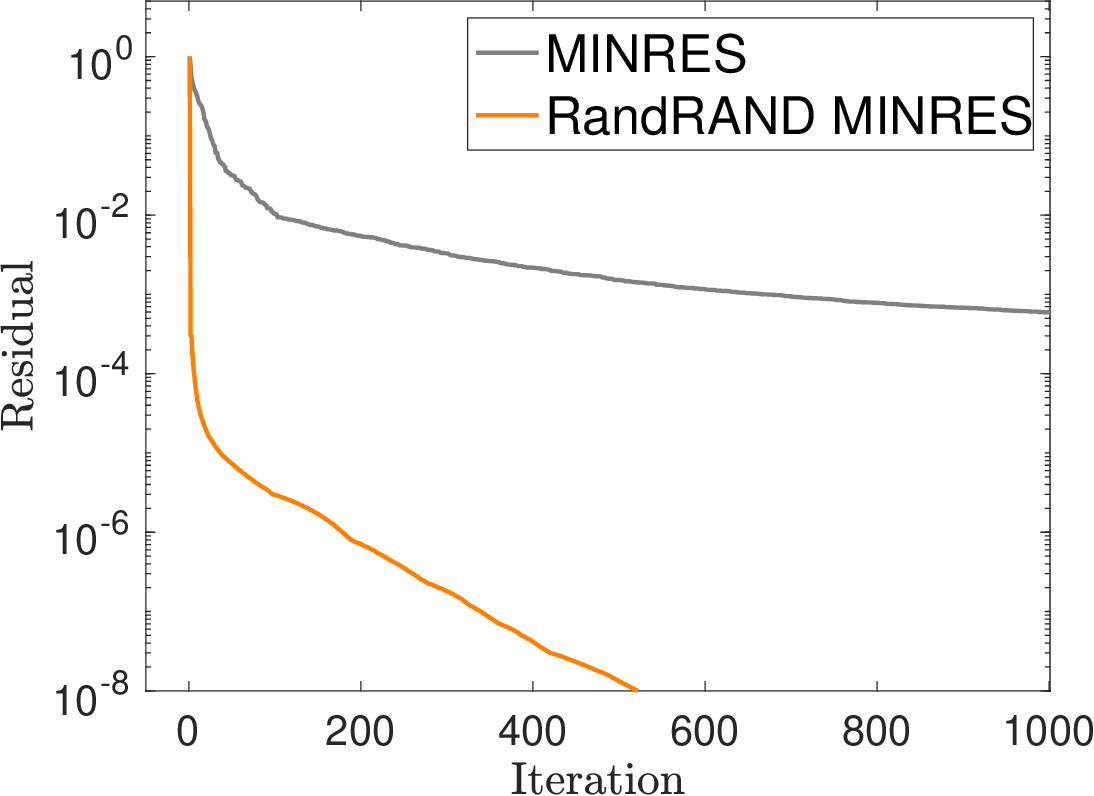} 
        \caption{MINRES error, $\gamma = 4.33$.}
        \label{fig:illustration4}
    \end{subfigure}
    \caption{\textbf{R-RandRAND preconditioning for ridge regression.} (a), (c) eigenvalues of $\bA$ and $\bA\bP$:  R-RandRAND drastically reduces the condition number; (b), (d) convergence of MINRES error for $\bA$ and $\bA\bP$:  unperconditioned MINRES exhibits dramatically slow convergence, while MINRES preconditioned with R-RandRAND exhibits fast convergence to machine precision.}
    \label{fig:main}
\end{figure}

    \Cref{fig:illustration1} compares the spectra of $\bA_\mu$ and $\bA_\mu \bP$, showing that the R-RandRAND preconditioner effectively deflates the top eigenvalues to below $2 \cdot 10^{-12}$, while preserving the smallest eigenvalue. \Cref{fig:illustration2} demonstrates the convergence of the relative residual error of MINRES solver. 
    With R-RandRAND preconditioning, MINRES achieves machine precision in fewer than $20$ iterations, whereas the convergence of unpreconditioned MINRES is extremely slow.
    Consider also a second test case with an rbf kernel hyperparameter $\gamma = 4.33$. 
    Here, $\bA$ is no longer nearly a regularized low-rank matrix, though it still has several significantly larger eigenvalues than the rest. 
    R-RandRAND effectively deflates these eigenvalues, reducing the condition number from approximately $4 \cdot 10^{13}$ to $10^4$ (see~\Cref{fig:illustration3}), and it drastically improves MINRES convergence (see~\Cref{fig:illustration4}). 
    A similar picture is also observed for other iterative methods, including CG.

    \subsection{Our Linear Solvers for SPD Systems}\label{spd_solvers}

    The potential of our approach can be realized on linear solvers for spd systems that result from integrating R-RandRAND and C-RandRAND into CG and MINRES.

    Assume that computing a matvec with $\bA_\mu$ takes $T_\mathrm{mv} \geq n$ time. 
    Let $d \geq 2$ represent the spectral dimension of $\bA_\mu$.\footnote{We have $\lambda_d \leq c(\lambda_n + \mu)$ for some constant $c$.} 
    RandRAND preconditioning allows one to compute the solution to~\cref{eq:initsys}, with a relative error $\mathrm{u}$, using 
    
    $$  \underbrace{\mathcal{O}(T_\mathrm{mv} d(q+1)~~+~~nd^2) }_{\text{compute and orthogonalize $\bA_\mu\bOmega$}} +~~ \underbrace{\mathcal{O}((T_\mathrm{mv} + nd) T_{n,d,q}  \log{\tfrac{1}{u}})}_{\text{solve preconditoned system}}~~\text{flops} $$
    and $2nd$ memory, where the coefficient  $T_{n,d,q} =  \left(\frac{n}{d} \right)^{\tfrac{1}{4q+4}} $  weakly depends on $n$ and $d$.     
    In the above complexity bound, the first term corresponds to constructing the preconditioner, and the second term corresponds to its usage in {CG or MINRES} solver. 
    This bound is independent of $\cond(\bA_\mu)$, and it can provide a drastic improvement over the cost
      \begin{equation}\label{eq:original_cost}
      T_\mathrm{mv} \cond(\bA_\mu)^{\frac{1}{2}} \log{\tfrac{1}{\mathrm{u}}}
      \end{equation}   
      of unpreconditioned solvers when $\bA_\mu$ is poorly conditioned. 
    Note that, for $q = \lfloor{\frac{1}{4}\log(\frac{n}{d})}\rfloor$, the coefficient $T_{n,d,q}$ can be replaced by a small constant.   
    Notably, our analysis shows that this can be done also for $ \lfloor{\frac{1}{4}\log(\frac{n}{d})}\rfloor > q \geq 1$, if $\lambda_d \leq  c (\frac{n}{d} )^{\frac{1}{2q}}(\lambda_n+\mu)$. 
    Moreover, in most practical settings, $T_{n,d,q}$ can be considered as a small constant already when $q=0$ or $1$.
    
   {Furthermore, RandRAND preconditioning can be implemented in a basis-less fashion without explicitly forming and operating with orthogonal basis for the range of  $\bA_\mu \bOmega$, allowing the use of fast randomized transforms and sparse randomized operators from RandNLA in efficient and numerically stable manner.}
   This approach can improve the dependence on the dimension $d$ in both the computational complexity and memory requirements. 
   As a result, the solution can be computed with high probability in\footnote{ assuming that $q \leq \lceil{\frac{1}{4}\log(\frac{n}{d})}\rceil$ and $T_\mathrm{mv} \geq n \log(n)$}:
   $$  \underbrace{{\mathcal{O}}(T_\mathrm{mv} d~~+~~d^3) }_{\text{compute and factorize sketch(es) of $\bA_\mu\bOmega$}} +~~ \underbrace{{\mathcal{O}}(T_\mathrm{mv} T_{n,d,q}  \log{\tfrac{1}{u}})}_{\text{solve preconditoned system}}\text{flops} $$   
   and $\mathcal{O}(n) \text{ memory}.$
    We here achieve a significant scaling improvement: from $nd^2$ to $d^3$ in terms of complexity, and from $nd$ to $n$ in terms of memory consumption.

    Employing preconditioning in basis-less representation can be essential when $T_\mathrm{mv} \ll nd$. 
    This can happen when $\bA_\mu$ is a sparse operator or is represented through a sparse factorization, and has only \emph{moderate} decay of dominant eigenvalues.
    In this case, preconditioning requires using a relatively large test matrix $\bOmega$ so that the explicit orthogonalization, application, and storage of the basis matrix for the range of $\bA_\mu \bOmega$ can be expensive relative to other operations and should be bypassed.  
    The basis-less representation is advantageous also when solving a sequence of linear systems, since it can allow effective recycling of $\bOmega$ along with the orthogonalization factors across iterations without the need to recompute $\bA_\mu \bOmega$. 
    For more details, see~\cref{basis_free_represent}.

   \subsection{Related Work} \label{related_work}

   \paragraph{Deflation Methods.}
    Deflation and augmentation methods have been extensively studied in the literature on Krylov methods. 
    The works in~\cite{saad2000deflated,nicolaides1987deflation,dostal1988conjugate} introduced deflation techniques using $\bA_\mu$-orthogonal projections for the CG and Lanczos methods. 
    Deflation based on $\ell_2$-orthogonal projections, which is connected to RandRAND, has been explored in studies such as~\cite{gutknecht2012spectral,coulaud2013deflation,bhattacharjee2025improved}. For a thorough historical perspective, we refer the reader to~\cite{gutknecht2012spectral,soodhalter2020survey}.

The R-RandRAND preconditioner is closely related to deflation methods for non-symmetric matrices, which decompose the solution into orthogonal components $\bPi \bx$ and $(\bI - \bPi)\bx$~\cite{coulaud2013deflation}. 
The C-RandRAND preconditioner, on the other hand, is connected to preconditioning strategies that correct the action of an operator using a preconditioner $\bP^{-1}$ of the form~\cite{coulaud2013deflation,tang2009comparison,giraud2006sensitivity}:  
\begin{equation}\label{eq:def_precond}
\bP^{-1} = (\bI - \bV \bV^\mathrm{T}) + \tau^{-1} \bV (\bV^\mathrm{T} \bB \bV)^{-1} \bV^\mathrm{T},
\end{equation}
where $\bB$ is the operator, and $\bV$ is an orthonormal deflation basis. 
Such preconditioners are designed to correct the rightmost part of the spectrum of $\bB$. 
The key insight is that the preconditioner~\cref{eq:def_precond} is algebraically equivalent to the inverse of the C-RandRAND preconditioner, when $\bB$ is taken as the inverse operator $\bA_\mu^{-1}$ and $\bV$ as an orthogonal basis for the range of $\bA_\mu \bOmega$:  
$$
(\bP^{-1})^{-1} = \left( (\bI - \bV \bV^\mathrm{T}) + \tau^{-1} \bV (\bV^\mathrm{T} \bB \bV)^{-1} \bV^\mathrm{T} \right)^{-1} = \left( \bI - \bPi + \tau^{-1} (\bPi \bA_\mu^{-1} \bPi)^{\dagger} \right)^{-1} = \bI - \bPi + \tau \bPi \bA_\mu^{-1} \bPi.
$$
This equivalence demonstrates that the effectiveness of the C-RandRAND preconditioner $\bP$ in correcting the dominant part of the spectrum of $\bA_\mu$ follows directly from the well-established effectiveness of $\bP^{-1}$ in correcting the rightmost part of the spectrum of $\bB = \bA_\mu^{-1}$.

  In classical deflation methods, the deflation subspace is typically chosen as the invariant subspace spanned by the eigenvectors associated with the extreme eigenvalues. However, computing and operating with exact eigenvectors is often computationally prohibitive. At the same time, using inexpensive approximate eigenbases can compromise the quality of the preconditioner~\cite{giraud2006sensitivity}. Moreover, the orthogonalization and maintenance of the basis can become computationally infeasible for large-scale problems.
RandRAND addresses these issues with RandNLA methods~\cite{mahoney2011randomized,halko2011finding,DM16_CACM,MD16_chapter}, providing a scalable and efficient alternative to classical approaches.

\paragraph{RandNLA Methods.}
{RandNLA applies to problems in large-scale and high-quality settings~\cite{yang2015implementing,Git16B_TR,randlapack_book_v2_arxiv}.
Existing randomized methods~\cite{el2014fast,rudi2017falkon,avron2017faster,lacotte2020effective,musco2017recursive,gower2015randomized,gonen2016solving,derezinski2024faster,derezinski2024solving,frangella2023randomized} for linear systems use sketching, sampling, and low-rank approximation. 
 However, they can be non-optimal in terms of quasi-optimality constants, numerical robustness or computational cost.  
Our work is the first to integrate projection-based spectral deflation into the RandNLA framework to address these limitations.}

Notice that the Nystr{\"o}m approximation $\hat{\bA}_\mu$ of $\bA_\mu$ associated with  test matrix $\bOmega$ can be represented in the following~form:
\begin{equation} \label{Nystrom_vs_CRandRAND}
\hat{\bA}_\mu := (\bA_\mu \bOmega) (\bOmega^\mathrm{T} \bA_\mu \bOmega)^{-1} (\bOmega^\mathrm{T} \bA_\mu) = \left( \left ( (\bA_\mu \bOmega)^\mathrm{T}\right)^\dagger (\bOmega^\mathrm{T} \bA_\mu \bOmega) (\bA_\mu \bOmega)^\dagger \right)^\dagger = \left( \bPi \bA_\mu^{-1} \bPi \right)^\dagger.
\end{equation} 
This implies that the C-RandRAND preconditioner $(\bI -\bPi) + \tau \bPi \bA_\mu^{-1} \bPi$ can be interpreted as $(\bI -\bPi) + \tau \hat{\bA}_\mu^\dagger$, which  establishes the link with the methods in~\cite{frangella2023randomized,derezinski2024faster,gonen2016solving} based on low-rank approximation. 
Among these methods, the preconditioner proposed in~\cite{frangella2023randomized} appears to be the closest to C-RandRAND. 
However, it differs in key aspects:  it uses the Nystr{\"o}m approximation $\hat{\bA}$ of the unregularized operator $\bA$, instead of $\bA_\mu$; it relies on the approximation error $\|\bA - \hat{\bA}\|$ and not the projection error;  it selects the coefficient $\tau$ based on the smallest eigenvalue of $\hat{\bA}$; and it requires computing  $\hat{\bA}^\dagger$ in SVD form, which can be computationally prohibitive for large-scale problems.  
The condition number bound from~\cite[Proposition 5.3]{frangella2023randomized} for Nystr{\"o}m approximation-based preconditioning differs from the bound~\eqref{eq:cndbound0} for RandRAND by a factor of $\frac{\mu + \|\bA - \hat{\bA}\|}{\|(\bI - \bPi) \bA_\mu\|}$. 
This factor can be large when the basis $\bA\bOmega$ approximates well $\bA$, but the Nystr{\"o}m approximation $\hat{\bA}$ is not a good approximation of $\bA$. 
This can happen for instance when $\bA$ has a heavy-tailed spectrum~\cite{gittens2013revisiting,gittens2016revisiting,Git16B_TR}. 

If we consider $\bOmega$ be a Gaussian OSE of dimension $l = \mathcal{O}(d)$, where $d$ is such that $\lambda_d$ is comparable to $\lambda_n$, then, according to~\cite{gittens2013revisiting}, we can expect $\|\bA - \hat{\bA}\| = \mathcal{O}\left(1 + \left(d^{-1}{\mathrm{sr}_d(\bA^{2q+1})}\right)^\frac{1}{2q+1}\right) \lambda_d$. 
At the same time, according to our analysis in~\cref{range_approx}, we have  $\|(\bI - \bPi) \bA_\mu\| = \mathcal{O}\left(1 + \left(d^{-1}{\mathrm{sr}_d(\bA^{2q} \bA_\mu^2)}\right)^\frac{1}{2q+2}\right) \lambda_d$. 
The stable rank can be as large as $\mathrm{sr}_d(\bA^{2q+1}) = \mathrm{sr}_d(\bA^{2q}\bA_\mu^2) = \mathcal{O}(n)$, even when the spectral rank $d$ of $\bA$ is small. 
In these cases, RandRAND achieves a bound lower by a factor of  
$$
\mathcal{O}\left(\left(\frac{n}{d}\right)^{\tfrac{2q+2}{2q+1}}\right),
$$  
which can be a drastic improvement. 
\Cref{validation_randrand_properties} provides a numerical validation of this fact, demonstrating that the error bound for Nystr{\"o}m approximation-based methods can be orders of magnitude larger than that of RandRAND.

Besides improving quasi-optimality constants, RandRAND allows a basis-less implementation, without the need to explicitly form and operate with the orthogonal basis for the range of $\bA_\mu \bOmega$. 
This approach can effectively exploit fast randomized transforms and sparse randomized operators from RandNLA to drastically reduce the complexity and memory consumption in scenarios where computing a low-rank approximation of $\bA$ in SVD form, as done in e.g. \cite{gonen2016solving,frangella2023randomized}, is computationally expensive relative to other computations. 

Several basis-less randomized approaches for solving linear systems have been proposed~\cite{derezinski2024solving,derezinski2024faster,derezinski2025approaching,derezinski2025fine,derezinski2025randomized}. 
The methods in~\cite{derezinski2024solving,derezinski2025fine,derezinski2025randomized} offer attractive advantages in terms of flop counts and data passes, but they are of Kaczmarz type and thus not directly comparable to RandRAND. The method in~\cite{derezinski2024faster} bears some resemblance to RandRAND, but it relies on the Nystr{\"o}m approximation and consequently yields non-optimal performance guarantees. A recent extension~\cite{derezinski2025approaching} achieves asymptotically optimal flop complexity. However, it proceeds with  nested layers of preconditioners, which in practice can incur large hidden constants in the complexity bounds and substantially hinder numerical robustness and efficiency. 
In contrast, the representations with orthogonal projectors enable RandRAND to ensure efficiency, numerical stability, and ease of implementation on high-performance architectures. Furthermore, RandRAND may be employed through matvecs, without access to matrix entries, making it applicable to a broader range of large-scale problems.

Another important distinction from existing randomized preconditioning techniques is that they are typically limited to spd operators $\bA_\mu$, whereas RandRAND naturally generalizes to indefinite and non-symmetric operators.

    \section{Randomized Range Approximation} 
\label{range_approx}

    {In this section, we present rigorous accuracy guarantees for the range of $\bA_\mu \bOmega$ in approximating the action of $\bA_\mu$, a key ingredient in RandRAND. We consider $\bOmega$ constructed as in~\cref{eq:testmatrix} with several standard choices of random matrix $\bX$. Establishing such guarantees requires extending the results of~\cite{halko2011finding,boutsidis2013improved,tropp2011improved,cohen2015optimal,gittens2016revisiting} to  the shifted subspace iteration.}

    \paragraph{Gaussian matrices.} 
    We start with Gaussian matrices.
    A matrix $\bX \in \mathbb{R}^{l \times n}$ is called Gaussian if its entries are iid normal random variables with mean $0$ and variance $l^{-1}$.  Gaussian matrices are advantageous because their theoretical analysis is both simple and tight. 
    Beyond characterization of the accuracy of algorithms involving Gaussian matrices, such analysis can also (sometimes~\cite{randnla_kdd24_TR}) help to understand the practical performance of algorithms that rely on other sketching operators. 
    This includes fast randomized transforms and multi-level embeddings for which obtaining precise theoretical guarantees can be difficult but which in practice offer comparable accuracy for a given sketching dimension.
    
    {The following lemma, adapted from~\cite[Theorem 8]{halko2011finding}, bounds $\|(\bI - \bPi) \bA_\mu\|$ when using Gaussian test matrix $\bOmega=\bX^\mathrm{T}$ without subspace iteration.}
    
    \begin{lemma}\label{thm:apriori_gen00}
    Draw a Gaussian matrix $\bX \in \mathbb{R}^{l \times n}$ using sketching dimension $l =2k+2\geq8$. Define test matrix $\bOmega = \bX^\mathrm{T}$.  Then the projection of $\bA_\mu$ onto the range of $\bA_\mu \bOmega$ satisfies the following error bound with probability $\geq 1-6^{-k}$:
    \begin{equation} \label{eq:gen_bound0}
        \|(\bI - \bPi) \bA_\mu\|  \leq   \left ( C_1 +   C_2  \sqrt{ \frac{\mathrm{sr}_{k}(\bA_\mu^\mathrm{T}\bA_\mu)}{k} } \right )  \sigma_k(\bA_\mu), 
    \end{equation} 
    where $C_1 = 18$ and $C_2 = 8\sqrt{2}$. The stable rank $\mathrm{sr}_{k}(\bA_\mu^\mathrm{T}\bA_\mu):=    \sum_{j\geq k} \frac{\sigma^2_j(\bA_\mu)}{\sigma^2_k(\bA_\mu)}   \leq n$  describes the decay of the trailing singular values  of  $\bA_\mu$. Furthermore, in expectation, $\|(\bI - \bPi) \bA_\mu\|$ satisfies~\cref{eq:gen_bound0}
    with smaller constants:  $C_1 = 2$ and $C_2 = \sqrt{2}\mathrm{e}$.
    \begin{proof}
        See~\cref{appendix_a}.
    \end{proof}
    \end{lemma}
    According to~\cref{thm:apriori_gen00}, we can guarantee with high probability the quasi-optimality of approximation basis $\bA_\mu\bX^\mathrm{T}$. 
    In particular, we have the following bound for $l =2k+2\geq6$:  
     $$ \frac{\mathbb{E}(\|(\bI-\bPi)\bA_\mu\|)}{\sigma_n(\bA_\mu)}   \leq 6 \left( \frac{n}{k} \right)^{\frac{1}{2}} \frac{\sigma_k(\bA_\mu)}{\sigma_n(\bA_\mu)} =  6 \left( \frac{n}{k} \right)^{\frac{1}{2}} \min_{\mathrm{rank}(\bPi^*) = k-1} \frac{\|(\bI-\bPi^*)\bA_\mu\|}{\sigma_n(\bA_\mu)} ,  $$
     where the quasi-optimality constant $6 \left( \frac{n}{k} \right)^{\frac{1}{2}}$ is independent of the properties of $\bA$ and $\mu$.
     
     Let us now present an extension of~\cref{thm:apriori_gen00} to the shifted subspace iteration. We restrict the analysis to symmetric systems  using the approximation basis $\bOmega = \bA_\mu \bA^q \bX^\mathrm{T}$. For non-symmetric cases similar guarantees should hold for the approximation basis $\bOmega = \bA_\mu (\bA^\mathrm{T} \bA)^q \bX^\mathrm{T}$. 
     
    \begin{theorem}\label{thm:apriori_gen2}
    Let $\bA$ be a symmetric matrix with singular values $\sigma_1 \geq \sigma_2 \geq \hdots \geq  \sigma_n$.
    Draw a Gaussian matrix $\bX \in \mathbb{R}^{l \times n}$ using sketching dimension $l =2k+2\geq8$. Define test matrix $\bOmega = \bA_\mu \bA^q \bX^\mathrm{T}$ with $q \geq 1$.   Then the projection of $\bA_\mu$ onto the range of $\bA_\mu \bOmega$ satisfies the following error bound with probability $\geq 1-6^{-k}$:
    \begin{equation} \label{eq:gen_bound}
        \|(\bI - \bPi) \bA_\mu\|  \leq 3 |\mu| +   \left ( 2C_1 +   2C_2  \sqrt{ \frac{\mathrm{sr}_{k}(\bA_\mu^2\bA^{2q})}{k} } \right )^{\frac{1}{q+1}}  \sigma_k^{\frac{q}{q+1}} (\sigma_k+ |\mu|)^{\frac{1}{q+1}}, 
    \end{equation} 
    where $C_1 = 18$ and $C_2 = 8\sqrt{2}$. The stable rank $\mathrm{sr}_{k}(\bA_\mu^2\bA^{2q}):=    \sum_{j\geq k} \frac{\sigma_j(\bA_\mu^2\bA^{2q})} {\sigma_k(\bA_\mu^2\bA^{2q})} \leq n$  describes the decay of the trailing singular values of  $\bA_\mu^2\bA^{2q}$. Furthermore, in expectation, $\|(\bI - \bPi) \bA_\mu\|$ satisfies~\cref{eq:gen_bound}
    with smaller constants:  $C_1 = 2$ and $C_2 = \sqrt{2}\mathrm{e}$.
    \begin{proof}
        See~\cref{appendix_a}. 
    \end{proof}
    \end{theorem}

    \noindent 
    The constants in~\cref{thm:apriori_gen2} can be refined if $\bA$ is spd and $\mu \geq 0$, as shown in the following result.
    
    \begin{proposition}\label{thm:apriori2}
    Let $\bA$ be a spd matrix and $\mu \geq 0 $. Let the sketching matrix be the same as in~\cref{thm:apriori_gen2}. Then with probability at least $\geq 1-6^{-k}$,
    \begin{equation} \label{eq:spd_bound}
        \|(\bI - \bPi) \bA_\mu\|  \leq  \mu +   \left ( C_1 +   C_2  \sqrt{ \frac{\mathrm{sr}_{k}(\bA_\mu^2\bA^{2q})}{k} } \right )^{\frac{1}{q+1}}  \lambda_k^{\frac{q}{q+1}} (\lambda_k+ \mu)^{\frac{1}{q+1}}, 
    \end{equation}    
    holds with constants $C_1 = 18$ and $C_2 = 8\sqrt{2}$.  Furthermore, in expectation, $\|(\bI - \bPi) \bA_\mu\|$ satisfies~\cref{eq:gen_bound}
    with constants $C_1 = 2$ and $C_2 = \sqrt{2}\mathrm{e}$.
    \begin{proof}
        See~\cref{appendix_a}.
    \end{proof}
    \end{proposition}

    \Cref{thm:apriori_gen2,thm:apriori2} imply that subspace iteration reduces the scaling of the quasi-optimality constant with dimension $n$.   
    For instance, for spd systems and $\mu \geq 0$, we have established the following guarantee: 
    $$ \frac{\mathbb{E}(\|(\bI-\bPi)\bA_\mu\|)}{\lambda_n+\mu}   \leq  1+ \left( \frac{n}{k} \right)^{\frac{1}{2q+2}} \left (6 \left( \frac{\lambda_k}{\lambda_n+\mu} \right)^{q} \left( \frac{\lambda_k+\mu}{\lambda_n+\mu} \right) \right)^{\frac{1}{q+1}}.   $$
    It follows that if  $d$ is larger than or equal to the count of eigenvalues of $\bA$ larger than $\mu$ {, i.e.,
    ${d \geq \#\{ i \mid \lambda_i \geq \mu \}},$}
    then the sketching dimension $l \geq 2d+4$ ensures that 
    $$ \frac{\mathbb{E}(\|(\bI-\bPi)\bA_\mu\|)}{\lambda_n+\mu}   \leq  1+ \left( 3.5 \, \frac{n}{d} \right)^{\frac{1}{2q+2}}.    $$
    Similarly, if $d$ is larger or equal to the count of eigenvalues of $\bA$ larger than $0.1\mu$, 
    then the sketching dimension $l \geq 2d+4$ ensures that 
     $$ \frac{\mathbb{E}(\|(\bI-\bPi)\bA_\mu\|)}{\lambda_n+\mu}   \leq  1+ \left( 2.6 \cdot 10^{-2q} \frac{n}{d} \right)^{\frac{1}{2q+2}}.    $$

    \paragraph{Fast randomized transforms.}
   For $n$ being a power of $2$, the Subsampled Randomized Hadamard Transform (SRHT) $\bX \in \mathbb{R}^{l \times n}$ is defined as a product of a diagonal matrix of iid Rademacher variables $\pm 1$, a Walsh-Hadamard matrix, and uniform subsampling matrix scaled by $l^{-\frac{1}{2}}$. 
   If $n$ is not a power of $2$,  SRHT can be employed via zero-padding the input. This is equivalent to taking $\bX$ as the first $n$ columns of an SRHT matrix of size $l \times s$, where $s$ is the smallest power of $2$ such that $n \leq s < 2n$. 
   In this paper, such an $\bX$ will be referred to also as SRHT.

   SRHT matrices can be computationally advantageous from the complexity standpoint.  
   They can be applied to vectors with only $s \log_2 s$ flops, when using the fast Hadamard transform, and with $2s \log_2 l$ flops when using the procedure form~\cite{ailon2009fast}. 
   This is in stark contrast to the application cost of $2 n l$ for unstructured matrices. 

    In terms of accuracy, SRHTs offer guarantees for the range approximation that are asymptotically comparable to those provided by Gaussian matrices.    
    In particular, \cref{thm:apriori3} suggests that to achieve the bound~\cref{eq:gen_bound} for $\|(\bI - \bPi) \bA_\mu\|$, SRHT matrices require a sketching dimension $l$ that is only logarithmically larger than that for Gaussian matrices. 
    Moreover, our empirical results indicate that, in practice, using a larger sketching size for SRHT matrices is often not needed for achieving accuracy comparable to that of Gaussian matrices.

    \begin{theorem}[Deviation bound for SRHT sketching based on Theorem 2.1 from~\cite{boutsidis2013improved} and \cite{tropp2011improved}]\label{thm:apriori3} 
    Let $s$ be the smallest power of $2$ such that $n \leq s < 2n$, and let $0<\delta<1$ be a coefficient describing failure probability.   Draw an SRHT matrix $\bX \in \mathbb{R}^{l \times n}$, with  $l \geq 3\log(s/\delta)^2$ rows. Let $k \geq 3$ be the maximum integer such that $$3 k \log(s/\delta) \leq l~\text{~and~}
 18 \left ( \sqrt{k} + \sqrt{8 \log(s/\delta)} \right )^2 \log(k/\delta)  < l. $$
    Then with probability at least $1-\delta$ the bound~\cref{eq:spd_bound} holds if $\bA$ is spd, the bound~\cref{eq:gen_bound} holds if $\bA$ is symmetric, and the bound~\cref{eq:gen_bound0} holds if $\bA$ is a general matrix,  with constants $C_1 = 5$ and $C_2 = 1$. 
    \begin{proof}
        See~\cref{appendix_a}.
    \end{proof}
    \end{theorem}

    \paragraph{Sparse embeddings.}
    In this paper, matrix $\bX \in \mathbb{R}^{l \times n}$ is called a sparse random operator with density $\gamma/l$ if: (i) its columns are generated independently, (ii) in each column, $\gamma$ entries are uniformly selected and assigned values $\frac{1}{\sqrt{\gamma}}$ or $-\frac{1}{\sqrt{\gamma}}$ with equal probabilities, and (iii) all non-selected entries are set to zero. 
    This choice of $\bX$ offers great computational advantages similar to those of SRHT, as the application cost of a sparse random operator is just $n \gamma$ flops.

    Sparse operators satisfy a similar guarantee for the range approximation as Gaussian and SRHT matrices. Specifically,~\cref{thm:apriori4} shows that a sparse embedding with logarithmic density requires a sketching dimension only logarithmically larger than that of a Gaussian matrix to guarantee comparable accuracy for the range approximation. Moreover, recent works~\cite{chenakkod2024optimal1,chenakkod2024optimal,chenakkod2025optimal} have developed improved sparse OSEs with almost optimal bounds on the required sketching dimension:~$l = \mathcal{O}(k)$.

    \begin{theorem}[Deviation bound for sparse operator based on Theorem 8 from \cite{cohen2015optimal}]\label{thm:apriori4} 
    Let $0<\delta<1$ be a parameter describing failure probability. Let $\bX$ be a sparse random operator with $l = \mathcal{O}(k \log(k/\delta))$ rows and $\gamma = \mathcal{O}( \log(k/\delta))$ non-zero entries per column. With probability at least $1-\delta$ we have that the bound~\cref{eq:spd_bound} holds if $\bA$ is spd, the bound~\cref{eq:spd_bound} holds if $\bA$ is symmetric, and the bound~\cref{eq:gen_bound0} holds if $\bA$ is a general matrix,  with constants $C_1 = 8$ and $C_2 = 2$. 
    \begin{proof}
        See~\cref{appendix_a}.
    \end{proof}
    \end{theorem}
       
   \paragraph{Multi-level embeddings.} 
   OSEs are composable and allow for the combining of the low application cost of fast and sparse transforms with the optimal dimension of Gaussian matrices.  
   Multi-sketching has been considered in several studies~\cite{avron2017faster,balabanov2019randomized,higgins2023analysis}.  
   In particular, we can construct \(\bX\) as  
    \begin{equation} \label{eq:multilevelose}
        \bX = \bX_2 \bX_1,
    \end{equation} 
    where \(\bX_1 \in \mathbb{R}^{l_1 \times n}\) is an SRHT or sparse random operator with \(l_1 = \mathcal{O}(k \log(k/\delta))\) rows and \(\gamma = \mathcal{O}(\log(k/\delta))\) nonzero entries per column, and \(\bX_2 \in \mathbb{R}^{l \times l_1}\) is a Gaussian matrix with small number (say $l =2k$) of rows. It can be shown that in this case we have a guarantee for the range approximation similar to~\cref{thm:apriori_gen2,thm:apriori2}.
    The multi-level OSEs of the form~\cref{eq:multilevelose} can be particularly useful for the basis-less RandRAND preconditioning described in~\cref{basis_free_represent}, where both the cost of applying \(\bOmega\) and its size significantly impact the overall computational costs.

    \paragraph{Data-aware embeddings.}  \label{sprse_emb}
   
        {Beyond  OSE operators $\bX$, RandRAND is also compatible with data-aware embeddings, such as those constructed via leverage score sampling~\cite{alaoui2015fast,gittens2016revisiting,woodruff2014sketching}, the sampling based on randomly pivoted Cholesky procedure~\cite{chen2025randomly,diaz2023robust}, and Krylov subspaces associated with a random starting matrix $\bTheta$~\cite{musco2015randomized,tropp2023randomized}. While a detailed analysis of data-aware embeddings is beyond the scope of this manuscript, we briefly highlight key aspects below.}

        {For symmetric matrices, the block Krylov basis is defined as $\mathrm{range}(\bX^\mathrm{T})  = \mathcal{K}_r(\bA,\bTheta) := \mathrm{range}([\bTheta, \bA\bTheta, \hdots, \bA^r\bTheta])$ and can be shown to capture well the action of $\bA$~\cite{musco2015randomized,tropp2023randomized}. We have $\mathrm{range}(\bA_\mu \bA^{q} \bX^\mathrm{T}) \subseteq  \bA_\mu \mathcal{K}_{r+q}(\bA,\bTheta)$ so the basis for  $\mathcal{K}_{r+q+1}(\bA,\bTheta)$ can be used instead of $\bOmega = \bA^{q} \bX^\mathrm{T}$, with preserving (or improving) condition number bounds. Similar consideration also holds for non-symmetric matrices.} 
         
        {Sampling is particularly attractive in kernel methods, as it avoids explicit matrix-matrix products with the full operator $\bA$. In this context, one can build on the results of~\cite{alaoui2015fast,gittens2016revisiting,chen2025randomly} to control the quality of the sampled columns as range approximation basis. 
         Moreover, in many problems, uniform column sampling is sufficient. We adopt this approach in our kernel regression experiments in~\cref{sec:basis-less-pde-appn}.}

        {Another promising approach left for future research is to combine RandRAND with sparse embeddings that integrate OSEs with leverage score sampling such as the LEverage Score Sparsified (LESS) embeddings from~\cite{derezinski2021sparse}.}

    \section{RandRAND preconditioners} 
\label{sec:allprecond}

In this section, we describe our RandRAND preconditioners.
Assume that we are given a test matrix $\bOmega$ such that $\bA_\mu\bOmega$ captures the action of $\bA_\mu$ sufficiently well. 
This matrix can be obtained with the randomized range approximation from~\cref{range_approx}.

\subsection{R-RandRAND Preconditioner} 
\label{sec:precond}     
    Consider the R-RandRAND preconditioner for spd $\bA_\mu$ that reduces the action space of $\bA_\mu$:
    \begin{equation}\label{eq:randrandprecond}
        \bP = \bA^{-1}_\mu\left ( (\bI -\bPi)\bA_\mu (\bI -\bPi) + \tau \bPi \right ), 
    \end{equation}
    and satisfies $\bA_\mu \bP \bPi = \tau \bPi$.

    In this case, we can compute the solution to the system $\bA_\mu \bx = \bb$ in two stages. 
    In the first stage, we solve the deflated system
    \begin{equation}\label{eq:APdef2}
        \bB \by = \bb,
     \end{equation}   
     where $\bB = \bA_\mu \bP = (\bI -\bPi)\bA_\mu (\bI -\bPi) + \tau \bPi$, with an iterative method for spd systems such as CG or MINRES. 
     Then, in the second stage, we retrieve the solution to the original system via 
    \begin{equation}\label{eq:APdef3}
        \bx = \bP \by = (\bA_\mu^{-1} \bPi) ( \tau \by - \bA_\mu (\bI - \bPi)\by ) + (\bI - \bPi)\by = \bOmega (\bA_\mu \bOmega)^\dagger ( \tau \by - \bA_\mu (\bI - \bPi)\by ) + (\bI - \bPi)\by,
     \end{equation}   
    where we used the fact that $\bA_\mu^{-1} \bPi = \bA_\mu^{-1} \bA_\mu \bOmega (\bA_\mu \bOmega)^\dagger = \bOmega (\bA_\mu \bOmega)^\dagger$. The matvecs with $(\bA_\mu \bOmega)^\dagger$ and $\bPi$ can be computed through orthogonalization of $\bA_\mu \bOmega$ or by using (randomized) iterative methods discussed in details in~\cref{ls_solves}.

    The following \Cref{thm:cnd_bound} shows that R-RandRAND preconditioning effectively deflates the range of $\bA_\mu \bOmega$ from the range of $\bA_\mu$. 
   
    \begin{proposition} \label{thm:cnd_bound}    
    Let $\bA_\mu$ be spd. Let $\bE = (\bI - \bPi) \bA_\mu (\bI - \bPi)$. Let $\rho\geq 1$ be some user-specified parameter. If 
    $\rho^{-1} \lambda_{\min}(\bA_\mu) \leq \tau \leq \rho \| \bE\|,$
    then the preconditioned operator $\bA_\mu \bP$ in~\cref{eq:APdef2} satisfies 
        \begin{equation} \label{eq:cnd_bound}
         \mathrm{cond}(\bA_\mu \bP) \leq  \rho \frac{\|\bE\|}{\lambda_{\mathrm{min}}(\bA_\mu)}.
        \end{equation}
        %Furthermore, it holds that  
        %\begin{equation}\label{eq:cnd_boundE}
        %    \lambda_{n-l+1}(\bE) \geq \lambda_{\min}(\bA_\mu)
        %\end{equation}.
        \begin{proof}
            See~\cref{appendix_b}.
        \end{proof}
    \end{proposition}

    \noindent 
    In practice, the coefficient $\tau$ satisfying $\lambda_{\min}(\bA_\mu) \leq \tau \leq  \|\bE\|$ in~\cref{thm:cnd_bound} can be obtained with a few power iterations on matrix $\bE$, and requires negligible computational cost. 

    The condition number bound~\cref{eq:cnd_bound} can be directly integrated with the bounds for the randomized range approximation from~\cref{range_approx}.   

    \begin{proposition} \label{thm:cnd_bound_spec}
         Let $\bA_\mu$ be spd. Draw a Gaussian matrix $\bX \in \mathbb{R}^{l \times n}$ using sketching dimension $l = 2ck+4$ for some constant $1 <  c \leq 2$ and an integer $k \geq 1$. Let $\bA_\mu \bP = \bE + \tau \bPi$ be the R-RandRAND preconditioned operator~\cref{eq:APdef2} where $\bPi$ is orthogonal projector onto a basis $\bA_\mu \bA^q \bX^\mathrm{T}$ and $\bE = (\bI - \bPi) \bA_\mu (\bI - \bPi)$. Let $\tau$ satisfy $\lambda_{\min}(\bA_\mu) \leq \tau \leq  \|\bE\|$. Then we have with probability at least $1-6^{-k}$ for $q=0$:
        \begin{equation}\label{eq:condbound1}
        \mathrm{cond}(\bA_\mu \bP) \leq  1 +  \left (C_3 \frac{n}{(c-1)k} \right )^{\frac{1}{2}} \cond_{k}(\bA_\mu^2)^{\frac{1}{2}},
        \end{equation}
        and for $q \geq 1$:
        \begin{equation}\label{eq:condbound2}
         \mathrm{cond}(\bA_\mu \bP) \leq  1 +   \left (C_3 \frac{n}{(c-1)k} \right )^{\frac{1}{2q+2}} \cond_{k, \mu^{2q+2}}(\bA_\mu^2 \bA^{2q})^{\frac{1}{2q+2}},
        \end{equation}
        where $C_3 = 2 \cdot 18^2$. Furthermore, in expectation,~\cref{eq:condbound1,eq:condbound2} hold with smaller constant: $C_3 = 4 \mathrm{e}^2$.
        \begin{proof}
            See~\cref{appendix_b}.
        \end{proof}
    \end{proposition}
    
    Note that the $k$-stable condition numbers $\cond_{k}(\bA_\mu^2)$ and $\cond_{k, \mu^{2q+2}}(\bA_\mu^2 \bA^{2q})$ from~\cref{thm:cnd_bound_spec} satisfy 
    $$ \cond_{k}(\bA_\mu^2)^{\frac{1}{2}} \leq  \frac{\lambda_k+\mu}{\lambda_n+\mu} ~ \text{ and } \cond_{k, \mu^{2q+2}}(\bA_\mu^2 \bA^{2q})^{\frac{1}{2q+2}} \leq  2 \frac{\lambda_k+\mu}{\lambda_n+\mu}.$$
    It follows that R-RandRAND preconditioning with $l = \mathcal{O}(k)$ achieves a near-optimal deflation of the first $k$ eigenvalues  already when using the power $q=0$ or $q=1$.   
    Furthermore, for $q=0$, the $k$-stable condition number behaves similarly to the standard averaged condition number, and it is always greater than or equal to $1$; and for $q\geq 1$, the $k$-stable condition number can be significantly smaller than $1$ and result in a constant bound on $\mathrm{cond}(\bB)$, when the trailing eigenvalues decay sufficiently rapidly.  
   
   Similar results also hold for SRHT and sparse matrices. In particular, it follows from~\cref{thm:apriori3,thm:apriori4} that the bounds~\cref{eq:condbound1,eq:condbound2} with coefficients $C_3=32$ and $c=2$ hold with probability $1-\delta$, when $\bX$ is an SRHT matrix with $l = \mathcal{O}(k \log(n/\delta))$ rows, or a sparse random matrix with  $l = \mathcal{O}(k \log(k/\delta))$ rows and $\gamma = \mathcal{O}( \log(k/\delta))$ nonzeros per~column.

\subsubsection{Left R-RandRAND Preconditioner}

Assuming that $\bA_\mu$ is spd, we can apply the transpose of $\bP$ in~\cref{eq:randrandprecond} as a left preconditioner. In particular, we can compute the solution to $\bA_\mu \bx = \bb$ via solving 
\begin{equation}\label{eq:APdef2_left}
   \bB\bx= \bP^\mathrm{T}\bb,
\end{equation}
where the operator $\bB = (\bI-\bPi)\bA_\mu(\bI-\bPi) + \tau\bPi$ is as in \eqref{eq:APdef2}, but the right-hand side is updated to
$$
   \bP^\mathrm{T}\bb= ((\bI-\bPi)\bA_\mu(\bI-\bPi)\ + \tau\bPi)\bA_\mu^{-1} \bb.
$$
As before, for efficiency, we must use the fact that $\bPi\bA_\mu^{-1}\bb$ can be computed by $\bigl((\bA_\mu\bOmega)^\dagger\bigr)^\mathrm{T}\bOmega^\mathrm{T}\bb$.

A downside of left preconditioning compared to right preconditioning is that the residual of \eqref{eq:APdef2_left} no longer corresponds directly to the residual of the original system $\bA_\mu \bx = \bb$.

\subsubsection{Split R-RandRAND Preconditioner} \label{sec:precond_split}

The R-RandRAND preconditioning can be used not only in its classical role, as a preconditioner for $\bA_{\mu}$, but also as a preconditioning method that splits the right hand side into orthogonal parts.  
Specifically, we can first compute a solution to the singular reduced system 
\begin{equation}\label{eq:APdef2_reduced}
\bE \by = (\bI -\bPi) \bb,
\end{equation}
where $\bE = (\bI - \bPi)\bA_\mu(\bI - \bPi)$, and then retrieve $\bx$ via the relation  
$$
\bx =  (\bA_\mu^{-1} \bPi) (\bb - \bA_{\mu}(\bI-\bPi) \by) + (\bI-\bPi) \by = \bOmega (\bA_\mu \bOmega)^\dagger(\bb - \bA_{\mu}(\bI-\bPi) \by) + (\bI-\bPi) \by. 
$$
The advantage of this approach is the reduction of the number of times the projector $\bPi$ needs to be applied inside Krylov iterative solver,\footnote{This is since when $\bu = (\bI-\bPi)\bu$, the matvec $\bu \gets \bE \bu$ can be computed as $(\bI -\bPi) \bA_\mu\bu$ and requires only one application of $\bPi$.} which can be important in settings where matvecs with $\bPi$ are expensive compared to other operations.

For spd systems $\bA_\mu$, one can show that the effective condition number $\mathrm{cond}_{\mathrm{eff}}(\bE)$, defined as the largest eigenvalue to the smallest non-zero eigenvalue, satisfies for $\tau > 0$:
    $$
        \mathrm{cond}_{\mathrm{eff}}(\bE) := \frac{\lambda_{1}{(\bE)}}{\lambda_{n-l+1}(\bE)} \leq \mathrm{cond}(\bE + \tau \bPi),
    $$
which implies the same guarantees of convergence of Krylov solvers as those for the right and left R-RandRAND.

\subsection{C-RandRAND Preconditioner} \label{sec:c_precond}
As an alternative to reducing the action space of $\bA_\mu$, as in R-RandRAND, we can correct the action of $\bA_\mu$ via the C-RandRAND preconditioner,
 \begin{equation} \label{eq:c_randrand_def}
    \bP = (\bI - \bPi) + \tau \bPi \bA_\mu^{-1} \bPi = (\bI - \bPi) + \tau \bPi \bOmega (\bA_\mu \bOmega)^\dagger  ,
    \end{equation}
that satisfies $ (\bA_\mu \bP) \bPi  = \tau \bA_\mu\bPi\bA_\mu^{-1}\bPi$.
Since $\bP$ in~\cref{eq:c_randrand_def} is spd whenever $\bA_\mu$ is spd, it can be directly integrated into classical preconditioned solvers for spd systems. 
In particular, the PCG algorithm is equivalent to applying CG to the preconditioned system:
\begin{equation}\label{eq:APdef2_2}
  \bP^{\frac{1}{2}} \bA_\mu\bP^{\frac{1}{2}}
  =
  \bP^{\frac{1}{2}} \bb.
\end{equation}

\Cref{thm:cnd_bound2} shows that C-RandRAND preconditioning effectively deflates the range of $\bA_\mu\bOmega$ from the range of $\bA_\mu$, yielding guarantees similar to R-RandRAND up to small constant factors. The coefficient $\rho$ minimizing the bound~\cref{eq:cnd_bound2} is 
\begin{equation}\label{eq:optimal_rho}
\rho = ( F \|\bE\| \lambda^{-1}_{\mathrm{max}}(\bPi \bA_\mu \bPi \bA_\mu^{-1} \bPi))^{-\frac{1}{2}},
\end{equation}
where $F = \lambda^{-1}_{\min}(\bA_\mu) = \|\bA_\mu^{-1}\|$. The quantities $\|\bE\|$ and $\lambda_{\mathrm{max}}(\bPi \bA_\mu \bPi \bA_\mu^{-1} \bPi)$ can be efficiently estimated with a few power iterations.  \blue{In~\cref{appendix_d} we present the refined bound~\cref{eq:cnd_bound2_ext} that can be (slightly) tighter than~\cref{eq:cnd_bound2} in scenarios where $\|\bPi \bA_\mu (\bI - \bPi)\|/\|\bE\|$ and $\|\bPi \bA_\mu^{-1} (\bI - \bPi)\|/\|\bF\|$ are small, which often occurs in practice. In our numerical experiments in~\cref{emp_eval}, we use $\rho$ minimizing this bound.}

If a sufficiently good lower bound for $\lambda_{\min}(\bA_\mu)$ is available, then we can select $F$ in~\cref{eq:optimal_rho} and~\cref{eq:cnd_bound2_ext}  based on this bound. 
For instance, when $\lambda_{\min}(\bA)\leq\mu$, which implies $\mu\leq\lambda_{\min}(\bA_\mu)\leq 2\mu$, we may simply set $F=\mu^{-1}$ and get a nearly optimal coefficient $\rho$.  
When $\lambda_{\min}(\bA)\geq\mu$ or this relation is unknown, estimating $\lambda_{\min}(\bA_\mu)$ directly may be difficult. In such scenarios, one can either select the coefficient $\rho$ in~\cref{eq:c_randrand_tau} as a constant (say $0.5$) or determine it using an empirical surrogate for $F$ like the one provided in~\cref{surrogate}. 

\blue{Finally, we note two alternative choices for the coefficient $\tau$ in C-RandRAND: 
\begin{equation}\label{eq:nys_tau}
    \tau = \|\bA_\mu^{-1} \bPi \|^{-1} = \| \bOmega (\bA_\mu \bOmega)^\dagger\|^{-1} \text{ and }  \tau = \|\bPi \bA_\mu^{-1} \bPi \|^{-1} = \| \bPi \bOmega (\bA_\mu \bOmega)^\dagger\|^{-1} ,
\end{equation}
which may be easier to compute than~\cref{eq:c_randrand_tau}. The first choice in~\cref{eq:nys_tau} satisfies 
$$\tau \lambda_{\max}(\bPi \bA_\mu \bPi \bA_\mu^{-1}\bPi) \leq \lambda_{l}(\bA_\mu) + \|\bPi \bA_\mu (\bI - \bPi)\|~~~\text{and}~\tau \geq \lambda_n(\bA_\mu),$$ and thus yields a quasi-optimal condition number, as follows from~\cref{thm:cnd_bound2}. The second choice of $\tau$ in~\cref{eq:nys_tau} represents the smallest eigenvalue of the Nystr{\"o}m approximation of $\bA_\mu$.\footnote{With this $\tau$, when $\mu = 0$, C-RandRAND mathematically recovers the Nystr{\"o}m preconditioner from~\cite{frangella2023randomized}.} This choice does not provide as strong theoretical guarantees compared to other presented choices of the coefficient $\tau$  but often performs well in practice (see~\cref{emp_eval}).}

\begin{proposition} 
\label{thm:cnd_bound2} 
    Let $\bP$ be a C-RandRAND preconditioner~\cref{eq:c_randrand_def} for spd matrix $\bA_\mu$. Denote $\bE = (\bI-\bPi)\bA_\mu(\bI-\bPi)$ and $\bF = (\bI-\bPi)\bA^{-1}_\mu(\bI-\bPi)$.  Let $\rho>0$ be some user-specified parameter. If $\tau$ is equal to
    \begin{equation} \label{eq:c_randrand_tau}
     \tau =  \rho \|\bE\|\lambda_{\mathrm{max}}^{-1}(\bPi \bA_\mu \bPi \bA_\mu^{-1} \bPi),
    \end{equation}
    then the eigenvalues of the preconditioned operator $\bP^{\frac{1}{2}} \bA_\mu \bP^{\frac{1}{2}}$ are bounded by
        \begin{equation} \label{eq:c_randrand_eigenvalues}
            (\|\bF\| +  \tau^{-1})^{-1} \leq \lambda_i(\bP^{\frac{1}{2}} \bA_\mu \bP^{\frac{1}{2}}) \leq (1+\rho)\|\bE\|.
        \end{equation}
    This in particular results in the following condition number bounds:
    \begin{equation} \label{eq:cnd_bound2}
    \begin{split}
    \cond(\bP^{\frac{1}{2}} \bA_\mu \bP^{\frac{1}{2}}) &\leq (1+\rho) \frac{\|\bE\|}{\lambda_{\mathrm{min}}(\bA_\mu)}  + (1+\rho^{-1}) \lambda_{\mathrm{max}}(\bPi \bA_\mu \bPi \bA_\mu^{-1} \bPi) \\
    %& \leq (1+\rho) \|\bE\|\|\bF\|  + (1+\rho^{-1}) \left (1 + \|\bE\|^{\frac{1}{2}}\|\bF\|^{\frac{1}{2}} \right )^2 \\
    & \leq (1+\rho) \frac{\|\bE\|}{\lambda_{\mathrm{min}}(\bA_\mu)} + (1+\rho^{-1}) \left (1 + \frac{\|\bE\|^{\frac{1}{2}}}{\lambda^{\frac{1}{2}}_{\mathrm{min}}(\bA_\mu)} \right )^2 .
    \end{split}
    \end{equation}
        \begin{proof}
            See~\cref{appendix_b}.
        \end{proof}
\end{proposition}

The combination of \cref{thm:cnd_bound2} with the randomized range approximation leads to condition number guarantees that are 
identical (up to constants) to those in the R-RandRAND scenario described in~\cref{thm:cnd_bound_spec}.

\subsection{G-RandRAND Preconditioner} \label{g_randrand}
The left G-RandRAND preconditioner,    
\begin{equation} 
\label{eq:g_randrand}
          \bP = (\bI - \bPi) + \tau \left  ((\bA_\mu ^{-1}\bPi )^\mathrm{T} \bA_\mu^{-1} \bPi  \right )^{\frac{1}{2}} = (\bI - \bPi) + \tau \left  (({\bA_\mu \bOmega)^\dagger}^\mathrm{T} (\bOmega^\mathrm{T} \bOmega) (\bA_\mu \bOmega)^\dagger  \right )^{\frac{1}{2}}  ,
\end{equation}
generalizes the C-RandRAND approach, as it can be applied to general systems. 
It is spd and admits the square root $\bP^\frac{1}{2} = (\bI - \bPi) + \tau^\frac{1}{2} \left  ((\bA_\mu ^{-1}\bPi )^\mathrm{T} \bA_\mu^{-1} \bPi  \right )^{\frac{1}{4}}$.
Thus, it can be employed in Krylov methods for symmetric indefinite systems, such as MINRES. 
For non-symmetric systems, $\bP$ in~\cref{eq:g_randrand} should be applied only as a left preconditioner.
If a right G-RandRAND preconditioner is required, it can be taken either as $\bA_{\mu}^{-1} \bP \bA_\mu$ or as the left G-RandRAND preconditioner for the transpose system $\bA^\mathrm{T}$.

The following result follows directly from the condition number bound~\cref{eq:cnd_bound2} for C-RandRAND.
\begin{proposition} \label{thm:cnd_bound3}
     Let $\rho>0$ be some user-specified parameter. Take $\tau$ as
     \begin{equation} \label{eq:g_randrand_tau}
     \tau = \rho^{\frac{1}{2}}  {\|(\bI - \bPi) \bA_\mu\|}{\|\bA_\mu^{-1} \bPi \bA_\mu \|^{-1}} .
    \end{equation}
     Then the preconditioned G-RandRAND operator $\bP \bA_\mu$ satisfies 
    \begin{equation} \label{eq:c_randrand_eigenvalues2}
            (\|(\bI-\bPi)\bA_\mu^{-1}\|^2 +  \tau^{-2})^{-1} \leq \sigma_i(\bP \bA_\mu )^2 \leq (1+\rho)\|(\bI-\bPi)\bA_\mu\|^2.
        \end{equation}
    This in particular results in the following condition number bounds:
    \begin{equation} \label{eq:cnd_bound22}
    \begin{split}
    \cond(\bP \bA_\mu)^2 &\leq (1+\rho) \frac{\|(\bI-\bPi)\bA_\mu\|^2}{\sigma_{\mathrm{min}}(\bA_\mu)^2}  + (1+\rho^{-1}) \|\bA_\mu^{-1} \bPi \bA_\mu \|^2 \\
    %& \leq (1+\rho) \|\bE\|\|\bF\|  + (1+\rho^{-1}) \left (1 + \|\bE\|^{\frac{1}{2}}\|\bF\|^{\frac{1}{2}} \right )^2 \\
    & \leq (1+\rho) \left (\frac{\|(\bI-\bPi)\bA_\mu\|}{\sigma_{\mathrm{min}}(\bA_\mu)} \right )^2  + (1+\rho^{-1}) \left (1 + \frac{\|(\bI-\bPi)\bA_\mu\|}{\sigma_{\mathrm{min}}(\bA_\mu)} \right )^2.
    \end{split}
    \end{equation}
        \begin{proof}
            See~\cref{appendix_b}.
        \end{proof}
\end{proposition}
The coefficient $\rho$ in~\cref{eq:g_randrand_tau} can be selected a priori (e.g., as $0.5$) or based on a value that minimizes the condition number bound~\cref{eq:cnd_bound22}:
\begin{equation} \label{eq:rho_g_rand_rand}
\rho = F^{-1} \|(\bI-\bPi)\bA_\mu\|^{-1} \|\bA_\mu^{-1} \bPi \bA_\mu \|,
\end{equation}
where $F = \sigma^{-1}_{\min}(\bA_\mu) = \|\bA_\mu^{-1}\|$.

The quantities $\|(\bI-\bPi)\bA_\mu\|$ and $\|\bA_\mu^{-1} \bPi \bA_\mu\|$ can be efficiently estimated with a few power iterations. The coefficient $F$ can be selected based on quasi-optimal lower bounds of the smallest singular value.

\Cref{thm:cnd_bound3} can be directly integrated with the randomized range finders. 
In the following proposition, we present a quasi-optimality guarantee of G-RandRAND deflation for symmetric matrices. 
\begin{proposition} 
\label{thm:cnd_bound_spec_ext}
  Let $\bA_\mu$ be symmetric. Define $\mu_{\mathrm{rel}} = |\mu|/\sigma_{\mathrm{min}}(\bA_\mu)$. Draw a Gaussian matrix $\bX \in \mathbb{R}^{l \times n}$ with sketching dimension $l = 2ck+4$, where $1 < c \leq 2$ and $k \geq 1$. Let $\bP$ be the G-RandRAND preconditioner from \cref{eq:g_randrand}, where $\bPi$ is the orthogonal projector onto the basis $\bA_\mu\,\bA^q\,\bX^\mathrm{T}$. Let $\tau$ satisfy \cref{eq:g_randrand_tau} with $\rho=1$. Then we have with probability at least $1-6^{-k}$ for $q=0$:
        \begin{equation}\label{eq:condbound1_g_randrand}
        \cond(\bP \bA_\mu) \leq 2+6\mu_{\mathrm{rel}} + 4\left( C_3 \frac{n}{(c-1)k}\right)^{\frac{1}{2}} \cond_k(\bA_\mu^2)^\frac{1}{2} ,
        \end{equation}
        and for $q \geq 1$:
        \begin{equation}\label{eq:condbound2_g_randrand}
         \mathrm{cond}(\bP \bA_\mu) \leq 2+  6\mu_{\mathrm{rel}}    + (4+ 4\mu_{\mathrm{rel}}) \left ( { C_3 \frac{n}{(c-1)k}} \right )^\frac{1}{2q+2}  \cond_{k,|\mu|^{2q+2}}(\bA_\mu^2\bA^{2q})^\frac{1}{2q+2} ,
        \end{equation}
        where $C_3 = 2 \cdot 18^2$. Furthermore, in expectation,~\cref{eq:condbound1_g_randrand,eq:condbound2_g_randrand} hold with smaller constant: $C_3 = 4 \mathrm{e}^2$.

        \begin{proof}
            See~\cref{appendix_b}.
        \end{proof}
       
\end{proposition} 

\noindent 
\blue{When the relative shift $\mu_{\mathrm{rel}} = |\mu|/\sigma_{\min}(\bA_\mu)$ is $\mathcal{O}(1)$, \Cref{thm:cnd_bound_spec_ext} guarantees the same bounds, up to constant factors, as those established for R-RandRAND and C-RandRAND in the spd case. Similar guarantees also extend to SRHT and sparse OSEs $\bX$ of size $l = \mathcal{O}(k \log(n/\delta))$.} Furthermore, similarly to the spd scenario, the $k$-stable condition numbers $\cond_{k}(\bA_\mu^2)$ and $\cond_{k, |\mu|^{2q+2}}(\bA_\mu^2 \bA^{2q})$ appearing in \cref{thm:cnd_bound_spec_ext} admit the bounds 
    $$ \cond_{k}(\bA_\mu^2)^{\frac{1}{2}} \leq  \frac{\sigma_k(\bA_\mu)}{\sigma_\mathrm{min}(\bA_\mu)} ~ \text{ and } \cond_{k, |\mu|^{2q+2}}(\bA_\mu^2 \bA^{2q})^{\frac{1}{2q+2}} \leq  2 \frac{\sigma_k(\bA)+|\mu|}{\sigma_\mathrm{min}(\bA_\mu)}.$$
Here, $\cond_{k}(\bA_\mu^2)^{\frac{1}{2}}$ again behaves as a standard averaged condition number and is always at least $1$.  
\blue{In contrast, $\cond_{k, |\mu|^{2q+2}}(\bA_\mu^2 \bA^{2q})^{\frac{1}{2q+2}}$ may fall below $1$ resulting in a constant bound on $\mathrm{\bP \bA_\mu}$. Specifically, let $N$ be the number of singular values $\sigma_i(\bA)$ such that $\sigma_i(\bA)^{2q} \leq \alpha |\mu|^{2q} $ for some (relatively small) $\alpha \leq 1$. Then we have the following bound 
$$\cond_{k, |\mu|^{2q+2}}(\bA_\mu^2 \bA^{2q})^{\frac{1}{2q+2}} \leq  \frac{4 \max{\left ( F_1  |\mu|,  F_2 (\sigma_k+|\mu|)\right )}}{\sigma_\mathrm{min}(\bA_\mu)}, $$
where $F_1 = \left ( \frac{4\alpha N}{n-k+1} \right )^{\frac{1}{2q+2}}$ and $F_2 = \left (1- \frac{N}{n-k+1}  \right  )^{\frac{1}{2q+2}}$.
Consequently, when the spectrum of $\bA$ decays sufficiently fast so that $N$ is close to $n-k$, the $k$-stable condition number $\cond_{k, |\mu|^{2q+2}}(\bA_\mu^2 \bA^{2q})^{\frac{1}{2q+2}}$ can become significantly smaller than $\frac{\sigma_k(\bA)+|\mu|}{\sigma_\mathrm{min}(\bA_\mu)}$.}
\begin{remark} \label{rmk:g_rand_rand_smallest}
When $\bA_\mu$ is indefinite, its smallest singular value $\sigma_{\min}(\bA_\mu)$ may be well below $|\mu|$. G-RandRAND has the ability to lift it to $\sigma_{\min}(\bP \bA) \approx |\mu|$ (along with deflating the dominant singular values). Specifically,  when the sketch size $l$ is sufficiently large and $q\geq 1$ so that the projection error satisfies $\|(\bI - \bPi)\bA\| \leq \beta|\mu|$ for some $\beta \leq 1$,\footnote{In symmetric case, we have $\|(\bI - \bPi)\bA \| \leq \|(\bI - \bPi^*)\bA \| + |\mu| \|(\bI - \bPi)\bPi^*\|$, where $\bPi^*$ denotes the orthogonal projector onto $\mathrm{range}(\bA^q\bOmega)$. Using the results from~\cref{range_approx} with zero shift and power $q-1$, one obtains that $\|(\bI - \bPi^*)\bA \|$ is, with high probability, bounded by $\mathcal{O}\left( \left( \frac{\mathrm{sr}_k(\bA^{2q})}{\sqrt{k}} \right)^{1/q} \sigma_k(\bA) \right)$, where $k = l/2 - 1$ for Gaussian sketching matrices. Consequently, $\|(\bI - \bPi)\bA\|$ can fall below $|\mu|$ when $\sigma_k(\bA) \leq \mathcal{O}\left( \left( \frac{\mathrm{sr}_k(\bA^{2q})}{\sqrt{k}} \right)^{-1/q} |\mu| \right)$.
} it holds
$$
\|(\bI - \bPi)\bA_\mu^{-1}\| =|\mu|^{-1} \| (\bI-\bPi) - (\bI - \bPi)\bA \bA_\mu^{-1}\| \leq |\mu|^{-1} +\beta \sigma_{\min}(\bA_\mu)^{-1}.
$$
This bound can be combined with the inequality
$$
\tau^{-1} \leq \rho^{-1/2} \left( \|(\bI - \bPi)\bA_\mu^{-1}\| + \|(\bI - \bPi)\bA_\mu\|^{-1} \right),
$$
and substituted into~\cref{eq:c_randrand_eigenvalues2} to obtain a desired bound of $\sigma_{\min}(\bP \bA)$. In particular, when $\beta \sigma_{\min}(\bA_\mu)^{-1}$ is comparable to $|\mu|^{-1}$, G-RandRAND achieves a constant condition number bound for $\bP \bA$. 

Furthermore, in the symmetric case, the quantity $\|(\bI - \bPi)\bA_\mu^{-1}\|$ can be estimated with the following efficiently computable bound
$$
\|(\bI - \bPi)\bA_\mu^{-1}\| \leq (1-\beta)^{-1}|\mu|^{-1} (1+ \|(\bI - \bPi)\bA_\mu^{-1} \bA \bPi\|) = (1-\beta)^{-1} \left (|\mu|^{-1} + \|(\bI - \bPi)\bA_\mu^{-1}\bPi\| \right ) .
$$

\end{remark}

\begin{remark} \label{rmk:g_randrand_vs_c_randrand} 
\blue{Besides G-RandRAND, a preconditioner for the indefinite matrix $\bA_\mu$ can also be constructed by taking the square root of the C-RandRAND preconditioner associated with $\bA_\mu^2$ with zero shift. However, G-RandRAND offers several advantages over C-RandRAND.  First, for the same computational cost, G-RandRAND allows either doubling the sketch size or increasing the power parameter $q$. Second, it is significantly more numerically stable, as it avoids squaring the condition number of the operator. Finally, according to~\cref{eq:condbound2_g_randrand},  when $q \geq 1$, G-RandRAND can effectively exploit the decay of the trailing singular values, unlike C-RandRAND on $\bA_\mu^2$ with zero shift. See~\cref{emp_eval} for empirical comparisons.} 
\end{remark}

\blue{\subsection{Additional RandRAND Preconditioners}}
\blue{Other choices for RandRAND-based precondiitoners are possible. Examples include
$$
\bP = (\bI-\bPi) + \tau (\bPi \bA_\mu \bPi)^\dagger,  
\qquad  
\bP = (\bI-\bPi) + \tau \left( ((\bA_\mu \bPi)^\mathrm{T} \bA_\mu\bPi)^\dagger \right)^\frac{1}{2}.  
$$
The advantage of these variants is that they remain computationally feasible for general approximation subspaces, not only those of the form $\bA_\mu \bOmega$. However, their construction and application (in basis-less form) may incur higher costs and lower numerical stability than the preconditioners discussed in~\cref{sec:precond,sec:c_precond,g_randrand}. The analysis of such alternatives is left outside the scope of this work.}
    \section{Explicit and Basis-less Implementation} 
\label{orthogonalization}

In this section, we discuss implementation details, including basis construction aspects and numerical stability.

\subsection{Explicit Basis Orthogonalization}
\label{basis_explicit_represent}

RandRAND preconditioners can be implemented through a QR (or SVD) factorization of $\bA_\mu \bOmega$. 
Let 
    \begin{equation} \label{eq:Rdef}
     \bA_\mu \bOmega = \bQ \bR,
    \end{equation}
    where $\bQ$ is orthonormal and $\bR$ is upper triangular or other easily invertible matrix. 
    Then, the products with the preconditioners and preconditioned operators can be computed using the relations
    \begin{equation} \label{eq:projector}
        \bPi \bu = \bQ \bQ^\mathrm{T} \bu ~~~\text{ and } \bA_\mu^{-1} \bPi \bu = \bOmega \bR^{-1} \bQ^\mathrm{T} \bu.
    \end{equation}
    For instance, we can use the following representations of the R-RandRAND preconditioned operator~\cref{eq:APdef2}: 
    \begin{equation} \label{eq:r_randrand_explicit}
          \bA \bP = (\bI - \bQ\bQ^\mathrm{T}) \bA_\mu (\bI - \bQ\bQ^\mathrm{T})  + \tau \bQ \bQ^\mathrm{T} , 
    \end{equation}  
    and the C-RandRAND~\cref{eq:c_randrand_def} and G-RandRAND~\cref{eq:g_randrand} preconditioners:
    \begin{equation} \label{eq:g_randrand_explicit}
          \bP = \bI - \bQ \bQ^\mathrm{T} + \tau \bQ \bG \bQ^\mathrm{T}, 
    \end{equation}
    where $\bG = \bR^{-\mathrm{T}} (\bOmega^{\mathrm{T}} \bA_\mu \bOmega) \bR^{-1} $ for C-RandRAND and  $\bG = \left(\bR^{-\mathrm{T}} ( \bOmega^\mathrm{T} \bOmega) \bR^{-1} \right )^\frac{1}{2}$ for G-RandRAND. 
    In the latter scenarios, the matrix $\bG$ is $l \times l$ and can be efficiently precomputed and applied through an SVD or Cholesky factorization for better numerical stability. 
     
    \paragraph{Computational cost.} 
    Assuming that $\bA$ is symmetric and that each matvec with $\bA$ takes $T_\mathrm{mv}$ flops, constructing the preconditioner via an explicit QR factorization requires\footnote{excluding orthogonalization costs in subspace iteration for computing $\bOmega$} $$T_\mathrm{mv}l(q+1) + nl^2 + \mathcal{O}(l^3) \text{ flops and } nl~\text{memory,} $$
    and applying the preconditioned operator to a vector requires $T_\mathrm{mv} + \mathcal{O}(nl)$ flops. 
    This approach should be considered when the cost of QR factorization and maintaining and applying the Q factor is minor compared to other operations (e.g., when the sketching dimension $l$ is small and $\bA_\mu$ is dense).

    \subsection{Basis-less Implementation} \label{basis_free_represent}
     In general, explicit orthogonalization of $\bA_\mu \bOmega$ and matvecs with the Q factor can be computationally expensive and may dominate both complexity and memory usage. 
     In such cases, it can be more efficient to compute the transform $\bR$ with a fast Q-less factorization and operate with $\bQ$ in~\cref{eq:projector,eq:r_randrand_explicit,eq:g_randrand_explicit} implicitly through matvecs $\bu \mapsto \bQ \bu$ and $\bu \mapsto \bQ^\mathrm{T} \bu$, using
      \begin{equation} \label{eq:implicit}
     \bQ \bu  = \bA_\mu \bOmega \bR^{-1} \bu
    \end{equation}
    and
    \begin{equation} \label{eq:implicit2}
    \bQ^{\mathrm{T}}\bu = \bR^{-\mathrm{T}}  \bOmega^\mathrm{T} \bA_\mu^\mathrm{T} \bu.
    \end{equation}

    \paragraph{Computational cost.} 
    
    For better scaling of the computational cost with $l$, we can choose the sketching matrix $\bX$, which defines $\bOmega$, to be a fast randomized transform such as SRHT, sparse OSE, or multi-level OSE, which require only $\mathcal{O}(n \log n)$ memory and flops to apply to a vector.
    For the rest of this section, we assume that $\bA$ is symmetric; the non-symmetric case is analogous. In this setting, a matvec with $\bA_\mu \bOmega$ or $(\bA_\mu \bOmega)^{\mathrm{T}}$ can be performed by separate application of $\bA_\mu \bA^q = (\bA_\mu \bA^q)^\mathrm{T}$ at a cost of $(q+1) T_\mathrm{mv}$ flops, and $\bX$ or $\bX^\mathrm{T}$ at a cost of $\mathcal{O}(n \log n)$ flops. This can be substantially cheaper than the $n l$ flops required for multiplying a vector with an explicit dense $n \times l$ matrix.
    Furthermore, the use of fast-to-apply sketching matrices $\bX$ also reduces the cost of constructing the preconditioner. In the basis-less approach, constructing the preconditioner amounts to computing the factor $\bR$. In turn, the methods below allow to compute $\bR$ with only
    $$
    \mathcal{O}\bigl(T_\mathrm{mv}  l  q + n  l \log n + l^3\bigr) ~\text{flops and } 
    \mathcal{O}(n \log n)~\text{memory}.
    $$

\paragraph{Fast Q-less Cholesky QR factorizations.}

The $\bR$ factor can be computed with a Cholesky factorization of the Gram matrix $\bOmega^\mathrm{T}  \bA_\mu \bA_\mu  \bOmega$, as depicted in~\cref{alg:chol_qr}.
Using fast or sparse OSE $\bOmega$, the product $\bOmega^{\mathrm{T}} \bA_\mu \bA_\mu \bOmega = \bX \bA^q \bA_\mu \bA_\mu \bA^q \bX^\mathrm{T}$ in step 1 of~\cref{alg:chol_qr} can be computed columnwise from right to left in $2(q+1)T_\mathrm{mv}l + \mathcal{O}(n  l\log n)$ time. 

However,~\cref{alg:chol_qr} may introduce significant rounding errors, leading to orthogonality loss of $\bQ$. 
Specifically, in finite precision arithmetic, the orthogonality measure $\|\bI - \bQ^\mathrm{T} \bQ \|$ in Cholesky QR scales as $ \mathrm{u} \,  \mathrm{cond}(\bA_\mu \bOmega)^2$, where $\mathrm{u} \,$ is the unit roundoff~\cite{yamamoto2015roundoff}, and it can become dramatically large if $\bA_\mu \bOmega$ is not sufficiently well conditioned. 
This issue can be alleviated by employing second-level preconditioning, as discussed in the following. 

\begin{algorithm} 
        \caption{Fast Q-less Cholesky QR} \label{alg:chol_qr}
        \begin{algorithmic} [1]
            \Statex{\textbf{function} $\bR =\mathtt{fast{\_}q{\_}less{\_}chol{\_}qr}(\bA_\mu, \bOmega)$} 
            \State Compute the Gram matrix using the structure of $\bOmega$: $\bG \mapsfrom  \bOmega^\mathrm{T}  \bA_\mu \bA_\mu  \bOmega$. 
            \State Return Cholesky factor: $\bR \mapsfrom \mathrm{chol}(\bG)$.
        \end{algorithmic}
\end{algorithm}

\Cref{alg:prec_chol_qr} describes a Q-less orthogonalization of $\bA_\mu \bOmega$ that uses a second-level randomized preconditioner $\bR_{\mathrm{sk}}^{-1}$ such that $\bA_\mu \bOmega\bR_{\mathrm{sk}}^{-1}$ is well-conditioned with high probability. The transform $\bR_{\mathrm{sk}}$ is taken as the R factor from a QR factorization of the sketched matrix $\bTheta \bA_\mu \bOmega$, where $\bTheta$ is an additional fast randomized transform of small size. 
This approach is connected to RandNLA methods for solving least-squares problems~\cite{rokhlin2008fast,AMT10,MSM14_SISC} and explicit orthogonalization~\cite{balabanov2020randomized,balabanov2022randomized,melnichenko2023choleskyqr,higgins2023analysis}.

It can be shown that, if $\bTheta$ satisfies the $\varepsilon$-embedding property for some $\varepsilon<1$,
\begin{equation} \label{eq:epsilon_embedding}
(1-\varepsilon) \|\bA_\mu \bOmega \bx\|^2 
\le  \|\bTheta \bA_\mu \bOmega \bx\|^2 
\le (1+\varepsilon) \|\bA_\mu \bOmega \bx\|^2,
\end{equation}
then $\cond(\bA_\mu \bOmega \bR_{\mathrm{sk}}^{-1}) \leq \sqrt{\frac{1+\varepsilon}{1-\varepsilon}}$. 
This can guarantee that the orthogonality measure $\| \bI - \bQ^\mathrm{T}\bQ \|$ associated with the transform $\bR$ computed via~\cref{alg:prec_chol_qr} scales as $ \mathrm{u} \,  \cond(\bA_\mu \bOmega)$, ensuring numerical stability even for ill-conditioned operators.

It is well known that SRHT matrices with $l_{\mathrm{sk}} = {\mathcal{O}}\bigl(l \log{n} \bigr)$ rows and sparse embeddings with $l_{\mathrm{sk}} = {\mathcal{O}}\bigl(l \log{n}  \bigr)$ rows and $\gamma_{\mathrm{sk}} = {\mathcal{O}}(\log{n})$ nonzero entries per column, satisfy the $\varepsilon$-embedding property for $\varepsilon = \frac{1}{2}$ with high probability~\cite{cohen2016nearly,woodruff2014sketching,halko2011finding,boutsidis2013improved}.
When $\bOmega$ is taken as $\bA^q\bX$ and the matrix products are computed from right to left, \Cref{alg:prec_chol_qr} runs in
$
3(q+1)T_\mathrm{mv} l 
+ {\mathcal{O}} \bigl (nl\log{n} + l^3 \log{n} \bigr)
$ time. Furthermore, the term $l^3 \log{n}$ can be reduced to $l^3$ by selecting $\bTheta$ as a multi-level embedding or an optimal sparse OSE, requiring sketching dimension $l_{\mathrm{sk}} = {\mathcal{O}}\bigl(l \bigr)$.
    
    \begin{algorithm} 
        \caption{Fast Q-less preconditioned Cholesky QR} \label{alg:prec_chol_qr}
        \begin{algorithmic} [1]
            \Statex{\textbf{function} $\bR =\mathtt{fast{\_}q{\_}less{\_}precond{\_}chol{\_}qr}(\bA_\mu, \bOmega)$}
            \State Draw sketching matrix $\bTheta$ and compute the sketch using the structures of $\bTheta, \bOmega$: $\bS \mapsfrom  \bTheta \bA_\mu \bOmega$. 
            \State Compute QR factorization of the sketch: $[\sim, \bR_{\mathrm{sk}}] \mapsfrom \mathrm{qr}(\bS)$. 
            \State Compute the preconditioned Gram matrix: $\bG \mapsfrom \bR_{\mathrm{sk}}^\mathrm{-T} \bOmega^\mathrm{T}  \bA_\mu \bA_\mu  \bOmega \bR_{\mathrm{sk}}^{-1}$.
            \State Compute Cholesky factor of $\bG$: $\bR_{\mathrm{chol}} \mapsfrom \mathrm{chol}(\bG)$.
            \State Return: $\bR \mapsfrom  \bR_{\mathrm{sk}} \bR_{\mathrm{chol}}$
        \end{algorithmic}
    \end{algorithm}

    \begin{remark}
        Considering large powers $q$ in the above fast Q-less QR methods can lead to severe numerical stability issues, due to the dependence of rounding errors on $\mathrm{cond}(\bA_\mu \bOmega) = \mathrm{cond}(\bA_\mu \bA^q \bX^\mathrm{T})$. 
        We suggest using them only for powers like $q=0$ or $q=1$.
    \end{remark}
    \begin{remark}
        When $\bA_\mu \bOmega$ is ill-conditioned, it can be truncated by replacing $\bOmega$ with $ \bOmega \bT$, where $\bT$ is orthogonal and truncates the singular values of $\bA \bOmega$ smaller than some tolerance. 
        Such a transform $\bT$ can be efficiently obtained with a low-rank approximation of the two-sided sketch $\bTheta \bA \bOmega$. 
    \end{remark}
    \begin{remark} \label{rmk:operator_shift}
        The unregularized matrix $\bA$ may have a high condition number, which can severely impact numerical stability, even when $\bA_\mu$ is not too badly conditioned. 
        In such cases, it may be pertinent to apply a shift to $\bA$ and $\mu$:   
        $\bA \mapsfrom \bA + \alpha \bI$ and $\mu \mapsfrom \mu - \alpha$, 
        where the shift $\alpha$ ensures that $\bA$ is sufficiently well conditioned. 
        This transformation leaves the regularized system $(\bA+\mu\bI)\bx = \bb$ unchanged. 
    \end{remark}

\subsection{Numerical Stability} \label{num_stab}

\paragraph{Re-orthogonalization.}

In finite precision arithmetic, it can be important to ensure numerical orthogonality between the two components of RandRAND. 
This can help to preserve crucial properties of the preconditioner and the preconditioned operator (e.g., numerical symmetry or positive definiteness), and it can be achieved through \textit{re-orthogonalization.}

Re-orthogonalization can be implemented by expressing the projection $(\bI - \bPi)\bu$ onto the complement of $\mathrm{range}(\bA_\mu \bOmega)$ as
$$
(\bI - \bPi)\bu = (\bI - \bPi)(\bI - \bPi)\bu.
$$
For example, re-orthogonalization modifies the R-RandRAND preconditioned operator to
\begin{equation} \label{eq:r_randrand_reorth}
\bA \bP
= (\bI - \bPi)(\bI - \bPi)\bA_\mu(\bI - \bPi)(\bI - \bPi)
+
\tau\bigl(\bI - (\bI - \bPi)(\bI - \bPi)\bigr).
\end{equation}
Whereas, for  C-RandRAND it can be incorporated as follows 
\begin{equation} \label{eq:c_randrand_reorth}
\bP
= (\bI - \bPi)(\bI - \bPi)
+
\tau\bigl(\bI - (\bI - \bPi)(\bI - \bPi)\bigr)\bA_\mu^{-1} \bPi.
\end{equation}
%
%Re-orthogonalization can also be important to ensure that the computed projections $\bv = \bA_\mu^{-1} \bPi \bu$ are such that $\bA_\mu \bv$ lie in $\mathrm{range}(\bA_\mu \bOmega)$. 
%This can be achieved by using the following identity:
%$$\bv = \bv -  (\bI- \bA_\mu^{-1} \bPi \bA_\mu)(\bI- \bA_\mu^{-1} \bPi \bA_\mu)\bv.$$

In some scenarios, re-orthogonalization can entail significant computational overhead. 
This cost can be mitigated by applying the procedure only every few iterations of the iterative solver.

\paragraph{Projectors.}

In practice, the basis matrix $\bQ$ is orthogonalized to some finite precision. 
For instance, in the basis-less approach, the measure $\|\bI - \bQ^\mathrm{T}\bQ\|$ can be on the order of $\mathrm{u} \,\cond(\bA_\mu \bOmega)^2$ or $\mathrm{u} \,\cond(\bA_\mu \bOmega)$, depending on the fast QR scheme used. 
Consequently, computing the projections as
$$
\bPi \bu = \bQ \bQ^\mathrm{T} \bu
\quad\text{and}\quad
\bA_\mu^{-1} \bPi \bu = \bOmega \bR^{-1} \bQ^\mathrm{T} \bu
$$
can results in additional approximation errors, which can impact the effectiveness of the preconditioner. 
If needed, these errors can be reduced by refining the projector. 
One way to achieve this is to represent $\bPi$ and $\bA_\mu^{-1} \bPi$ as
\begin{equation}\label{eq:accurate_proj}
\bPi \bu = \bQ \bQ^\dagger \bu = \bQ (\bQ^\mathrm{T} \bQ )^{-1} \bQ^\mathrm{T} \bu \approx \bQ \bT \bQ^\mathrm{T} \bu
\quad\text{and}\quad
\bA_\mu^{-1} \bPi \bu \approx \bOmega \bR^{-\mathrm{T}} \bT \bQ^\mathrm{T} \bu,
\end{equation}
where $\bT$ approximates $(\bQ^\mathrm{T}\bQ)^{-1}$, based on the fact that $\bQ$ is nearly orthonormal. 
A simple choice is the Neumann-series approximation $\bT = 2\bI - \bQ^\mathrm{T}\bQ$, which can reduce the errors associated with computing $\bPi \bu$ and $\bA_\mu^{-1}\bPi \bu$ from the order of $\|\bI - \bQ^\mathrm{T}\bQ\|$ to the order of $\|\bI - \bQ^\mathrm{T}\bQ\|^2$. An alternative way is compute $\bQ^\dagger \bu$ through a few refinement iterations, resulting in the following expressions for $\bPi$ and $\bA_\mu^{-1} \bPi$:
\begin{equation}\label{eq:accurate_proj_2}
\bPi \bu \approx \bQ \bT(\bu,\bQ^\mathrm{T}\bu),
\qquad
\bA_\mu^{-1} \bPi \bu \approx \bOmega \bR^{-\mathrm{T}} \bT(\bu,\bQ^\mathrm{T}\bu),
\end{equation}
where $\bT(\bu,\bQ^\mathrm{T}\bu) = \bv_k$ is given by the recurrence $\bv_{1} = \bQ^\mathrm{T}\bu$, 
$$\bv_{i+1} = \bv_i+ \bQ^\mathrm{T}(\bu - \bQ\bv_i ).$$
In some scenarios this approach can be more numerically robust than using~\cref{eq:accurate_proj}.

Similarly to the orthogonalization, the refinement of the projector can have significant impact on computational cost. 
For instance, the approximation~\cref{eq:accurate_proj} with $\bT = 2\bI - \bQ^\mathrm{T}\bQ$ essentially doubles the cost of applying $\bPi$.

    \section{Adaptations} \label{adaptations}

\blue{RandRAND preconditioners are inherently flexible and can be effectively adapted to a variety of practical scenarios and computational architectures. In this section, we describe several such adaptations.}

\subsection{Multiple Shifting Parameters} \label{mult_shifts}

One often needs to compute solutions to the system~\cref{eq:initsys} across a range of shifting parameters $\mu$.   
In statistics, for example, analyzing the regularization path provides insights into the trade-off between bias and variance, identifying regions where a model avoids underfitting or overfitting. 
A related setting occurs when computing eigenvalues of $\bA$. 
The standard technique for this task is shift-and-invert, where for multiple shifts $\mu$ (often chosen near $-\lambda_i$), one computes the dominant eigenvalue(s) of $\bA_\mu^{-1}$ via power iterations or other methods requiring linear solves with $\bA_\mu$. \blue{Solving systems with multiple shifting parameters is also needed when computing matrix functions $f(\bA_\mu) \bb$ in integral form.}
Regularization and \blue{shifting} sweeps can be expensive, especially when evaluating a wide range of $\mu$ values; and so it is therefore important to recycle computations when building preconditioners. 
This scenario is discussed next. 
    
We can precompute the sketch $\bA \bOmega$ and then, for each value of $\mu$, obtain the deflation basis $\bQ$ and the transform $\bR$, with a QR of $\bA \bOmega+ \mu\bOmega$:  
        $$[\bQ, \bR] \leftarrow  \mathrm{qr}(\bA \bOmega+ \mu\bOmega).$$
        In some cases, computing multiple QR factorizations of $n \times l$ matrices can be relatively expensive. 
        \blue{This issue can be circumvented with the following approach. First, we precompute the joint basis matrix $[\bQ_\mathrm{joint}, \sim] \leftarrow  \mathrm{qr}([\bOmega~\bA \bOmega])$ and precompute $\bQ_\mathrm{joint}^\mathrm{T}\bA \bOmega$ and $\bQ_\mathrm{joint}^\mathrm{T} \bOmega$.  Then for each $\mu$ we compute QR factorization of an $2l \times l$ reduced matrix: $[\bQ_\mathrm{red}, \bR] \leftarrow  \mathrm{qr}(\bQ_\mathrm{joint}^\mathrm{T}\bA \bOmega+ \mu\bQ_\mathrm{joint}^\mathrm{T}\bOmega)$ and use\footnote{The explicit matrix product $\bQ_\mathrm{joint} \bQ_\mathrm{red}$ can be avoided. We can operate with  $\bQ$ implicitly via matvecs $\bQ \bv \leftarrow \bQ_\mathrm{joint} (\bQ_\mathrm{red} \bv)$.}
        $$ \bQ  =  \bQ_\mathrm{joint} \bQ_\mathrm{red}. $$
        }

        In the basis-less scenario from~\cref{basis_free_represent}, computing an explicit QR of the (joint) basis is considered to be computationally intractable. To address this point we can
        exploit the parameter separation of the Gram matrix:
            $$ (\bA \bOmega + \mu \bOmega)^\mathrm{T}(\bA \bOmega + \mu \bOmega) = \bOmega^\mathrm{T} \bA^\mathrm{T} \bA \bOmega + 2 \mu \bOmega^\mathrm{T} \bA \bOmega + \mu^2\bOmega^\mathrm{T}  \bOmega.$$
        This decomposition allows the computation of $\bR$ from small, parameter-independent projections $\bOmega^\mathrm{T} \bA^\mathrm{T} \bA \bOmega$, $\bOmega^\mathrm{T} \bA \bOmega$, and $\bOmega^\mathrm{T} \bOmega$:
            $$ \bR \leftarrow \mathrm{chol} \left( \bOmega^\mathrm{T} \bA^\mathrm{T} \bA \bOmega + 2 \mu \bOmega^\mathrm{T} \bA \bOmega + \mu^2 \bOmega^\mathrm{T}  \bOmega \right).$$
        However, as was already mentioned in~\cref{basis_free_represent}, in finite precision arithmetic, the Cholesky QR leads to orthogonality loss of the order of $\mathrm{u}\, \mathrm{cond}(\bA_\mu \bOmega)^2$ leading to stability issues when $\mathrm{cond}(\bA_\mu \bOmega) \geq  \mathcal{O}(\mathrm{u}^{-\frac{1}{2}})$. 
    
        For spd operators and $\mu \geq 0$, the numerical stability issues of the Q-less Cholesky QR can be alleviated by using a shifted Cholesky preconditioning technique~\cite{fukaya2020shifted}. 
        In this approach, we compute an R factor $\bR_\alpha$ for $\bA_\alpha \bOmega$ for a fixed shift $\alpha = c\mathrm{u}^{\frac{1}{2}}\|\bA\|$. 
        For each regularization parameter $\mu \leq \alpha$, the orthogonalization transform $\bR$ for  $\bA_\mu \bOmega$ is then obtained as:
            $$\bR = \bR_{\mathrm{chol}} \bR_{\alpha}, $$
        where $\bR_{\mathrm{chol}}$ is the Cholesky \blue{factor} of the preconditioned Gram matrix $(\bOmega \bR_{\alpha}^{-1})^\mathrm{T} \bA_\mu^\mathrm{T} \bA_\mu (\bOmega \bR_{\alpha}^{-1})$:
            $$ \bR_{\mathrm{chol}} \leftarrow \mathrm{chol} \left( (\bOmega \bR_{\alpha}^{-1})^\mathrm{T} \bA_\alpha^\mathrm{T} \bA_\alpha (\bOmega \bR_{\alpha}^{-1}) + 2 (\mu-\alpha) (\bOmega \bR_{\alpha}^{-1})^\mathrm{T} \bA_\alpha (\bOmega \bR_{\alpha}^{-1}) + (\mu^2+\alpha^2) (\bOmega \bR_{\alpha}^{-1})^\mathrm{T}  (\bOmega \bR_{\alpha}^{-1}) \right ). $$  
        As before, the projections $(\bOmega \bR_{\alpha}^{-1})^\mathrm{T} \bA_\alpha^\mathrm{T} \bA_\alpha (\bOmega \bR_{\alpha}^{-1})$, $(\bOmega \bR_{\alpha}^{-1})^\mathrm{T} \bA_\alpha (\bOmega \bR_{\alpha}^{-1})$, $(\bOmega \bR_{\alpha}^{-1})^\mathrm{T}  (\bOmega \bR_{\alpha}^{-1}) $ can be precomputed and used for multiple $\mu$.

\subsection{Representation via Least-Squares Solves} \label{ls_solves}

\paragraph{Second-level Preconditioning.}
    
    \blue{In~\cref{orthogonalization} and~\cref{mult_shifts},} 
    we considered implementations  based on orthogonalization of $ \bA_\mu \bOmega$.
    On the other hand, RandRAND preconditioning is compatible with transforms $\bR$
   which orthogonalize $\bA_\mu$ only approximately and hence may be cheaper to compute than the exact orthogonalization factors.

    In such cases, matvecs with preconditioner $\bP$ and the preconditioned operators can be obtained by using the relations 
    $$ \bPi \bu = \bQ \bQ^\dagger\bu,~~~\text{and } \bA_\mu^{-1}\bPi \bu = \bOmega \bR^{-1} \bQ^\dagger \bu, $$
    where the matvecs with $\bQ^{\dagger}$ are obtained via a least-squares solve:
    \begin{equation} \label{eq:ls_representation} 
        \bQ^{\dagger}\bu = \arg \min_\bv \|\bQ \bv - \bu \|.
    \end{equation}
    Here, a trade-off arises between the cost of constructing the preconditioner (i.e., computing $\bR$) and the cost of applying the preconditioned operator to  vectors (i.e., computing  $\bQ^{\dagger} \bu$). 
    When $\bQ$ is well-conditioned, the least-squares problem in~\eqref{eq:ls_representation} can be efficiently solved using a few iterations of iterative solvers such as LSQR or LSRN. 
    
    In the basis-less representation, the identity~\cref{eq:ls_representation} becomes  
    \begin{equation} \label{eq:ls_representation2}
        \bQ^{\dagger}\bu = \arg \min_\bv \| \bA_\mu \bOmega \bR^{-1} \bv - \bu \|. 
    \end{equation}  
    Thus, the matrix $\bR^{-1}$ can here be viewed as a second-level preconditioner used for solving least squares problems with the left-hand side matrix $\bA_\mu \bOmega$. 
    This leads to an immediate opportunity for computational cost reduction and/or improvement of numerical stability, by applying iterative refinement methods for least-squares problems~\cite{carson2020three}, which compute preconditioner $\bR$ via a low-precision QR factorization of $\bA_\mu \bOmega$. 

    A particularly useful choice for $\bR$ in~\cref{eq:ls_representation2} is the RandNLA-based preconditioner from~\cite{rokhlin2008fast,AMT10,MSM14_SISC} for solving tall least-squares problems. 
    We have already employed this preconditioner to stabilize the fast Q-less Cholesky QR~(\cref{alg:prec_chol_qr}). 
    Recall that it proceeds with $\bR$ that orthogonalizes the sketched basis $\bTheta \bA_\mu \bOmega$:
    $$[\sim,\bR] \mapsfrom \mathrm{qr}(\bTheta\bA_\mu\bOmega),$$
    where $\bTheta$ is an additional small sketching matrix. 
    As discussed in~\cref{basis_free_represent}, if $\bTheta$ satisfies the $\varepsilon$-embedding property~\cref{eq:epsilon_embedding}, then $\cond(\bA_\mu \bOmega \bR^{-1}) \leq \sqrt{\frac{1+\varepsilon}{1-\varepsilon}}$. 
    This guarantees that computing the matvecs with $\bQ^{\dagger}$ via an iterative solution of~\cref{eq:ls_representation2} will require only $\mathcal{O}(\log{\frac{1}{\mathrm{u}}})$ matvecs with $\bA_\mu$ and $\bOmega$. 
    When $l$ is not too large, the cost of constructing $\bR$ is primarily determined by the cost of computing the two-sided sketch $\bTheta \bA \bOmega$. 
    This computation can be made highly efficient by using the structure of $\bOmega$, $\bTheta$, and $\bA$.

The formulation~\cref{eq:ls_representation2} also provides an effective mechanism for recycling computations when solving a sequence of similar systems. For instance, once a preconditioner $\bR$ for $\bA_\mu \bOmega$ has been precomputed for a given matrix $\bA$, it can be reused as a preconditioner for subsequent  systems, overcoming the need to recompute and orthogonalize $\bA_\mu \bOmega$ every time.

\paragraph{Factorization-free Preconditioning.}

\blue{In some situations, factoring the sketch $\bTheta\bA_\mu\bOmega$ (or $\bOmega^\mathrm{T} \bA_\mu^2\bOmega$)  may be computationally impractical. This can occur, for example, due to the cubic complexity in $l$ or due to poor compatibility of factorization routines with some modern computational architectures.}
In such cases, we can use a factorization-free approach that takes the $\bR$ factor as an identity and computes the solution $\bv = (\bA_\mu \bOmega)^\dagger \bu$ to the least-squares problem~\cref{eq:ls_representation2} via Newton-sketch method~\cite{pilanci2017newton} equivalent to the following second-level iterative refinement:
\begin{equation}\label{eq:iterative_refinement}
    \bv \gets \bv + \bD^{-1} (\bA_\mu\bOmega)^\mathrm{T} (\bu - \bA_\mu\bOmega \bv),
\end{equation}
where $\bD = (\bTheta\bA_\mu\bOmega)^\mathrm{T} (\bTheta\bA_\mu\bOmega)$.  
It can be shown~\cite{pilanci2017newton} that~\cref{eq:iterative_refinement} converges to the least-squares solution in $\mathcal{O}(\log{\frac{1}{\mathrm{u}}})$ steps if $\bTheta$ satisfies the $\varepsilon$-embedding property~\cref{eq:epsilon_embedding}.

It follows that the least-squares solution~\cref{eq:ls_representation2} can be obtained using $\mathcal{O}(\log{\frac{1}{\mathrm{u}}})$ matvecs with $\bA_\mu\bOmega$ and $(\bA_\mu\bOmega)^\mathrm{T}$, and  $\mathcal{O}(\log{\frac{1}{\mathrm{u}}})$ additional matvecs with $\bD^{-1}$. In turn, each application of $\bD^{-1}$ to a vector can be performed by an iterative solution of a $l \times l$ linear system with CG or MINRES. 
In the worst-case scenario, this would require $\mathcal{O}(\cond(\bTheta\bA_\mu\bOmega)\log{\frac{1}{\mathrm{u}}}) = \mathcal{O}(\cond(\bA_\mu \bOmega)\log{\frac{1}{\mathrm{u}}})$ matvecs with $\bD$, resulting in $\mathcal{O}(\cond(\bA_\mu \bOmega) l^2\log{\frac{1}{\mathrm{u}}})$ flops. This reduces the complexity of RandRAND solvers from cubic in the deflation dimension $d$ to ${\mathcal{O}}(d^2 \cond(\bA_\mu))$. Although, the scenario when $\cond(\bA_\mu) \ll d$ is rare, circumventing the cubic complexity in $d$ can be important from the theoretical standpoint. 
Furthermore, in practice, CG and MINRES should converge in significantly fewer iterations, and
one can apply additional regularization and preconditioning strategies to accelerate the application of $\bD^{-1}$.  

\blue{Since this RandRAND implementation involves only the computation of the two-sided sketch $\bTheta \bA \bOmega$ and avoids matrix factorizations, it can be particularly well suited to modern architectures, including GPUs and massively parallel environments.} An extensive exploration of the benefits of RandRAND preconditioning, using~\cref{eq:ls_representation,eq:ls_representation2,eq:iterative_refinement} with second-level preconditioners $\bR$ and iterative refinement will be carried out in follow-up work.

\subsection{Adaptive Selection of Sketching Dimension}

Here, we assume that the power parameter $q$ is specified \emph{a priori} to be $0$ or $1$.  
Larger values are unnecessary in most scenarios, since the quasi‑optimality coefficients scale as
$\mathcal{O}\left ( \left(\frac{n}{l}\right )^{\frac{1}{2q+2}}\right )$, which is already nearly constant for $q=0$ or~$1$.

In contrast, the sketching dimension $l$ can have a much larger impact on RandRAND performance and should be selected with care.
 In practice, it can be chosen adaptively: given $\mathbf{A}_\mu$, one can keep increasing $l$ by some factor until the subspace projection error either falls below the tolerance\footnote{For R-RandRAND and C-RandRAND the subspace projection error can be refined to $\|(\mathbf{I} - \mathbf{\Pi})\,\mathbf{A}_\mu(\mathbf{I} - \mathbf{\Pi})\|$ }
\begin{equation} \label{eq:ad_criterion}
\|(\mathbf{I} - \mathbf{\Pi})\,\mathbf{A}_\mu\|\;\le\;f \left(\frac{n}{l_{\max}}\right)^{\frac{1}{2q+2}} \lvert\mu\rvert,
\end{equation}
where $f$  is some factor (e.g., $f=100$),  or stagnates, or $l$ reaches the maximum allowed value $l_{\max}$. 
According to the results from~\cref{sec:allprecond} for spd systems with $\mu \geq 0$, the relation~\cref{eq:ad_criterion} implies that the condition number of the preconditioned matrix $\bB$ is bounded by 
$$
\cond(\mathbf{B}) \;\lesssim\;  f \left(\frac{n}{l_{\max}}\right)^{\frac{1}{2q+2}}.
$$
According to the results from~\cref{range_approx}, the tolerance~\cref{eq:ad_criterion} can be reached with high probability using $l = \mathcal{O}(k)$, where $k$ is such that $\sigma_k(\bA_\mu)$ is within a factor \blue{$\mathcal{O}(f \left(\frac{k}{l_{\max}}\right)^{\frac{1}{2q+2}} )$} of $\sigma_n(\bA_\mu)$.

When solving systems~\cref{eq:initsys} for multiple parameter  values $\mu$ via recycling methods, one can estimate 
$\|(\mathbf{I} - \mathbf{\Pi})\,\mathbf{A}_\mu\|$ with a few power iterations using the recycled sketch $\bA \bOmega$, and verify whether the approximation basis satisfies~\cref{eq:ad_criterion} or should be enlarged. 
In general, $\|(\mathbf{I} - \mathbf{\Pi})\,\mathbf{A}_\mu\|$ tends to be smaller for smaller $\lvert\mu\rvert$.

\section{RandRAND Solvers} \label{solvers}

In this section, we integrate RandRAND preconditioning into Krylov subspace methods~\cite{barrett1994templates,saad2003iterative} for solving the linear system~\cref{eq:initsys}. 
We focus on the (preconditioned) CG and MINRES algorithms for symmetric systems. 
The extensions to nonsymmetric methods such as GMRES are straightforward.

\subsection{R-RandRAND CG and MINRES}

\Cref{alg:RandRAND_CG,alg:RandRAND_CG_split} present R-RandRAND solvers based on CG and MINRES, using an explicit QR factorization of $\bA_\mu \bOmega$. 
We consider two contexts for R-RandRAND, described in~\cref{sec:precond}: right preconditioning of $\bA_\mu$; and splitting the right-hand side into orthogonal components. 
In the algorithms, the input $\bx_0$ specifies the initial guess, and $\eta$ denotes the residual tolerance. Similar algorithms can also be formulated for the left preconditioning.   

\begin{remark} \label{reorthog}
As discussed in~\cref{orthogonalization}, reorthogonalization can improve numerical stability. 
This is particularly beneficial for the split R-RandRAND,~\cref{alg:RandRAND_CG_split}. 
In this case, reorthogonalization can be incorporated by considering
$ \bb_p \mapsfrom (\bI - \bPi)(\bI - \bPi)\bb $
in step 2, and 
$ \bE \bu = (\bI - \bPi)(\bI - \bPi) \bA_\mu \bu $
in step 6.
Note that if matvecs with $\bPi$ are computationally expensive, reorthogonalization can be applied only every few iterations of CG/MINRES.
\end{remark}

\Cref{alg:BL_RandRAND_CG,} outlines a basis-less variant of R-RandRAND that avoids explicitly computing or applying the $\bQ$ factor of $\bA_\mu \bOmega$.

    \begin{algorithm} 
     \caption{R-RandRAND CG/MINRES} \label{alg:RandRAND_CG}
    \begin{algorithmic} [1]
    \Statex{\textbf{function} $\bx =\mathtt{r{\_}randrand{\_}solver}(\mu, \bA, \bb, \bOmega, \bx_0, \eta)$}
    \State Compute QR factors: $[\bQ, \bR] = \mathrm{qr}(\bA\bOmega + \mu \bOmega)$. Set $\bPi = \bQ \bQ^\mathrm{T}$ (or use more numerically robust projector from~\cref{num_stab}).
    \State Estimate $\tau = \|(\bI - \bPi) \bA_\mu (\bI - \bPi)\|$ with few power iterations. 
    \State Let $\bb \mapsfrom \bb - \bA_\mu \bx_0$.
    \State Initialize residual $\br \mapsfrom \bb$, conjugate direction for CG, or the Lanczos state for MINRES.
    \While{$\|\br\| > \eta$ \textbf{and} not reached max of iterations}        
    \State Perform a CG/MINRES iteration on preconditioned system $\bA \bP \by = \bb$, computing matvecs with $\bA \bP$ via
    $$\bA \bP \bu  := (\bI - \bPi) \bA_\mu (\bI - \bPi) \bu  + \tau \bPi \bu  .$$
  \EndWhile   
    \State Recover solution to original system $\bx \mapsfrom  (\bOmega \bR^{-1} \bQ^\mathrm{T}) ( \tau \by - \bA_\mu (\bI - \bPi)\by ) + (\bI - \bPi)\by$.
    \State (Restarting) Update tolerance $\eta$ and initial guess $\bx_0 \mapsfrom \bx$ and repeat steps 3 to 9. 
    \State Return $\bx$
    \end{algorithmic}
\end{algorithm}

    \begin{algorithm} 
     \caption{Split R-RandRAND CG/MINRES} \label{alg:RandRAND_CG_split}
    \begin{algorithmic} [1]
    \Statex{\textbf{function} $\bx =\mathtt{split_{\_}r{\_}randrand{\_}solver}(\mu, \bA, \bb, \bOmega, \bx_0, \eta)$}
    \State Compute QR factors: $[\bQ, \bR] = \mathrm{qr}(\bA\bOmega + \mu \bOmega)$. Set $\bPi = \bQ \bQ^\mathrm{T}$ (or use more numerically robust projector from~\cref{num_stab}).
    \State Let $\bb \mapsfrom \bb - \bA_\mu \bx_0$.
    \State Compute projection $\bb_p \mapsfrom (\bI-\bPi)\bb$.
    \State Initialize residual $\br \mapsfrom \bb_p$, conjugate direction for CG, or the Lanczos state for MINRES.
    \While{$\|\br\| > \eta$ \textbf{and} not reached max of iterations}        
    \State Perform a CG/MINRES iteration on preconditioned system $\bE \by = \bb_p$, computing matvecs with $\bE$ via
    $$\bE \bu  := (\bI - \bPi) \bA_\mu \bu. $$
  \EndWhile   
    \State Recover solution to original system $\bx \mapsfrom  (\bOmega \bR^{-1} \bQ^\mathrm{T}) ( \bb - \bA_\mu \by ) + \by$.
    \State (Restarting) Update tolerance $\eta$ and initial guess $\bx_0 \mapsfrom \bx$ and repeat steps 3 to 9. 
    \State Return $\bx$
    \end{algorithmic}
\end{algorithm}

\begin{algorithm} 
     \caption{Basis-less R-RandRAND CG/MINRES} \label{alg:BL_RandRAND_CG}
    \begin{algorithmic} [1]
    \Statex{\textbf{function} $\bx =\mathtt{bl{\_}r{\_}randrand{\_}solver}(\mu, \bA, \bb, \bOmega, \bx_0, \eta)$}
    \State Compute R factor via $$\bR \mapsfrom \mathtt{fast{\_}q{\_}less{\_}chol{\_}qr}(\bA+\mu\bI,\bOmega) \text{,~~~or~~~} \bR \mapsfrom \mathtt{fast{\_}q{\_}less{\_}precond{\_}chol{\_}qr}(\bA+\mu\bI,\bOmega)$$
    from~\cref{basis_free_represent}.
    \State Set up routines for computing matvecs with $\bQ$ and $\bQ^\mathrm{T}$ using expressions 
      \begin{equation} 
      \bQ \bu := \bA_\mu \bOmega \bR^{-1} \bu,~~~ \bQ^\mathrm{T} \bu := \bR^{-\mathrm{T}} \bOmega ^{\mathrm{T}}  \bA_\mu \bu. 
      \end{equation}
    \State Using routines from step 2, set up matvec operations with projector $\bPi$ via 
   $\bPi \bu := \bQ \bQ^{\mathrm{T}} \bu $. Possibly incorporate reorthogonalization or refining of the projector from~\cref{num_stab}.
     \State Perform steps 2-10 of~\cref{alg:RandRAND_CG} or~\cref{alg:RandRAND_CG_split}.
    \end{algorithmic}
\end{algorithm}

\begin{remark} \label{refined_projector}
    As discussed in~\cref{orthogonalization}, R-RandRAND allows using a refined projector $\bPi = \bQ (2\bI - \bQ^\mathrm{T}\bQ)\bQ^\mathrm{T}$ in step 3 of~\cref{alg:BL_RandRAND_CG}. 
    This choice improves numerical stability, but it incurs roughly twice the cost. 
    In practice, it can be sufficient to apply the refined projector only once at the final step when recovering the solution $\bx$ from $\by$.
\end{remark}

\subsection{C-RandRAND PCG and  MINRES}
The C-RandRAND preconditioner can be incorporated directly into PCG and preconditioned MINRES, as depicted in~\Cref{alg:RandRAND_PCG}.

 \begin{algorithm} 
     \caption{C-RandRAND PCG/MINRES} \label{alg:RandRAND_PCG}
    \begin{algorithmic} [1]
    \Statex{\textbf{function} $\bx =\mathtt{c{\_}randrand{\_}solver}(\mu, \bA, \bb, \bOmega, \bx_0, \eta)$}
    \State Compute QR factors and sketch: $[\bQ, \bR] = \mathrm{qr}(\bA\bOmega + \mu \bOmega)$ and  $\bOmega^\mathrm{T} \bA_\mu \bOmega = \bOmega^\mathrm{T}\bA\bOmega + \mu \bOmega^\mathrm{T}\bOmega$.
    \State Set projectors $\bPi = \bQ \bQ^\mathrm{T}$ and $\widetilde{\bPi} =  \bQ \bG \bQ^\mathrm{T}$, where $\bG = \bR^{-\mathrm{T}} (\bOmega^\mathrm{T} \bA_\mu \bOmega) \bR^{-1}$ (or use more numerically robust projectors from~\cref{num_stab}).
    \State Estimate $\|(\bI - \bPi) \bA_\mu (\bI - \bPi)\|$ and $\lambda_{\mathrm{max}}(\bPi \bA_\mu \widetilde{\bPi})$ with few power iterations. 
    \State Choose $\rho$ a priori (e.g., $\rho = 0.5$ or $1$), or compute it using~\cref{eq:optimal_rho}~or~\cref{eq:cnd_bound2_ext}, and set $$\tau = \rho \|(\bI - \bPi) \bA_\mu (\bI - \bPi)\| \lambda^{-1}_{\mathrm{max}}(\bPi \bA_\mu \widetilde{\bPi}).$$
    \State Set $\bb \mapsfrom \bb - \bA_\mu \bx_0$.
    \State Initialize residual $\br \mapsfrom \bb$, conjugate direction for PCG, or the Lanczos state for MINRES.
    \While{$\|\br\| > \eta$ \textbf{and} not reached max of iterations}        
    \State Perform a PCG/MINRES iteration on $\bA_\mu \bx = \bb$ using the preconditioner
    $$\bP  = \bI - \bPi  + \tau \widetilde{\bPi}.$$
  \EndWhile   
    \State (Restarting) Update tolerance $\eta$, set $\bx_0 \mapsfrom \bx$, and repeat steps 5–10. 
    \State Return $\bx$
    \end{algorithmic}
\end{algorithm}

\begin{remark}

For improved numerical stability, it is advisable in step~2 of~\cref{alg:RandRAND_PCG} to express $\bG$ as $ \bG =(\bR_{\mathrm{sk}}  \bR^{-1})^\mathrm{T} (\bR_{\mathrm{sk}}  \bR^{-1})$, where $\bR_{\mathrm{sk}}$ is a Cholesky or SVD-based factor such that $\bR_{\mathrm{sk}}^\mathrm{T} \bR_{\mathrm{sk}} = \bOmega^\mathrm{T} \bA_\mu \bOmega$. This formulation helps ensure that the preconditioner remains numerically positive definite.
\end{remark} 

\Cref{alg:BL_RandRAND_PCG} describes a basis-less variant of C-RandRAND preconditioning. 

\begin{remark}
As in the R-RandRAND setting, one can improve numerical stability of~\cref{alg:RandRAND_PCG,alg:BL_RandRAND_PCG} by reorthogonalization of the components of the preconditioner. This can be achieved by computing matvecs with $(\bI-\bPi)$ and $\widetilde{\bPi}$ as follows: $$(\bI-\bPi)\bu = (\bI-\bPi)(\bI-\bPi)\bu,~~~ \widetilde{\bPi}\bu  = \bigl(\bI - (\bI - \bPi)(\bI - \bPi)\bigr)\bA_\mu^{-1} \bPi \bu. $$ 
Furthermore, the projector $\bPi$ can be represented in a refined form: $\bPi = \bQ  (2\bI - \bQ^\mathrm{T} \bQ) \bQ^\mathrm{T}$. As discussed in~\cref{orthogonalization}, the refined projectors can mitigate rounding errors introduced during the computation of $\bR$, albeit at a higher computational cost.
\end{remark}

\begin{algorithm} 
     \caption{Basis-less C-RandRAND PCG/MINRES} \label{alg:BL_RandRAND_PCG}
    \begin{algorithmic} [1]
    \Statex{\textbf{function} $\bx =\mathtt{bl{\_}c{\_}randrand{\_}solver}(\mu, \bA, \bb, \bOmega, \bx_0, \eta)$}
    \State Compute R factor via $$\bR \mapsfrom \mathtt{fast{\_}q{\_}less{\_}chol{\_}qr}(\bA+\mu\bI,\bOmega) \text{~~~or~~~} \bR \mapsfrom \mathtt{fast{\_}q{\_}less{\_}precond{\_}chol{\_}qr}(\bA+\mu\bI,\bOmega)$$
    from~\cref{basis_free_represent}, and compute the sketch $\bOmega^\mathrm{T} \bA_\mu\bOmega = \bOmega^\mathrm{T} \bA\bOmega + \bOmega^\mathrm{T} \bOmega$.
    \State Set up routines for matvecs with $\bQ$ and $\bQ^\mathrm{T}$, using 
      \begin{equation} 
      \bQ \bu := \bA_\mu \bOmega \bR^{-1} \bu,~~~ \bQ^\mathrm{T} \bu := \bR^{-\mathrm{T}} \bOmega ^{\mathrm{T}}  \bA_\mu \bu. 
      \end{equation}
    \State Using routines from step 2, set up matvec operations with the projectors $\bPi$ and $\widetilde{\bPi} = \bPi \bA_\mu^{-1} \bPi$ as follows: 
   $$\bPi \bu := \bQ \bQ^{\mathrm{T}} \bu, \text{~~~} \widetilde{\bPi} \bu :=  \bQ \bG \bQ^\mathrm{T} \bu  $$
   where $\bG = \bR^{-\mathrm{T}} (\bOmega^\mathrm{T} \bA_\mu \bOmega) \bR^{-1}$.
    \State Perform steps 3-11 of~\cref{alg:RandRAND_PCG}.
    \end{algorithmic}
\end{algorithm}

\subsection{G-RandRAND MINRES}
Recall that the key advantage of G-RandRAND is that it provides rigorous quasi-optimality guarantees of the preconditioner, even when $\bA_\mu$ is indefinite.

G-RandRAND preconditioning can be incorporated into MINRES in nearly the same way as C-RandRAND, following the structure of~\cref{alg:RandRAND_PCG,alg:BL_RandRAND_PCG}. 
The only differences are:
\begin{itemize}
\item The operator $\bA_\mu$ can be indefinite.
\item The projector $\widetilde{\bPi}$ uses different core matrix $\bG = \left(\bR^{-\mathrm{T}} ( \bOmega^\mathrm{T} \bOmega) \bR^{-1} \right )^\frac{1}{2}. $
\item The parameter $\tau$ in step 4 is now given by
$$\tau =  \rho^{\frac{1}{2}}  {\|(\bI - \bPi) \bA_\mu\|}{\|\bA_\mu^{-1} \bPi \bA_\mu \|^{-1}} = \rho^{\frac{1}{2}}  {\|(\bI - \bPi) \bA_\mu\|}{\| \widetilde{\bPi}\bA_\mu \|^{-1}} .$$ Thus, instead of estimating $\|(\bI - \bPi) \bA_\mu (\bI - \bPi)\|$ and $\lambda_{\max}(\bPi \bA_\mu \widetilde{\bPi})$ as in C-RandRAND, one needs to estimate $\|(\bI - \bPi) \bA_\mu\|$ and $\|\widetilde{\bPi} \bA_\mu\|$.
\end{itemize}

\subsection{Convergence Guarantees for SPD Systems}
We now establish rigorous convergence guarantees for spd systems, focusing on the R-RandRAND preconditioner. Similar bounds hold for C-RandRAND and G-RandRAND, up to small constant factors.

By combining the condition number estimates from~\cref{sec:precond} with standard convergence results for CG and MINRES, we obtain the following result.

\begin{proposition} \label{thm:conv_bound}
Let $\bA_\mu$ be spd and $\mu \geq 0$. Let $\by_t$ be the approximate solution to the preconditioned R-RandRAND system $\bB \by = \bb$ obtained after $t$ iterations of CG  or MINRES (see~\cref{alg:RandRAND_CG,alg:BL_RandRAND_CG}). Let $\bx_t = \bP \by_t$. Assume that the parameter $\tau$ used in R-RandRAND satisfies $\lambda_{\mathrm{min}}(\bA_\mu) \leq \tau \leq \|\bE\|$. Then, we have 
\begin{equation} \label{eq:Berror}
    \begin{split}
    \|\by - \by_t \|_{ \bB } &\leq 2 \mathrm{e}^{-2t/T} \|\by - \by_0 \|_{\bB } \text{~~~ for CG,}\\
    \|\bA(\bx - \bx_t) \| &\leq 2 \mathrm{e}^{-2t/T} \|\bA( \bx - \bx_0) \| \text{~~~ for MINRES,}
    \end{split} 
\end{equation} 
where $T = \left(\frac{\|(\bI - \bPi)\bA_\mu(\bI - \bPi)\|}{\lambda_n+\mu} \right)^{\frac{1}{2}}$. 
%\begin{proof}
%    This result follows immediately by substituting the condition number bound~\cref{eq:cnd_bound} into the standard convergence guarantees of MINRES and CG.
%\end{proof}
\end{proposition}

The coefficient $T$ in~\cref{thm:conv_bound} represents the number of iterations needed to reduce the R-RandRAND error by a factor of $\mathrm{e}^{2} \approx 7.4$. It follows that the number of iterations in~\cref{alg:RandRAND_CG} needed to compute the solution to $\bA_\mu \bx = \bb$ with relative residual error $\mathrm{u}$ is bounded by 
\begin{equation} 
    t \leq \tfrac{1}{2} \, T \left( \log{2} + \log{\tfrac{1}{\mathrm{u}}} \right). 
\end{equation}

According to~\cref{thm:apriori2}, the coefficient $T$ in~\cref{thm:conv_bound} should be small when the sketching dimension $l$ is comparable to the effective spectral dimension. This guarantee is formalized below.
\begin{proposition} \label{thm:conv_bound2}
Consider R-RandRAND preconditioning with a Gaussian sketching matrix $\bX$ of size $l \geq 2.4d+4$, where $d$ is such that $\lambda_d \leq \lambda_n + \mu$.   Then the constant $T$ in~\cref{thm:conv_bound} satisfies
$$\mathbb{E}(T) \leq 1 +   \left (148\tfrac{n}{d} \right)^{\frac{1}{4q+4}}.$$
In addition, with probability $1-6^{-d}$,
\begin{equation} \label{eq:T_bound}
T \leq 1 +   \left (3240\tfrac{n}{d} \right)^{\frac{1}{4q+4}}.
\end{equation}
\end{proposition}

Furthermore,~\Cref{thm:apriori2} allows to extend the result from~\cref{thm:conv_bound2} to scenarios with varying ratios between $\lambda_d$ and $\mu$, as well as different decays of the trailing spectrum. 
For instance, we can guarantee that with high probability,  $T \leq \mathcal{O}(1)$ when $d$ satisfies $\lambda_d \leq \mathcal{O}\left(\left(\frac{n}{d}\right)^{-\frac{1}{2q}} \mu\right)$ for $q \geq 1$. 
Table 4.3 shows the bounds for $\mathbb{E}(T)$ under selected ratios of $n$ to $d$, ratios of $\lambda_d$ to $\mu$, and different powers $q$.

 \begin{table}[ht]
\centering
\small
\setlength{\tabcolsep}{3pt}
\renewcommand{\arraystretch}{1.0}

\begin{subtable}{.46\textwidth}
\centering
\begin{tabular}{lcccc}
\toprule
 & $\tfrac{n}{d} \leq 10^2$ & $\tfrac{n}{d} \leq 10^3$ & $\tfrac{n}{d} \leq 10^4$ & $\tfrac{n}{d} \leq 10^5$ \\
\midrule
$\tfrac{\lambda_d}{\mu} \leq 10$   & 20.5 & 36.4 & 64.7 & 115 \\
$\tfrac{\lambda_d}{\mu} \leq 10^0$ & 8.8  & 15.6 & 27.6 & 50  \\
$\tfrac{\lambda_d}{\mu} \leq 10^{-1}$ & 6.5 & 11.5 & 20.5 & 36.4 \\
\bottomrule
\end{tabular}
\caption{$q=0$}
\end{subtable}
\hspace{0.04\textwidth}
\begin{subtable}{.46\textwidth}
\centering
\begin{tabular}{lcccc}
\toprule
 & $\tfrac{n}{d} \leq 10^2$ & $\tfrac{n}{d} \leq 10^3$ & $\tfrac{n}{d} \leq 10^4$ & $\tfrac{n}{d} \leq 10^5$ \\
\midrule
$\tfrac{\lambda_d}{\mu} \leq 10$   & 6.1 & 8.1 & 10.7 & 14.3 \\
$\tfrac{\lambda_d}{\mu} \leq 10^0$ & 2.3 & 3.0 & 4.0  & 5.3  \\
$\tfrac{\lambda_d}{\mu} \leq 10^{-1}$ & 1.2 & 1.5 & 2.0  & 2.6  \\
\bottomrule
\end{tabular}
\caption{$q=1$}
\end{subtable}

\caption{Upper bounds on $\mathbb{E}(\frac{1}{2}T)$, the expected number of iterations required to reduce the error by a factor of $\mathrm{e}$ when using CG or MINRES with RandRAND preconditioning based on Gaussian sketching with sketching dimension $l \geq 2.4d + 4$. Scenarios vary with $n/d$, $\lambda_d/\mu$, and $q \in \{0,1\}$.}
\label{table_iter}
\end{table}

Since applying a preconditioned R-RandRAND operator using an explicit basis $\bQ$ costs $T_{\mathrm{mv}} + \mathcal{O}(nl)$ flops, \cref{thm:conv_bound2} implies that solving an spd  system with CG or MINRES (excluding the preconditioner construction) requires
$$
\mathcal{O}\bigl((T_{\mathrm{mv}} + nl)\,T\,\log\tfrac{1}{\mathrm{u}}\bigr)
\;=\;
\mathcal{O}\left(\left(T_{\mathrm{mv}} + nd\right)\left(\tfrac{n}{d}\right)^{\frac{1}{4q+4}}
\log\tfrac{1}{\mathrm{u}}\right)
$$
flops.
We can establish similar guarantees for SRHT and sparse random matrices. According to \cref{thm:apriori3,thm:apriori4}, the bound for $T$ holds with probability $1-\delta$ if $\bX$ is an SRHT matrix with $l = \mathcal{O}(d \log(\tfrac{n}{\delta}))$ rows, or a sparse random matrix with $l = \mathcal{O}(d \log(\tfrac{n}{\delta}))$ rows and $\gamma = \mathcal{O}(\log(\frac{k}{\delta}))$ nonzeros per column. 

In the basis-less approach, applying a preconditioned R-RandRAND operator costs $(q+1)T_{\mathrm{mv}} + \mathcal{O}(n \log n)$ flops; hence solution of spd system with CG or MINRES (excluding preconditioner construction) requires
$$
\mathcal{O}\bigl((q\,T_{\mathrm{mv}} + n \log n)\,T\,\log\tfrac{1}{\mathrm{u}} \bigr)
\;=\;
\mathcal{O}\left(\bigl(qT_{\mathrm{mv}} + n \log n\bigr)\left(\tfrac{n}{d}\right)^{\frac{1}{4q+4}}
\log\tfrac{1}{\mathrm{u}}\right)
$$
time.
    \section{Empirical Evaluation}
\label{emp_eval}
In this section, we study the performance and numerical properties of RandRAND and compare them with randomized preconditioners based on Nystr{\"o}m low-rank approximation~\cite{frangella2023randomized,derezinski2024faster} across various problem families. In particular, 
in \Cref{validation_randrand_properties}, we evaluate RandRAND on synthetic systems;
in \Cref{sxn:basis-explicit-randrand} we apply RandRAND to a variety of linear systems arising in machine learning; and in \Cref{sec:basis-less-pde-appn}, we consider the performance of basis-less RandRAND preconditioners on an inverse PDE problem and kernel systems under memory constraints.  For the normal equations arising from tall least-squares problems, we also tested sketching-based low-rank preconditioners~\cite{lacotte2020effective}, but they performed substantially worse because they fail to preserve the smallest eigenvalue of the operator unless the sketching dimension is comparable to the effective trace dimension, which is computationally impractical in our setting.

Overall, RandRAND often outperforms the low-rank preconditioning in terms of condition numbers (and especially their bounds), iteration counts, and numerical stability in the basis-less setting.

A comprehensive evaluation of the computational benefits of RandRAND will be carried out in future work.

\subsubsection*{Experimental setup}

In all test cases, for R-RandRAND preconditioners, we select parameter $\tau$ using~\cref{eq:cnd_bound} with $\rho=1$. For C-RandRAND preconditioners, we select parameter $\tau$  based on minimization of the bound~\cref{eq:cnd_bound2_ext} from~\cref{appendix_d} either using the exact value of $\lambda_{\min}(\bA_\mu)$ or the estimate $\lambda_{\min}(\bA_\mu) \approx \mu$. In some experiments we also evaluate the performance of C-RandRAND using $\tau = \|\bPi \bA_\mu^{-1}\bPi \|^{-1} $.   For G-RandRAND preconditioners, we select parameter $\tau$ using~\cref{eq:g_randrand_tau} with  $\rho=0.5$.

A key implementation detail in all our experiments is incorporating restarting of the preconditioned iterative solver after each reduction of the residual norm by a factor of $\eta = 10^{-2}$, which helps to ensure numerical stability. 
We tested RandRAND and other randomized preconditioning approaches for both CG and MINRES. 
The overall trends were similar in both cases. 
However, because the MINRES convergence curves were more regular and representative, we chose to report those results.

\subsection{Validation of RandRAND properties} \label{validation_randrand_properties}

First, we validate the properties of RandRAND on synthetic systems, where we can control the condition number and the spectrum.

\paragraph{Systems with exponential+polynomial spectrum decay.}  

\begin{figure}
 \begin{subfigure}[b]{0.3\textwidth}
      \centering
      \includegraphics[width=\textwidth]{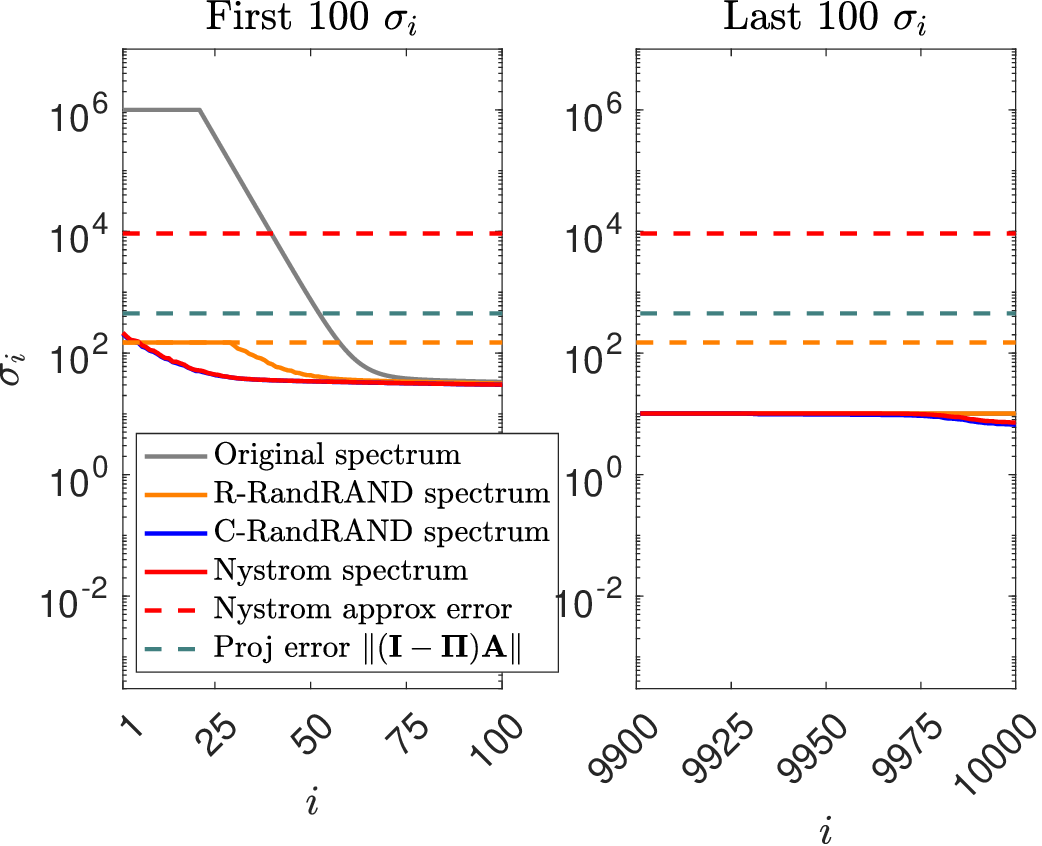}
      \subcaption{
      $\bA$ with tail $\sigma_i \sim i^{-1/4}$; $l = 100$}
      \label{fig:cond-num-sketch-size2a}
  \end{subfigure} 
   \begin{subfigure}[b]{0.3\textwidth}
      \centering
      \includegraphics[width=\textwidth]{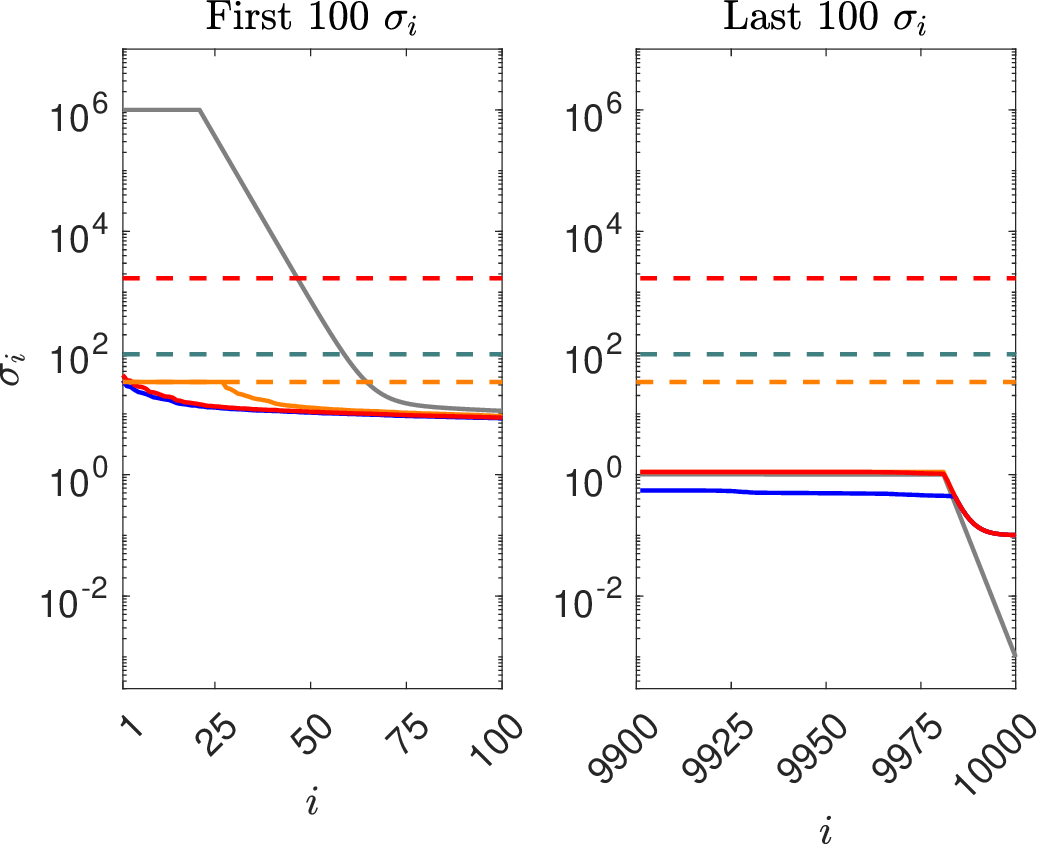}
      \subcaption{$\bA$ with tail $\sigma_i \sim i^{-1/2}$; $l = 100$}
      \label{fig:cond-num-sketch-size2b}
  \end{subfigure} 
  \begin{subfigure}[b]{0.3\textwidth}
      \centering
      \includegraphics[width=\textwidth]{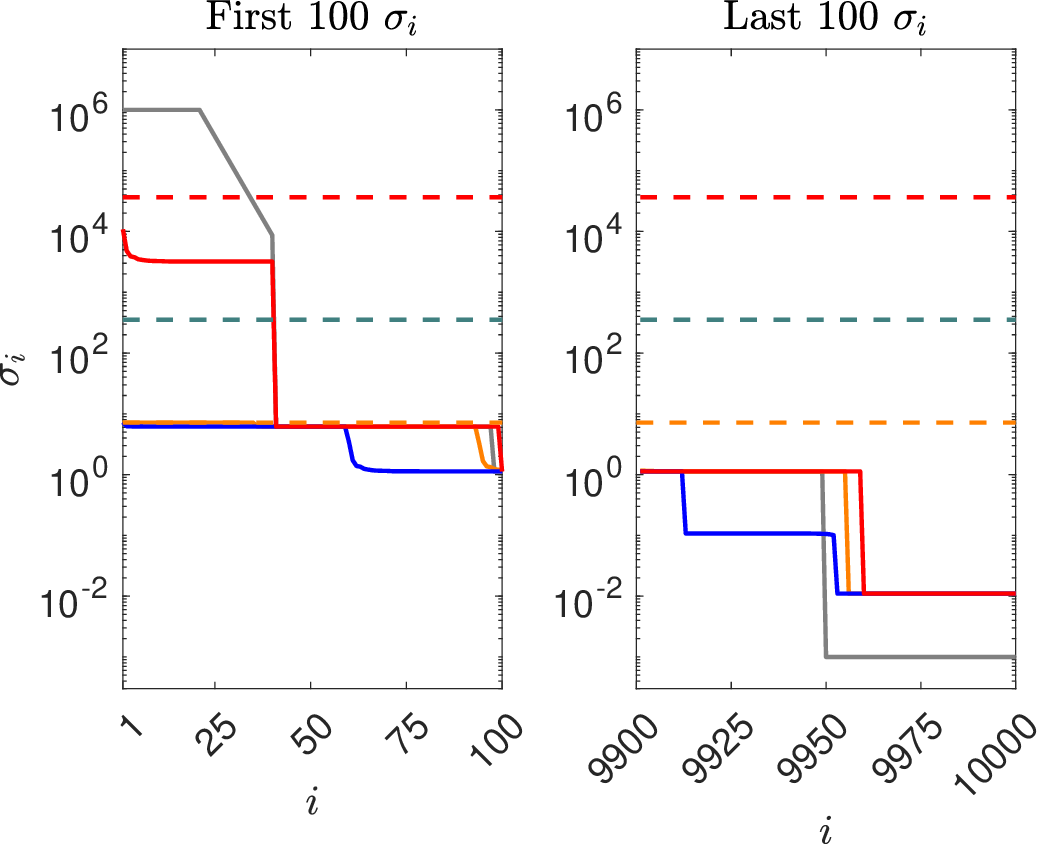}
      \subcaption{$\bA$ with stepwise $\sigma_i$; $l = 40$}
      \label{fig:cond-num-sketch-size2_step}
  \end{subfigure} 
  \caption{Eigenvalues of $\bA$ with exponential+polynomial or stepwise decay of the spectrum, and the eigenvalues of preconditioned $\bA_\mu$ with RandRAND or Nystr{\"o}m preconditioning using random Gaussian $\bOmega$ of size $l$.  }
  \label{fig:cond-num-sketch-size2}
\end{figure}

We fix the system size $n = 10^4$ and the regularization parameter $\mu = 0.1$. 
Then we generate two psd matrices $\bA$ of the form $\mathbf{A} = \mathbf{U} \mathbf{\Sigma} \mathbf{U}^\mathrm{T}$, where $\mathbf{U}$ is an orthogonalized Gaussian matrix and $\mathbf{\Sigma}$ is diagonal. 
The matrices $\mathbf{\Sigma}$ have few large entries, on the order of $10^6$, followed by an exponentially decaying gap and a power-law tail of the form $i^{-\alpha}$ with $\alpha = 1/4$ in the first scenario and $\alpha = 1/2$ in the second. 
We also add a spike of decreasing smallest singular values in the second scenario to test whether the methods preserve the smallest eigenvalues $\lambda_{n}(\mathbf{A})+\mu$. 
The largest and smallest $100$ eigenvalues of each tested $\mathbf{A}$ are shown in~\cref{fig:cond-num-sketch-size2a,fig:cond-num-sketch-size2b}. 

We then construct randomized preconditioners for $\mathbf{A}_\mu$ using Gaussian sketching matrices $\mathbf{\Omega}$ of varying sizes. 
We also build the randomized Nystr{\"o}m preconditioner following~\cite{frangella2023randomized} and use it as a baseline for evaluating effectiveness of RandRAND.  

\begin{figure}[htbp]

  \centering
  \begin{subfigure}[b]{0.3\textwidth}
      \centering
      \includegraphics[width=\textwidth]{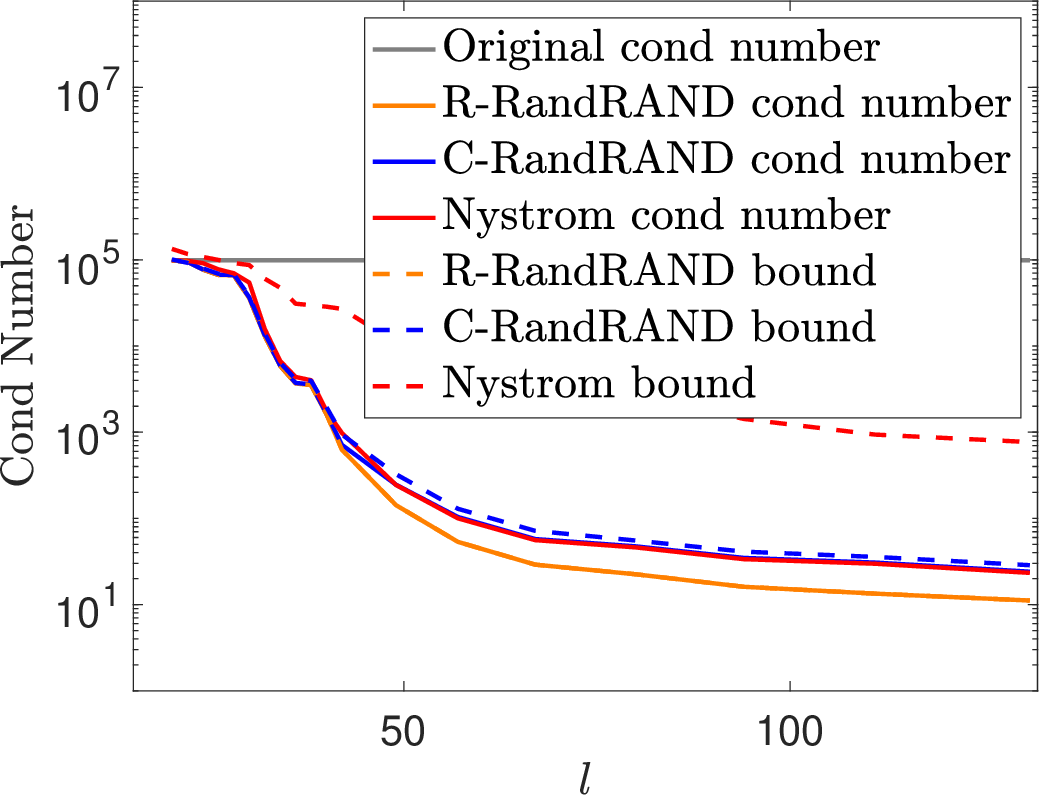}
       \subcaption{$\bA$ with tail $\sigma_i \sim i^{-\frac{1}{4}}$}
       \label{fig:cond_num_vs_sketch_size_1}
  \end{subfigure}
  \quad \quad \quad \quad 
  \begin{subfigure}[b]{0.3\textwidth}
      \centering
      \includegraphics[width=\textwidth]{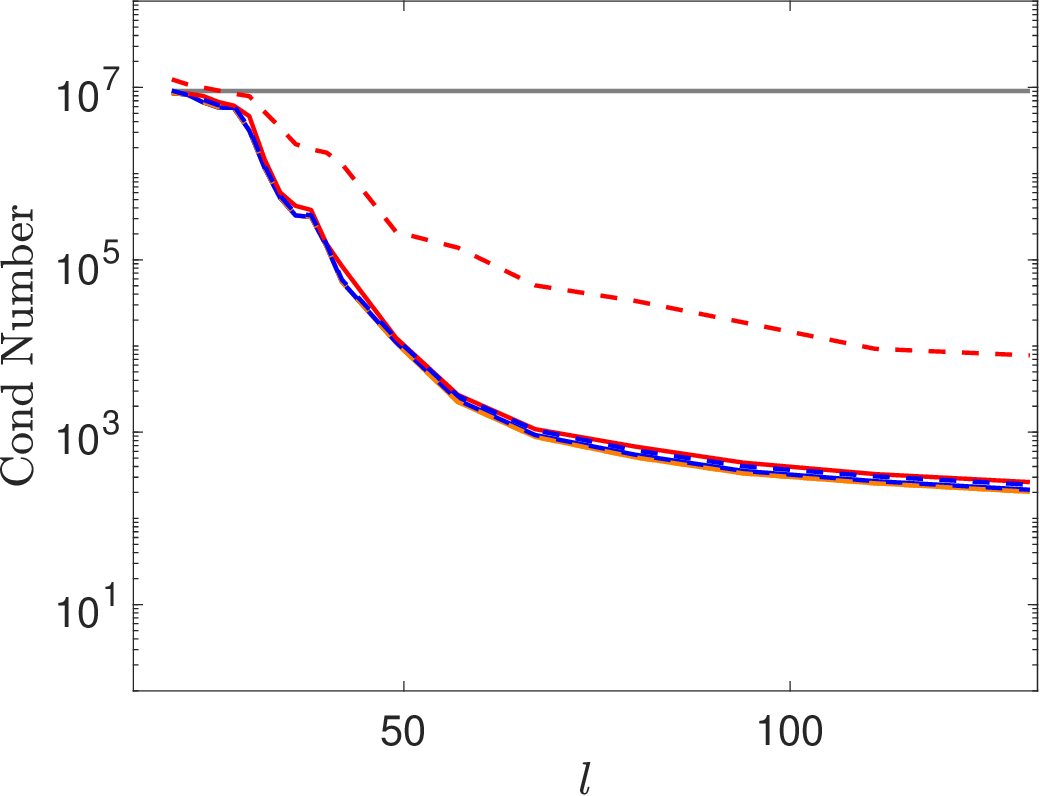}
      \subcaption{$\bA$ with tail $\sigma_i \sim i^{-\frac{1}{2}}$}
      \label{fig:cond_num_vs_sketch_size_2}
  \end{subfigure}
  \caption{Condition numbers $\mathrm{cond}(\bA_\mu)$ and $\mathrm{cond}(\bA_\mu\bP)$ for R-RandRAND and $\mathrm{cond}(\bP^{1/2}\bA_\mu\bP^{1/2})$, for  C-RandRAND and Nystr{\"o}m preconditioning and their bounds vs size of $\bOmega$. 
  }
  \label{fig:cond_num_vs_sketch_size}
\end{figure}

From~\cref{fig:cond_num_vs_sketch_size_1,fig:cond_num_vs_sketch_size_2}, we see that all considered randomized preconditioners yield condition number improvement by two orders of magnitude already when using $\mathbf{\Omega}$ of size $l = 40$ and by four orders of magnitude (for the second matrix and near this value for the first) when using $\mathbf{\Omega}$ of size $l=100$. 
Moreover, as the sketch size $l$ increases, the condition numbers of the preconditioned systems improve consistently. 

We see that, in addition to RandRAND's more efficient and flexible design, it generally achieves lower condition numbers than the randomized Nystr{\"o}m preconditioning from~\cite{frangella2023randomized}, given the same sketching dimension~$l$.\footnote{unless $\bA$ is a rank-$l$ matrix within spectral error $\mu$}
Specifically, as seen from~\cref{fig:ratio_vs_Nystrom} R-RandRAND yields condition numbers lower by a factor of up to $1.5$ for relatively small $l$ and up to $2.2$ for relatively large $l$. 
C-RandRAND yields condition numbers lower by up to a factor of $1.5$ for relatively small $l$ and $1.25$ for relatively large $l$. 
Moreover, over $100$ random trials with different samples of $\bOmega$ and $\bA$, we observed that the condition number gap between RandRAND and Nystr{\"o}m preconditioning can reach up to $3.3$. 
We also compare the empirical results with theoretical condition number bounds: RandRAND bounds (from~\cref{appendix_d}) closely match observed values, while the Nystr{\"o}m bound (from~\cite[Proposition 5.3]{frangella2023randomized}) is overly pessimistic, and in most cases is off by more than an order of magnitude. 
This validates (i) RandRAND's stronger theoretical guarantees and (ii) the possibility of more effective adaptive selection criterion for the sketch size for RandRAND.

\begin{table}[h!]
\centering
\small
\setlength{\tabcolsep}{2pt}
\renewcommand{\arraystretch}{0.9}
\begin{tabular}{l *{5}{c}*{5}{c}}
\toprule
 & \multicolumn{5}{c}{Ratio of C-RandRAND and Nystr{\"o}m} & \multicolumn{5}{c}{Ratio of R-RandRAND and Nystr{\"o}m} \\
\cmidrule(lr{0.4em}){2-6}\cmidrule(lr{0.4em}){7-11}
& $l=30$ & $l=34$ & $l=40$ & $l=80$ & $l=111$
& $l=30$ & $l=34$ & $l=40$ & $l=80$ & $l=111$ \\
\cmidrule(lr{0.4em}){2-6}\cmidrule(lr{0.4em}){7-11}
$i^{-\tfrac{1}{4}}$, sampled:
& 1.49 (2.3) & 1.13 (7.4) & 1.07 (14) & 0.97 (50) & 0.97 (26) 
& 1.51 (2.4) & 1.13 (8.0) & 1.09 (17) & 2.05 (122) & 2.22 (70) \\
\cmidrule(lr{0.4em}){2-6}\cmidrule(lr{0.4em}){7-11}
$i^{-\tfrac{1}{2}}$, sampled: 
& 1.47 (2.3) & 1.16 (6.2) & 1.06 (11.5) & 1.24 (54) & 1.20 (30) 
& 1.45 (2.5) & 1.15 (6.5) & 1.06 (12) & 1.33 (65) & 1.27 (36) \\
\cmidrule(lr{0.4em}){2-6}\cmidrule(lr{0.4em}){7-11}
$i^{-\tfrac{1}{2}}$, selected: 
& 1.14 (2.1) & 1.60 (5.7) & 3.14 (11.6) & 1.34 (55.5) & 1.28 (45) 
& 1.13 (2.2) & 1.60 (6.0) & 3.15 (12.2) & 1.46 (67.6) & 1.35 (54.7) \\
\bottomrule
\end{tabular}
\caption{Ratios between condition numbers (or their bounds) of operators with different power-law singular value tails, preconditioned by the  Nystr{\"o}m method and by RandRAND, for varying sketch dimensions $l$. The test matrix $\bOmega$ is a Gaussian matrix either randomly sampled or selected from $100$ random samples with maximum gain of RandRAND over Nystr{\"o}m preconditioning. }
\label{fig:ratio_vs_Nystrom}
\end{table}

For further insight, we took the preconditioners from the previous experiment associated with test matrices  of sizes $l=40$ and $l=100$ and evaluated the extreme eigenvalues  of preconditioned operators; results are depicted in~\cref{fig:cond-num-sketch-size2a,fig:cond-num-sketch-size2b}. 
We observe that the preconditioners provide a quasi-optimal deflation of the largest eigenvalues of $\mathbf{A}_\mu$ while maintaining the smallest eigenvalues.  
We additionally report the approximation errors $\|(\mathbf{I} - \mathbf{\Pi})\,\mathbf{A}_\mu\|$, $\|\mathbf{E}\| = \|(\mathbf{I} - \mathbf{\Pi})\,\mathbf{A}_\mu\,(\mathbf{I} - \mathbf{\Pi})\|$, and $\|\mathbf{A} - \hat{\mathbf{A}}\|$, where $\hat{\mathbf{A}}$ is the Nystr{\"o}m approximation of $\mathbf{A}_\mu$. According to~\cref{sec:allprecond}, $\|(\mathbf{I} - \mathbf{\Pi})\,\mathbf{A}_\mu\|$ and $\|\mathbf{E}\|$ provide bounds on the largest eigenvalue of the RandRAND-preconditioned matrix. By~\cite[Proposition 5.3]{frangella2023randomized}, $\|\mathbf{A} - \hat{\mathbf{A}}\|+\mu$ is a bound on the largest eigenvalue of the Nystr{\"o}m-preconditioned matrix. 
This experiment again confirms that the RandRAND bounds are robust, while the Nystr{\"o}m bound can be highly pessimistic.

\paragraph{System with step-wise spectrum decay.}  
To realize the benefits of RandRAND in terms of theoretical guarantees and robustness more starkly, we constructed $\bA = \bU\bSigma \bU^\mathrm{T}$ as in previous example, but this time we let the entries of $\mathbf{\Sigma}$ to have a stepwise decay shown in~\cref{fig:cond-num-sketch-size2_step}. 
We set regularization $\mu = 0.01$. 
Then we generated Gaussian $\bOmega$ with sketching dimension $l=40$ (this dimension aligns with the drop of eigenvalues of $\bSigma$).
As seen from~\cref{fig:cond-num-sketch-size2_step}, in this setting, Nystr{\"o}m’s condition number can be worse than RandRAND’s by more than three orders of magnitude. Moreover, Nystr{\"o}m’s largest eigenvalue can exceed $\|(\mathbf{I} - \mathbf{\Pi})\,\mathbf{A}_\mu\|$ by a factor of $30$. All this demonstrates that for Nystr{\"o}m preconditioning, the projection-based arguments alone can be insufficient, while for RandRAND, the projection-based arguments are enough to guarantee quasi-optimality and robustness of deflation.

\begin{figure}[htbp]
  \centering
  \begin{subfigure}[b]{0.3\textwidth}
      \centering
      \includegraphics[width=\textwidth]{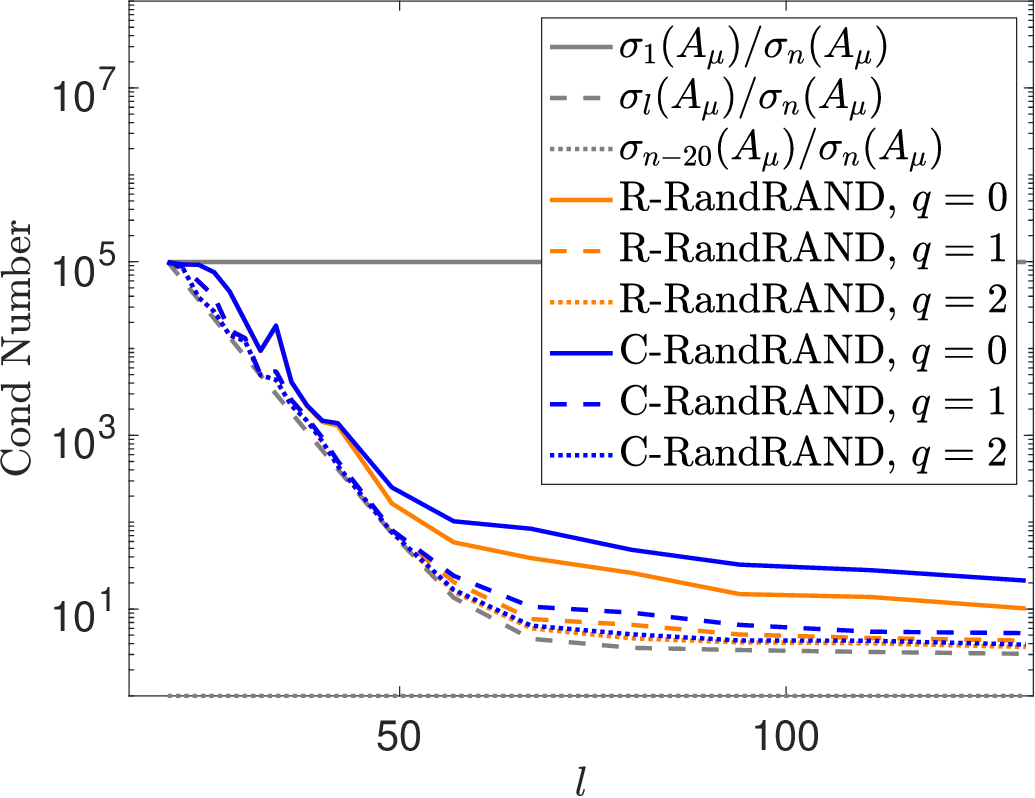}
       \subcaption{$\bA$ with tail $\sigma_i \sim i^{-1/4}$}
  \end{subfigure}
  \hfill
  \begin{subfigure}[b]{0.3\textwidth}
      \centering
      \includegraphics[width=\textwidth]{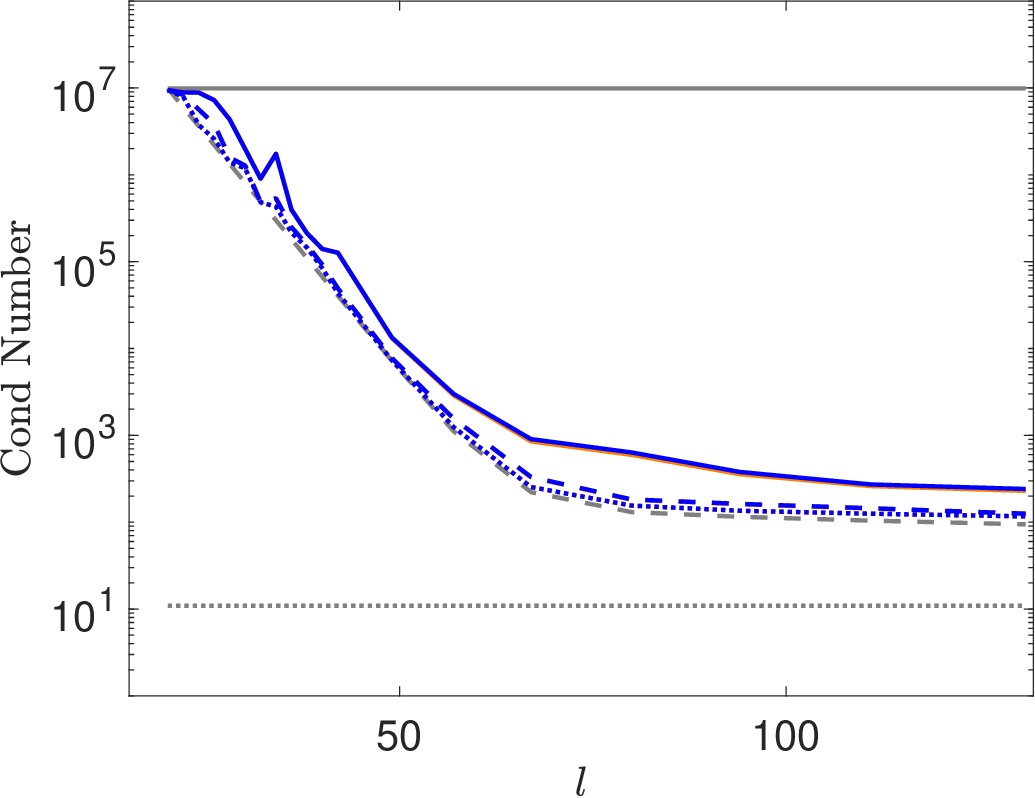}
      \subcaption{$\bA$ with tail $\sigma_i \sim i^{-1/2}$}
  \end{subfigure}
  \hfill
  \begin{subfigure}[b]{0.3\textwidth}
      \centering
      \includegraphics[width=\textwidth]{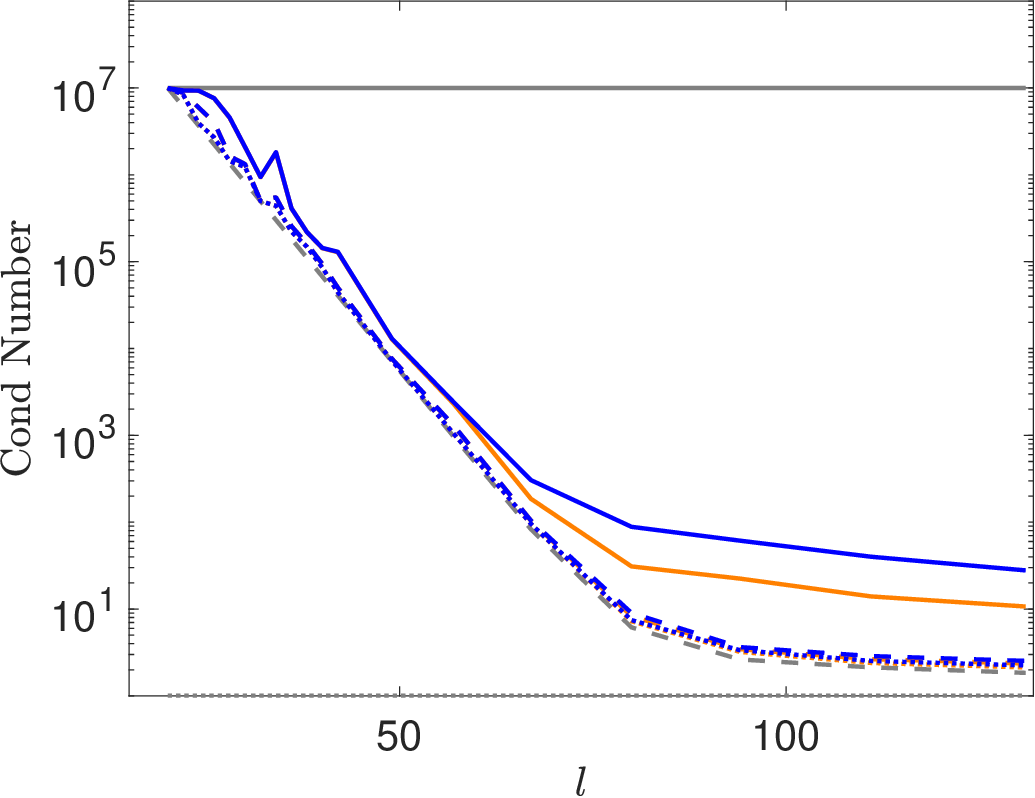}
      \subcaption{$\bA$ with tail $\sigma_i \sim i^{-3/2}$}
  \end{subfigure}
  \caption{Condition numbers $\mathrm{cond}(\bA_\mu\bP)$ for R-RandRAND and $\mathrm{cond}(\bP^{1/2}\bA_\mu\bP^{1/2})$ for  C-RandRAND vs size of $\bOmega$. The test matrix $\bOmega$ is taken as $\bA^q \bX^\mathrm{T}$, where $\bX$ is a Gaussian OSE.}
 \label{fig:cond_num_vs_sketch_size_q}
\end{figure}

\paragraph{RandRAND with subspace iteration.}
To verify the advantages of subspace iteration, we consider test matrices of the form $\bOmega = \bA^q \bX^\mathrm{T}$, and we vary the sketch size $l$ and the exponent $q$. 
The matrix $\bX$ is taken as a Gaussian OSE.
\Cref{fig:cond_num_vs_sketch_size_q} reports the condition numbers of the resulting preconditioned matrices $\bA_\mu$ with exponential+polynomial spectral decay from the previous example. 
We also test a third matrix $\bA_\mu$ that shares the exponential component of the second matrix but exhibits a faster decay of the trailing eigenvalues (of the form $i^{-3/2}$) and lacks the added spike of small eigenvalues. 
As a baseline, we use the optimal condition number that can be achieved by deflating an $l$-dimensional subspace, given by
$
\min_{\mathrm{rank}(\bPi^*) = l-1} \|(\bI - \bPi^*)\bA_\mu\|/\sigma_n(\bA_\mu) = \sigma_{l}(\bA_\mu)/\sigma_n(\bA_\mu)
$.

From~\cref{fig:cond_num_vs_sketch_size_q}, 
we observe that subspace iteration becomes highly useful once the sketch size reaches the stall region of the eigenvalue decay, where simply increasing $l$ yields only marginal gains. 
For example, setting $q = 1$ can improve the condition number by up to a factor of $10$ for C-RandRAND and $5$ for R-RandRAND compared to $q = 0$.
Moreover, using $q = 1$ or $2$ achieves near-optimal range deflation for the given $l$, and we do not see any practical benefit from using larger $q > 2$. This observation aligns well with our theoretical analysis from~\cref{range_approx,sec:allprecond}.

\paragraph{Indefinite Systems.} 

We turn to RandRAND in the context of preconditioning shifted indefinite systems. We reuse the matrices from the previous experiment, but flip the signs of the 21st to 25th eigenvalues and apply a negative shift, $\mu = -100$. 
We test G-RandRAND with either sketch size $l$ and power $q=1$, or sketch size $2l$ and power $q=0$. Both configurations require a similar number of matvecs with $\bA$ for construction. 
For comparison, we also consider a preconditioner defined as the square root of the C-RandRAND preconditioner for $\bA_\mu^2$, as described in~\cref{rmk:g_randrand_vs_c_randrand}. 
The sketch dimension for this preconditioner is chosen as $l$ for $q=0$ and $l/2$ for $q=1$, to ensure comparable construction cost to the G-RandRAND variants. 

\begin{figure}[htbp]
  \centering
  \begin{subfigure}[b]{0.3\textwidth}
      \centering
      \includegraphics[width=\textwidth]{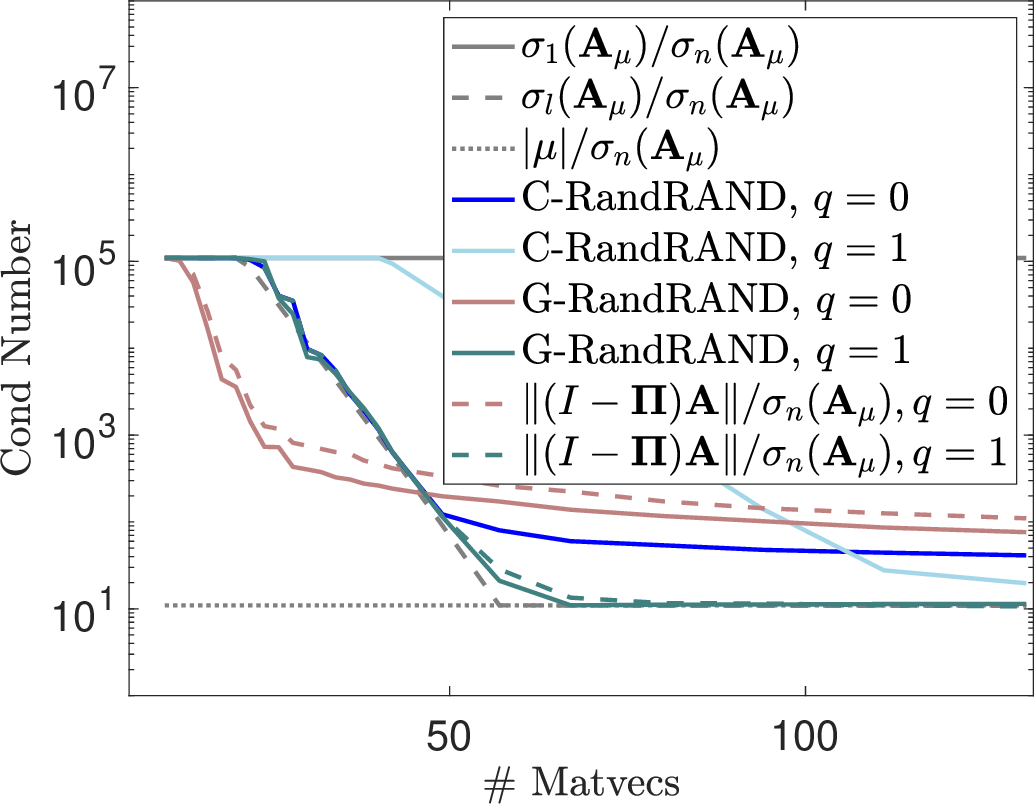}
      \subcaption{$\bA$ with tail $\sigma_i \sim i^{-1/4}$}
  \end{subfigure}
  \hfill
  \begin{subfigure}[b]{0.3\textwidth}
      \centering
      \includegraphics[width=\textwidth]{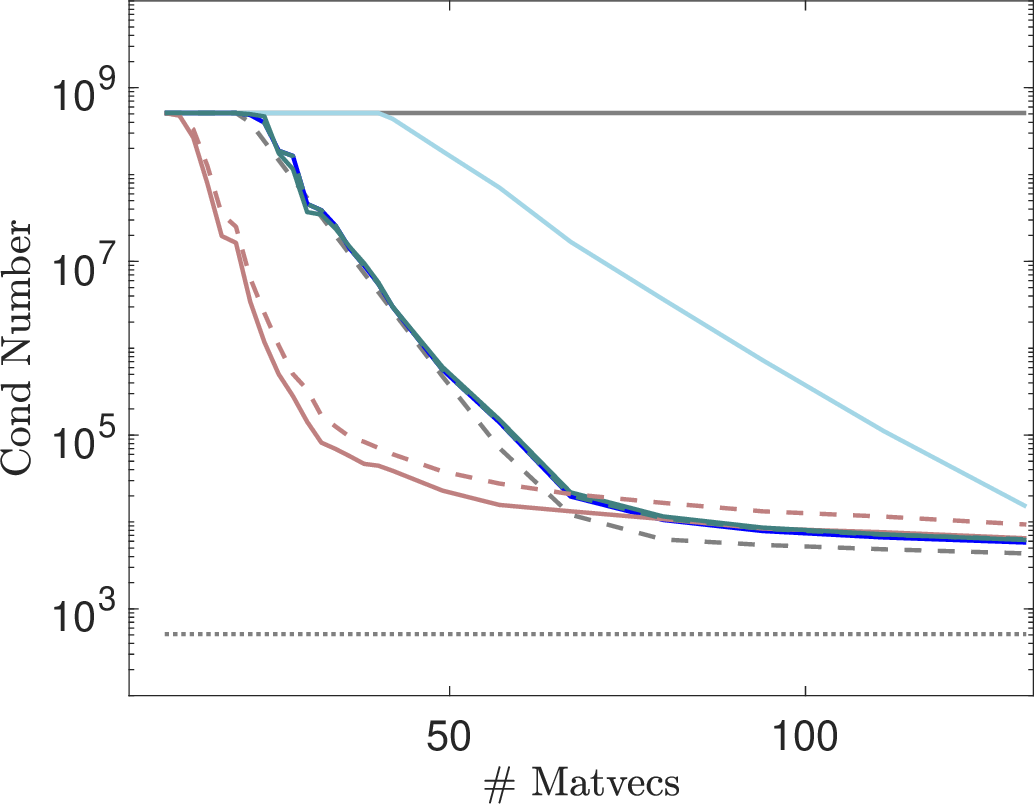}
      \subcaption{$\bA$ with tail $\sigma_i \sim i^{-1/2}$}
  \end{subfigure}
  \hfill
  \begin{subfigure}[b]{0.3\textwidth}
      \centering
      \includegraphics[width=\textwidth]{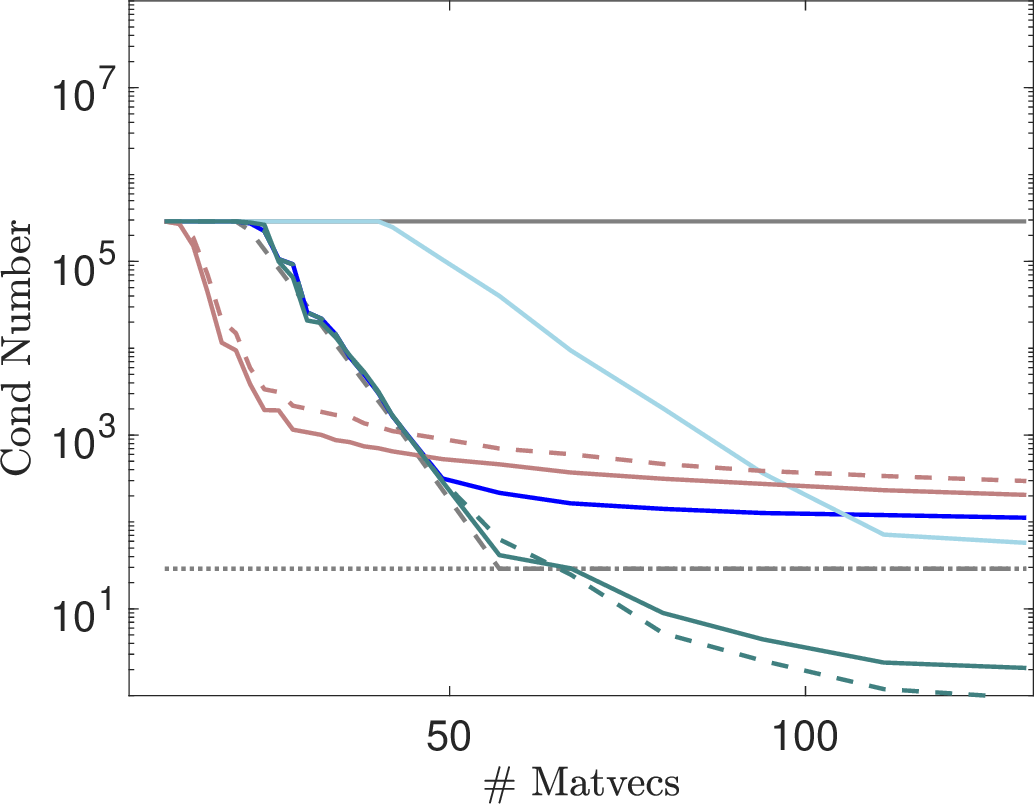}
      \subcaption{$\bA$ with tail $\sigma_i \sim i^{-3/2}$}
  \end{subfigure}
  \caption{
    Condition numbers $\mathrm{cond}(\bP^{1/2}\bA_\mu\bP^{1/2})$, where $\bP$ is a RandRAND preconditioner. The test matrix $\bOmega$ is taken as $\bA^q \bX^\mathrm{T}$, where $\bX$ is a Gaussian OSE. We test $\bP$ constructed with G-RandRAND using sketch dimension $l$ and $q=1$, G-RandRAND using sketch dimension $2l$ and $q=0$, and the square root of the C-RandRAND preconditioner for $\bA_\mu^2$ using sketch dimension $l$ and $q=0$, or sketch dimension $l/2$ and $q=1$.
  }
  \label{fig:cond_num_vs_sketch_size_q_indefinite}
\end{figure}

\Cref{fig:cond_num_vs_sketch_size_q_indefinite} shows the condition numbers of the deflated systems versus $\#$ of matvecs with $\bA$ needed for constructing the preconditioners. 
We observe that G-RandRAND can be highly effective, reducing the spectrum nearly to $\sigma_l(\bA_\mu)/\sigma_n(\bA_\mu)$, or even surpassing this value.  
As in the spd case, when $l$ lies in a range of rapid singular value decay, it is more advantageous to use a larger sketch size with $q=0$. However, when the spectral decay stalls, using $q=1$ with a smaller sketch size becomes more pertinent. 
Importantly, when $q=1$ and the sketch size becomes large enough so that $\|(\bI - \bPi) \bA\| < |\mu|$, G-RandRAND deflates not only the largest singular values but also the smallest ones down to $|\mu|$ (see~\cref{rmk:g_rand_rand_smallest})
and thus may surpass the condition number value $\sigma_l(\bA_\mu)/\sigma_n(\bA_\mu)$. 
Compared to the C-RandRAND preconditioner, G-RandRAND (with either $q=0$ or $q=1$) can achieve significantly smaller condition numbers for the same construction cost. 
For relatively small $l$, the difference with G-RandRAND using $q=0$ can reach several orders of magnitude, while for larger $l$, the difference with G-RandRAND using $q=1$ can reach up to $1.5$ orders of magnitude.

\subsection{Applications of Basis-explicit RandRAND}
\label{sxn:basis-explicit-randrand}

Next, we apply RandRAND to a variety of linear systems arising in machine learning, including ridge regression, kernel ridge regression, and portfolio optimization.

\subsubsection{Ridge regression}
Ridge regression is a method built upon ordinary least squares~\cite{hoerl1970ridge}. 
Mathematically, it consists in a least squares problem with Tikhonov regularization:
\begin{align*}
\min_{\bx \in\mathbb{R}^n}\frac{1}{2m}\|\bZ \bx -\bv\|^2 + \frac{\mu}{2}\|\bx\|^2
\end{align*}
equivalent to linear system~\cref{eq:initsys} with $\bA = \frac{1}{m}\bZ^\mathrm{T} \bZ$ and $\bb = \frac{1}{\sqrt{m}} \bZ^\mathrm{T} \bv$. 
Here, $\bZ \in \mathbb{R}^{m \times n}$ is the data matrix, and $\bv \in \mathbb{R}^n$ is a vector of labels. In our experiments we construct $\bv$ with a one-hot vector encoding. We set $v_i = 1$ if the $i$-th label is equal to the $1$-st, and $0$ otherwise. We solve the ridge regression problem on standard machine learning datasets from~\cite{vanschoren2014openml,chang2011libsvm}, summarized in~\cref{tab:rr_dataset}. We use experimental setup similar to the one used in~\cite{frangella2023randomized}. For the shuttle, yearMSD, and higgs datasets, we employ random Fourier features with kernel parameters $\gamma$ specified in~\cref{tab:rr_dataset}. For smallNORB, we use random ReLU features.
We solve the linear systems to a relative residual $10^{-8}$ with a max iteration count set of $5000$ as a hard upper cap. 

\begin{table}[H]
   \centering
   \begin{tabular}{@{}llrrrr@{}} \toprule
    Dataset & Features & $m$ & $n$ & $\gamma$  & $\mu$ \\ \midrule
    shuttle & Fourier & 43500 & 10000 & 8/9, 4/9, 4/3, 13/3 & 1e-8$/m$ \\
    yearMSD &  Fourier & 463715 & 15000 & 1/128, 1/20, 1/8 & 1e-5 \\
    higgs   & Fourier & 1000000& 10000 & 0.01, 0.02, 0.04 & 1e-4 \\
    smallNORB & ReLU & 464809& 10000 & -  & 6e-4 \\
    CIFAR10 & None & 50000 & 3072  & - & 1e-5 \\ \bottomrule
    %guillermo & None & 16000 & 4296 & - & 1e-4 \\ \bottomrule
   \end{tabular}
   \caption{Description of ridge regression datasets.}
   \label{tab:rr_dataset}
\end{table}

\begin{table}[H]
\small
\centering
\setlength{\tabcolsep}{3pt}
\renewcommand{\arraystretch}{1.0}
\begin{tabular}{@{}l ccc crrr @{\hspace{15pt}} rrr @{\hspace{15pt}} rrr @{\hspace{15pt}} rrr}
\toprule
Dataset, $\gamma$
& \multicolumn{3}{c}{Sketch size $l$}
& \shortstack{\# MINRES}
& \multicolumn{3}{c}{\hspace{-7.5pt} \# iter Nystr\"om \hspace{-15pt}}
& \multicolumn{3}{c}{\# iter C-RAND* \hspace{-15pt}}
& \multicolumn{3}{c}{\# iter C-RAND \hspace{-15pt}}
& \multicolumn{3}{c}{\# iter R-RAND} \\
\cmidrule(lr){2-4} \cmidrule(lr){6-8} \cmidrule(lr){9-11} \cmidrule(lr){12-14} \cmidrule(lr){15-17}
& $l_1$ & $l_2$ & $l_3$ & 
& $l_1$ & $l_2$ & $l_3$
& $l_1$ & $l_2$ & $l_3$
& $l_1$ & $l_2$ & $l_3$
& $l_1$ & $l_2$ & $l_3$ \\
\midrule
shuttle, $4/9$   & 400 & 200 & 100 & >5e3 & 15  & 21  & 80  & 28  & 40  & 88  & 22  & 34  & 93  & \textbf{11} & 15  & 61  \\
shuttle, $8/9$   & 400 & 200 & 100 & >5e3 & 18  & 47  & 579  & 33  & 64  & 549  & 30  & 68 & 500  & \textbf{12} & 39  & 446 \\
shuttle, $4/3$   & 400 & 200 & 100 & >5e3 & 26  & 185 & 1803 & 44    & 184 & 1784 & 42  & 203 & 1629 & \textbf{17} & 157  & 1540  \\
shuttle, $13/3$  & 800 & 400 & 200 & >5e3 & 109 & 1190 & >5e3 & 114 & 1183 & >5e3 & 114 & 1162& >5e3 & \textbf{91 }&  1132& >5e3 \\
yearMSD, $1/128$ & 400 & 200 & 100 & >5e3 & 212 & 305 & 627 & 195 & 313 & 655 & 189  & 296 & 672 & \textbf{175} & 279 & 572 \\
yearMSD, $1/20$  & 400 & 200 & 100 & 3696 & 243 & 341 & 604  & 243 & 340 & 608  & \textbf{224} & 334 & 580  & 231 & 326 & 558 \\
yearMSD, $1/8$   & 400 & 200 & 100 & 691 & 88 & 123 & 179  & 88  & 123 & 179  & 89 & 117 & 174  & \textbf{81} & 113 & 169 \\
higgs, $0.01$    & 400 & 200 & 100 & 928 & 43  & 78   & 82  & 66  & 55  & 61  & 52  & 70  & 78  & \textbf{23} & 43  & 60 \\
higgs, $0.02$    & 400 & 200 & 100 & 2643 & 99 & 167 & 234  & 118  & 155 &  223  & 96 & 158 & 222  & \textbf{77} & 158 & 214 \\
higgs, $0.04$    & 800 & 400 & 200 & >5e3 & 80  & 305 & 599 & 71 & 310  & 576 & 85 & 319  & 596  & \textbf{65} & 282 & 578 \\
smallNORB        & 800 & 400 & 200 & 543 & 31 & 53 & 92  & 43  &  57   & 87  & 39  & 63   & 87  & \textbf{28}  & 53   & 86 \\
CIFAR10          & 1000& 500 & 300 & >5e3 & 240 & 509 & 758 & 240 & 516 & 758 & {253} & 496 & 734  & \textbf{229}  & 486  & 745 \\
\bottomrule
\end{tabular}
\caption{Numbers of (preconditioned) MINRES iterations to achieve relative residual $10^{-8}$ with sketch dimensions $l_1,l_2,l_3$. C-RAND* uses $\tau=\|\bPi\bA^{-1}\bPi\|$, whereas C-RAND uses $\tau$ chosen by minimizing~\cref{eq:cnd_bound2_ext} taking $F=\mu^{-1}$.}
\label{tab:rr_results}
\end{table}

We evaluate RandRAND preconditioning with varying sketch dimensions $l$. This setting can be considered to validate both RandRAND with fixed and adaptively chosen $l$ (for different $l_{\max}$ and $f$ in the criterion~\cref{eq:ad_criterion}). We set $q=0$, as $q=1$ doubles the preconditioner construction cost without noticeably improving convergence for the tested problems. For comparison, we also report results for Nystr{\"o}m preconditioning from~\cite{frangella2023randomized}, using the same sketch sizes as RandRAND. We also tested Newton-Sketch preconditioning~\cite{lacotte2020effective,pilanci2017newton}, but it performed much worse due to its failure to preserve the smallest eigenvalue for the selected values of $l$, so we omit those results.

\begin{figure}[H]
   \begin{subfigure}{.253\textwidth}
     \centering
     \includegraphics[width=\linewidth]{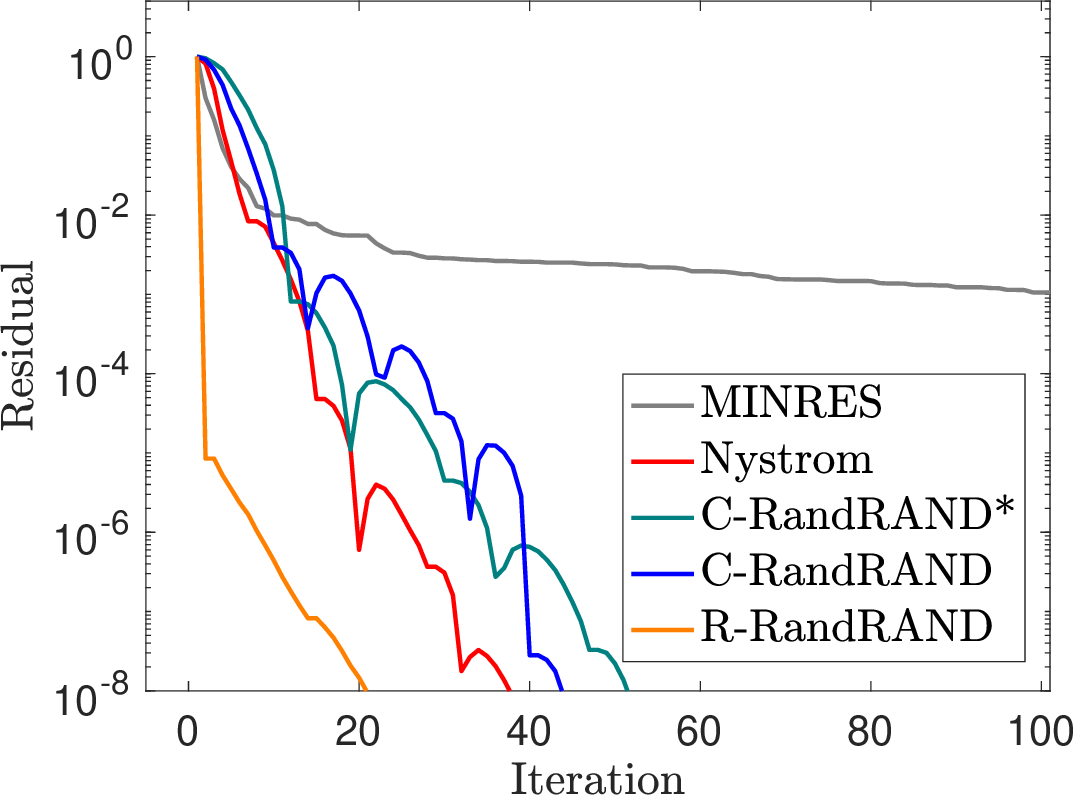}
     \caption{shuttle(8/9) $l=250$}
   \end{subfigure}%
   \hspace{0.0025\textwidth}
   \begin{subfigure}{.24\textwidth}
     \centering
     \includegraphics[width=\linewidth]{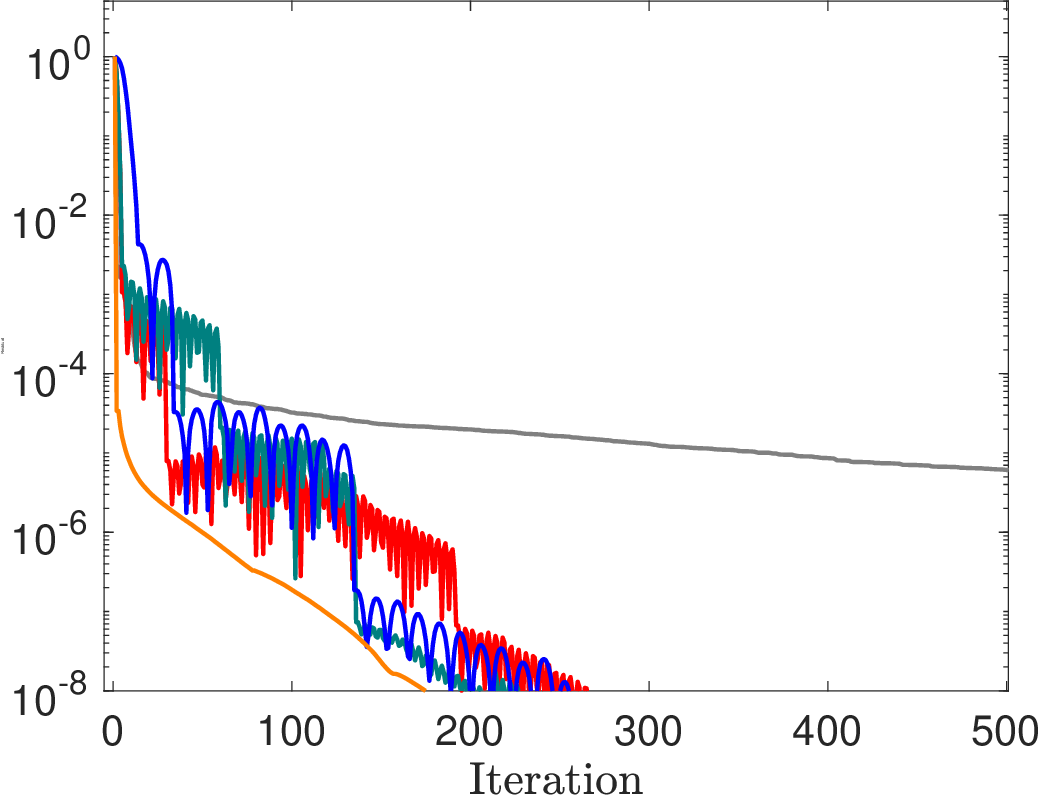}
     \caption{yearMSD(1/128) $l = 400 $}
   \end{subfigure}
   \hspace{0.0025\textwidth}
   \begin{subfigure}{.24\textwidth}
     \centering
     \includegraphics[width=\linewidth]{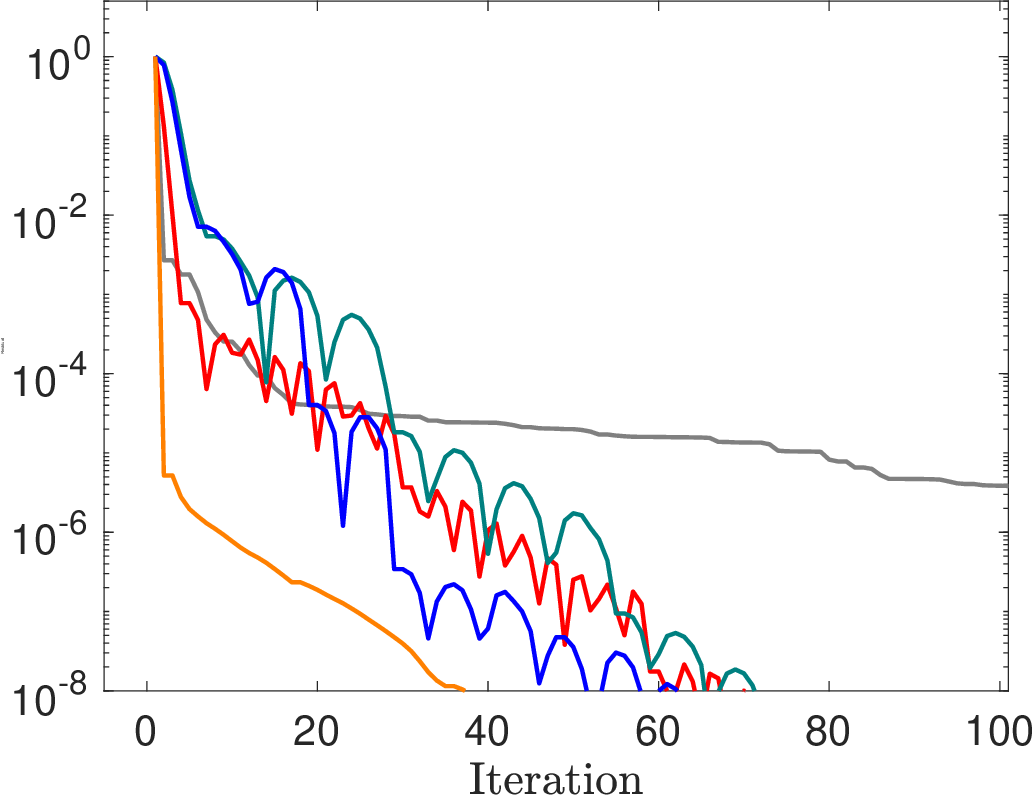}
     \caption{higgs(0.01) $l=250$}
   \end{subfigure}%
   \hspace{0.0025\textwidth}
   \begin{subfigure}{.24\textwidth}
     \centering
     \includegraphics[width=\linewidth]{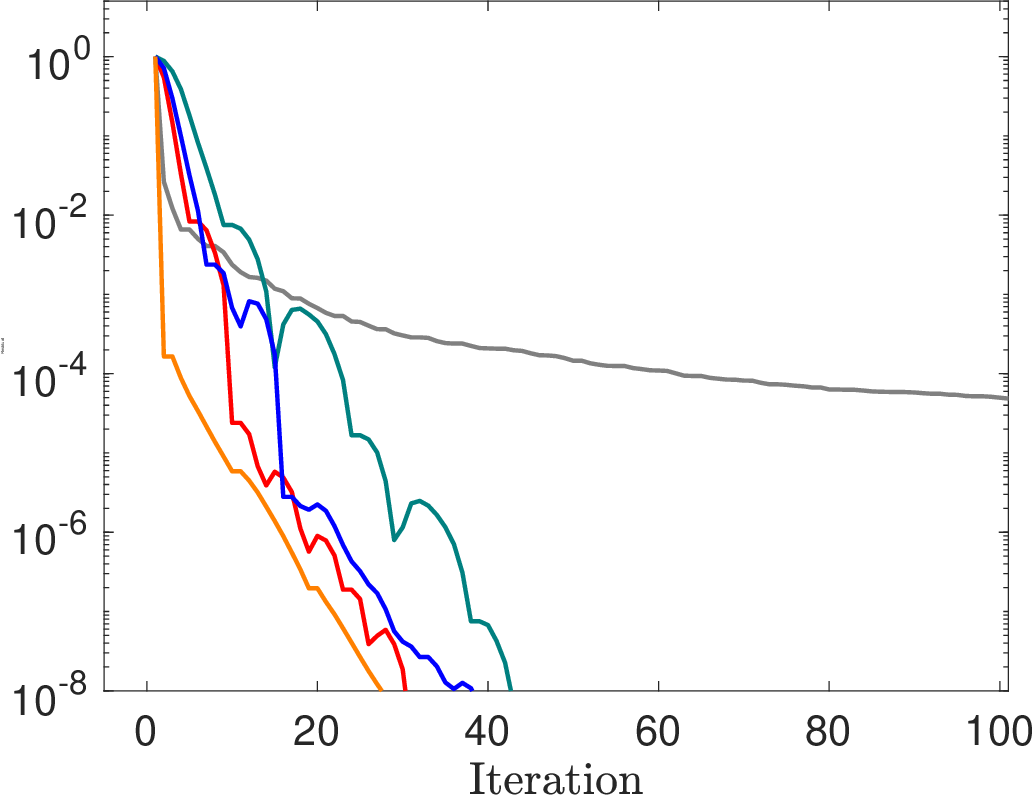}
     \caption{smallNORB $l=800$}
   \end{subfigure}%
   \caption{MINRES convergence for four different ridge regression datasets and sketch sizes.}
   \label{fig:rr_convergence_progress}
\end{figure}

RandRAND preconditioners drastically reduce the number of MINRES iterations compared to the unpreconditioned case. In almost all instances, R-RandRAND requires fewer iterations than other preconditioners. In some test cases the difference between R-RandRAND, C-RandRAND and Nystr\"om  is significant, and in others all three preconditioners perform similarly. For example, on the higgs dataset with $\gamma=0.01$, R-RandRAND with $l=400$ requires almost half as many MINRES iterations as the Nystr{\"o}m preconditioner. Figure~\ref{fig:rr_convergence_progress} illustrates residual convergence for four representative datasets and sketch sizes.
On the shuttle and higgs datasets, R-RandRAND outperforms the alternatives in the early iterations. It exhibits smoother residual decay that reflects greater robustness. This robustness is particularly evident on the yearMSD dataset: for high accuracy solution, R-RandRAND and Nystr{\"o}m require a similar number of iterations, but for moderate accuracy solution with relative residual $10^{-7}$, R-RandRAND needs about $120$ iterations, whereas Nystr{\"o}m requires more than $180$.

\subsubsection{Kernel Ridge Regression}\label{krr}

We now study the performance of our preconditioned MINRES on a closely related regression problem involving nonlinear kernels.
In kernel ridge regression, the kernel method maps the data into a nonlinear feature space and then performs linear regression in that lifted space~\cite{scholkopf2002learning}.
Let $\mathcal H$ be a reproducing kernel Hilbert space of functions on $X \subset \mathbb{R}^d$ with norm $|\cdot|$. The optimization problem reads
\begin{equation}\label{eq:kernel_reg}
\min_{f \in \mathcal H} \frac{1}{2n} \sum_{i=1}^n (f(u_i)-v_i)^2 + \frac{\mu}{2}|f|^2,
\end{equation}
where ${(u_i,v_i) \in \mathbb{R}^{d+1}}$ are data–label pairs and $\kappa : X \times X \to \mathbb{R}$ is the reproducing kernel of $\mathcal H$.
By the representer theorem~\cite{kimeldorf1970correspondence}, the minimizer has the form
$f^*(\cdot) = \sum_{i=1}^n x_i \kappa(u_i,\cdot)$,
where the coefficients $\bx = [x_1,\hdots, x_n]^\mathrm{T}$ satisfy
\begin{equation} \label{eq:kernel}
\left (\frac{1}{n}\bK + \mu \bI \right )\bx = \frac{1}{n} \by,
\end{equation}
with kernel matrix $\bK \in \mathbb{R}^{n \times n}$ given by $\bK_{ij} = \kappa(u_i,u_j)$ and $\by = [\kappa(u_1,v_1), \kappa(u_2,v_2)\dots,\kappa(u_n,v_n)]^\mathrm{T}$.
We focus on solving the full system~\eqref{eq:kernel} to sufficient precision with an iterative linear solver and preconditioning, rather than using approximate methods such as classical sketch-and-solve approaches that apply an inverse shifted low-rank approximation to the right-hand side~\cite{williams2000using} and Falkon~\cite{rudi2017falkon}. While such methods can be cheaper since they avoid matvecs with the full $\bK$, solving the exact system can yield lower test-set error~\cite{frangella2023randomized}.
For experiments, we adopt datasets considered in~\cite{gittens2016revisiting,frangella2023randomized} for evaluating kernel methods based on Nystr\"{o}m approximation. We employ the rbf kernel with bandwidth parameter $\gamma$. Dataset details are given in~\cref{tab:krr_dataset}.
The label vectors $[v_1,\dots,v_n]$ in~\eqref{eq:kernel_reg} are obtained using one-vs-all encoding as in~\cite{frangella2023randomized}: the $i$-th label vector has entries $v_j  = 1$ if the $j$-th label matches the $i$-th class, and $v_j = -1$ otherwise. This may result in a system with multiple right-hand sides. For wine, abalone and ijcnn1 datasets we take $\bb$ to be the first right-hand-side vector.  For MNIST we solve the system for multiple right-hand-sides with a block preconditioned MINRES\footnote{To ensure numerical stability, in block MINRES we use Krylov basis reorthogonalization}. All systems are solved to relative residual tolerance~\footnote{for MNIST defined as $\max_i \|\bA \bx_i - \bb_i\|/\|\bb_i\|$} specified in~\cref{tab:krr_dataset}.

\begin{table}[H]

   \centering
   \begin{tabular}{@{}llrrrr@{}} \toprule
    Dataset &  $n$ & $\gamma$ & $\mu$ & tol &right-hand-sides \\ \midrule
    wine    &  5197 & 2.1 & 1e-3/$n$, 1e-4/$n$ & 1e-6& 1  \\
    abalone &  3341 &  1  & 1e-3/$n$, 1e-4/$n$ & 1e-6& 1 \\
    MNIST  & 60000  &  1e-2    & 1e-7, 1e-8 & 1e-4 & 10 \\
    ijcnn1 & 49990  & 1    & 1e-6, 1e-7  &  1e-6 & 1 \\ \bottomrule
   \end{tabular}
   \caption{Description of kernel ridge regression datasets.}
   \label{tab:krr_dataset}
\end{table}

\begin{table}[H]
\centering
\setlength{\tabcolsep}{3pt} % tighter columns
\renewcommand{\arraystretch}{1.0}
\begin{tabular}{@{}l ccc @{\hspace{15pt}} rrr @{\hspace{15pt}} rrr @{\hspace{15pt}} rrr @{\hspace{15pt}} rrr@{}}
\toprule
Dataset, $\mu$ 
& \multicolumn{3}{c}{Sketch sizes $l$}
& \multicolumn{3}{c}{{\hspace{-7.5pt}} Nystr\"om {\hspace{-15pt}}}
& \multicolumn{3}{c}{C-RAND*{\hspace{-15pt}}}
& \multicolumn{3}{c}{C-RAND{\hspace{-15pt}}}
& \multicolumn{3}{c}{R-RAND} \\
\cmidrule(lr){2-4}\cmidrule(lr){5-7}\cmidrule(lr){8-10}\cmidrule(lr){11-13}\cmidrule(lr){14-16}
& $l_1$ & $l_2$ & $l_3$
& $l_1$ & $l_2$ & $l_3$
& $l_1$ & $l_2$ & $l_3$
& $l_1$ & $l_2$ & $l_3$
& $l_1$ & $l_2$ & $l_3$ \\
\midrule
wine, 1e-3/$n$     & 2000 & 1000 & 500 & 356 & 459 & 548 & 386 & 454 & 542 & 255  & 294 & 347 & \textbf{234} & 301 & 358 \\
wine, 1e-4/$n$     & 2000 & 1000 & 500 & 925 & 1189 & 1467 & 1046 & 1182  & 1464 & 478 & 541 & 652  & \textbf{435} & 551   & 666 \\
abalone, 1e-3/$n$  & 1000 & 500  & 250 & 197 & 259 & 248 & 224 & 295 & 393 & 147 & 190  & 254 & \textbf{136} & 173 & 236\\
abalone, 1e-4/$n$  & 1000 & 500  & 250 & 502 & 633 & 686 & 455 &  691 & 725 & 338 & 519 & 685 & \textbf{293} & 438 & 638  \\
MNIST, 1e-7      & 2000 & 1000 & 500 & 84 & 104 & 135 & 81 & 107 & 130 & 82 & 112 & 135 & \textbf{71} & 100 & 129 \\
MNIST, 1e-8      & 2000 & 1000 & 500 & 172 & 211 & 247 & 168 & 211 & 250 & 174 & 210 & 243 & \textbf{155} & 204 & 248 \\
ijcnn1, 1e-6     & 2000 & 1000 & 500 & \textbf{32} & 53 & 92 & 57 & 82 & 123 & 54 & 77 & 108 &  {39} & 58 & 88 \\
ijcnn1, 1e-7     & 2000 & 1000 & 500 & \textbf{80} & 155  & 255 & 123 & 213 & 302 & 117 & 180 & 276 & {85} & 146 & 243 \\
\bottomrule
\end{tabular}
\caption{Number of preconditioned MINRES iterations to achieve the target residual tolerance for kernel system~\cref{eq:kernel} using RandRAND and Nystr\"om preconditioning with sketching dimensions $l=l_1,l_2,l_3$. C-RAND* uses $\tau=\|\bPi\bA^{-1}\bPi\|$, whereas C-RAND uses $\tau$ chosen by minimizing~\cref{eq:cnd_bound2_ext} taking $F = \mu^{-1}$.}
\label{tab:krr_results}
\end{table}

\begin{figure}[H]
   \begin{subfigure}{.253\textwidth}
     \centering
     \includegraphics[width=\linewidth]{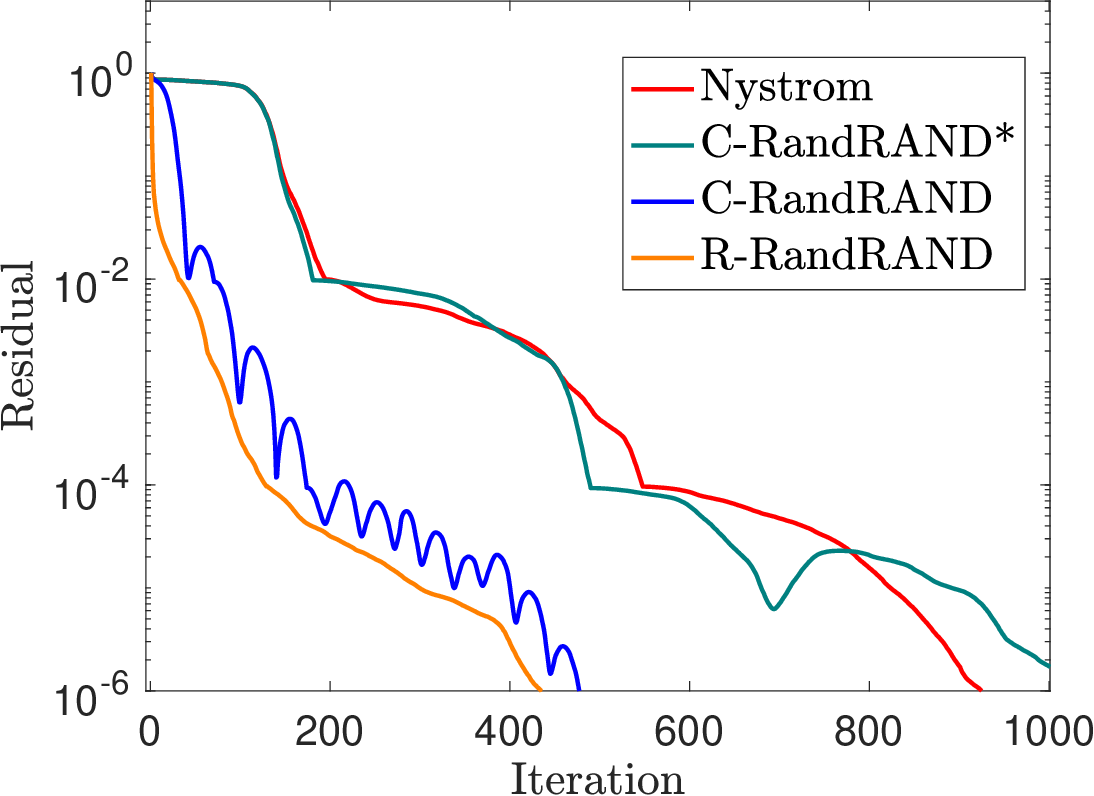}
     \caption{wine(1e-4$/n$) $l=2000$}
   \end{subfigure}%
   \hspace{.0025\textwidth}
   \begin{subfigure}{.24\textwidth}
     \centering
     \includegraphics[width=\linewidth]{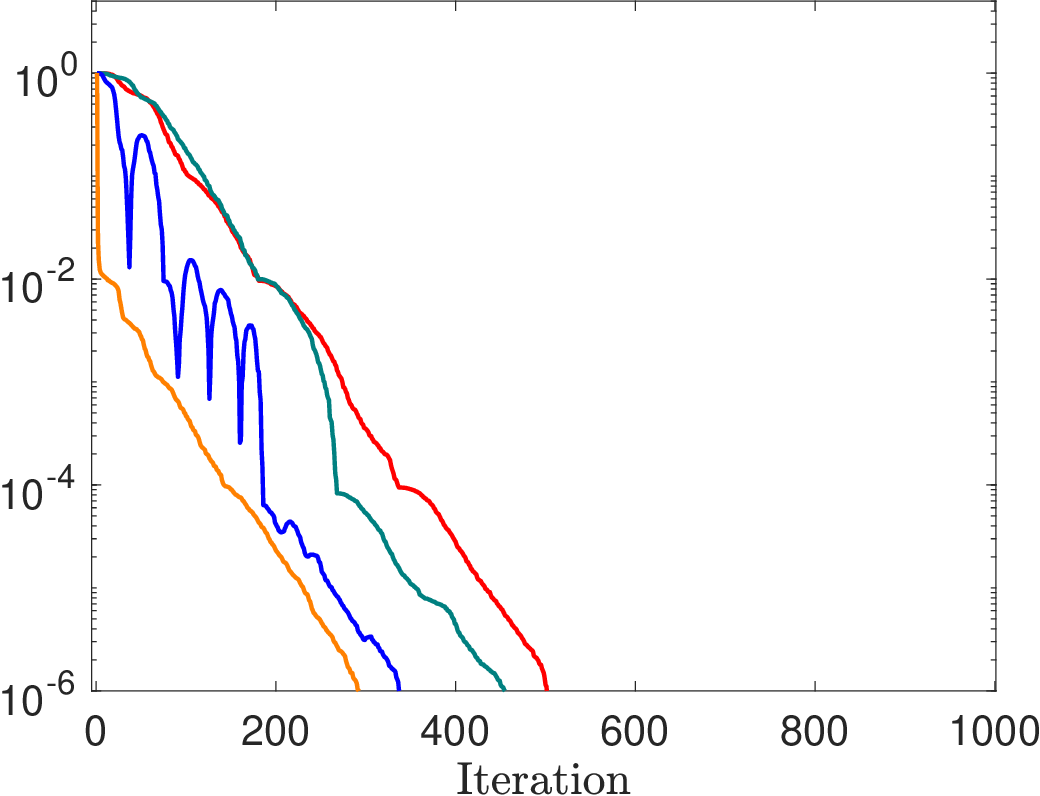}
     \caption{abalone(1e-4$/n$) $l = 1000 $}
   \end{subfigure}
   \hspace{.0025\textwidth}
   \begin{subfigure}{.24\textwidth}
     \centering
     \includegraphics[width=\linewidth]{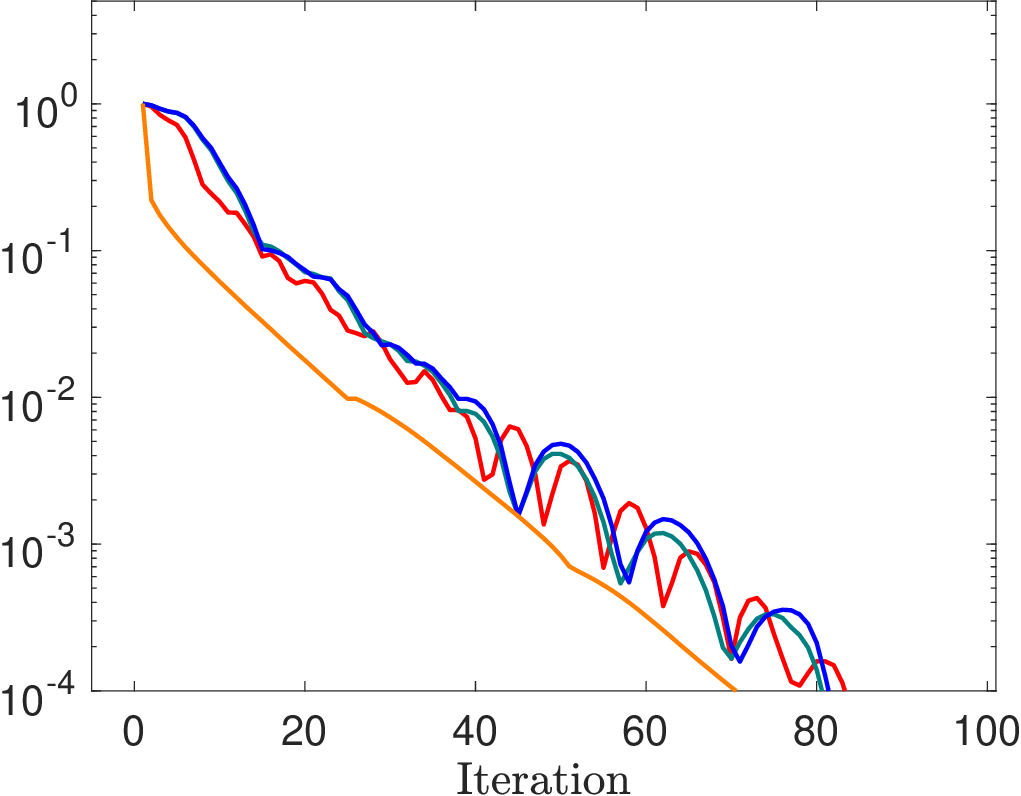}
     \caption{MNIST(1e-7) $l=2000$}
   \end{subfigure}%
   \hspace{.0025\textwidth}
   \begin{subfigure}{.24\textwidth}
     \centering
     \includegraphics[width=\linewidth]{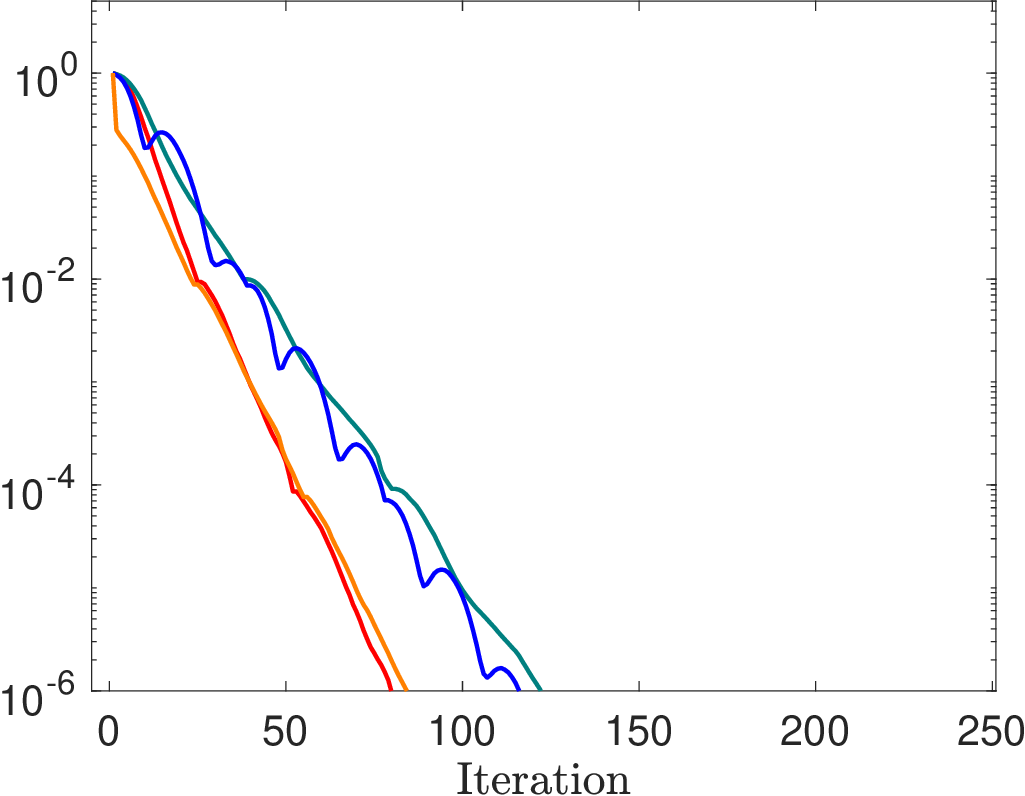}
     \caption{ijcnn1(1e-7) $l=2000$}
   \end{subfigure}%
   \caption{MINRES convergence for four different kernel ridge regression datasets and sketch sizes.}
   \label{fig:krr_convergence_progress}
\end{figure}

We test C-RandRAND and R-RandRAND preconditioners and compare the numbers of MINRES iterations against the Nystr{\"o}m preconditioner. We take $\bOmega = \bX^\mathrm{T}$ as a  uniform sampling matrix. This choice is here the most appropriate  since it enables application of $\bA_\mu \bOmega$ without forming the full $\bA_\mu$.
Analyzing the results in~\cref{tab:krr_results}, we make the following observations.
First, for the wine and abalone datasets, R-RandRAND and C-RandRAND significantly outperform Nystr\"om preconditioning by up to a factor of $2$.
Second, for the MNIST dataset, R-RandRAND consistently yields the lowest iteration counts, outperforming both Nystr\"om and C-RandRAND.
Finally, for ijcnn1, Nystr\"om achieves the lowest iteration counts for large deflation dimensions. This happens due to very rapid decay of the dominant eigenvalues of the kernel matrix, which leads to highly accurate low-rank approximation when using $l = 2000$.  Nevertheless, R-RandRAND remains competitive and clearly superior for smaller sketch sizes.
Similar to the ridge regression experiments, we also plot the residual convergence of MINRES over time in~\cref{fig:krr_convergence_progress} for four representative datasets. Across all cases, preconditioning substantially accelerates convergence compared to unpreconditioned MINRES. Moreover, in the scenarios where R-RandRAND and C-RandRAND show significant improvements over Nystr{\"o}m, the advantage is reflected in a consistently faster decrease of the residual.

\subsubsection{Power iteration on random features regression operators}

To test G-RandRAND preconditioning on real-world indefinite operators, we consider the problem of computing eigenpairs of random features operators associated with the shuttle dataset and rbf kernels with parameters $\gamma = 4/9$ and $13/3$. The target eigenvalues are located near $10^{-7}$, $10^{-9}$, and $10^{-11}$ and are approximated via power iteration. This setting involves negative shift parameters, $\mu = -10^{-7}$, $\mu = -10^{-9}$, and $\mu = -10^{-11}$, and leads to indefinite $\bA_{\mu}$. The right-hand side is chosen as an approximation to the dominant eigenvector of $\bA_{\mu}^{-1}$ obtained after two power iterations.

We solve the resulting indefinite shifted systems using G-RandRAND preconditioners. For comparison, we also employ preconditioners given by the square roots of the C-RandRAND preconditioners constructed for $\bA_\mu^2$ and zero shift, as described in~\cref{rmk:g_randrand_vs_c_randrand}. In all methods, the coefficient $\tau$ is selected according to~\cref{eq:g_randrand_tau,eq:c_randrand_tau}, with $\rho = 0.5$. The unpreconditioned MINRES method requires prohibitively many iterations and is therefore omitted in this test.

\Cref{tab:piter_results} reports the number of iterations required to reach a relative residual of $10^{-10}$. Two key observations follow.
First, G-RandRAND offers significantly greater numerical robustness than C-RandRAND on indefinite systems. In particular, C-RandRAND with $q=0$ fails when $|\mu| \leq 10^{-9}$. This numerical breakdown occurs due to squaring the operator condition number.  The second key observation is that G-RandRAND with $q=1$ dramatically outperforms the alternatives. For example, with $q=1$ and $l=400$, G-RandRAND requires only $20$ MINRES iterations for the operator with $(\gamma,\mu,\sigma_{\min}) = (4/9,-10^{-11},2.1 \cdot 10^{-13})$, compared to $335$ iterations for C-RandRAND with $q=1$ and $485$ for G-RandRAND with $q=0$. This performance gain arises from the combination of the ability of G-RandRAND to exploit the decay of singular values below $|\mu|$, thereby reducing the largest singular value to $\mathcal{O}(|\mu|)$(see~\cref{thm:cnd_bound3}), with the ability to  lift the smallest singular value to $\mathcal{O}(|\mu|)$ (see~\cref{rmk:g_rand_rand_smallest}).

\begin{table}[H]
\small 
\centering
{%
\setlength{\tabcolsep}{2pt} % <--- local reduction of column spacing
\renewcommand{\arraystretch}{1.0}
\begin{tabular}{@{}c
                cccc @{\hspace{15pt}}
                cccc @{\hspace{15pt}}
                cccc @{\hspace{15pt}}
                cccc @{\hspace{15pt}}
                cccc@{}}
\toprule
$\gamma,|\mu|$ 
& \multicolumn{4}{c}{Sketch sizes $l$}
& \multicolumn{4}{c}{C-RAND $q=0$} 
& \multicolumn{4}{c}{C-RAND $q=1$} 
& \multicolumn{4}{c}{G-RAND $q=0$} 
& \multicolumn{4}{c}{G-RAND $q=1$} \\
\cmidrule(lr){2-5}\cmidrule(lr){6-9}\cmidrule(lr){10-13}\cmidrule(lr){14-17}\cmidrule(lr){18-21}
& $l_1$ & $l_2$ & $l_3$ & $l_4$
& $l_1$ & $l_2$ & $l_3$ & $l_4$
& $l_1$ & $l_2$ & $l_3$ & $l_4$
& $l_1$ & $l_2$ & $l_3$ & $l_4$
& $l_1$ & $l_2$ & $l_3$ & $l_4$ \\
\midrule
4/9,~1e-7
& 50  & 100 & 200 & 400
& 215 & 166 & 164 & 151
& 149 & 138 & 129 & 124
& 452 & 323 & 215 & 223
& 162 &  38 &  10 &   \textbf{8} \\
4/9,~1e-9
& 50  & 100 & 200 & 400
& 1667 & >5e3 & >5e3 & >5e3
& 1535 & 233 & 226 & 205
& 2019 & 388 & 302 & 233
& 1808 &  29 &   8 &   \textbf{8} \\
4/9,~1e-11
& 200 & 400 & 800 & 1600
& >5e3 & >5e3 & >5e3 & >5e3
& 449 & 335 & 297 & 266
& 808 & 485 & 383 & 311
& 227 &  20 &  12 &  \textbf{12} \\
13/3,~1e-7
& 200 & 400 & 800 & 1600
& 1466 & 812 & 741 & 644
& 1297 & 600 & 537 & 501
& >5e3 & 1300 & 1092 & 611
& 1521 & 163 &  25 &   \textbf{9} \\
13/3,~1e-9
& 500 & 800 & 1200 & 1600
& >5e3 & >5e3 & >5e3 & >5e3
& 1513 & 995 & 884 & 897
& 2482 & 1635 & 1254 & 1130
& 1257 & 166 &  35 &  \textbf{20} \\
13/3,~1e-11
& 1000 & 1200 & 1600 & 
& >5e3 & >5e3 & >5e3 & 
& 1876 & 1536 & 1442 & 
& >5e3 & 2679 & 1804 & 
& 1034 & 421 & \textbf{69} &  \\
\bottomrule
\end{tabular}%
}
\caption{Convergence results with indefinite systems constructed with ridge regression operators for the shuttle dataset and negative $\mu$.}
\label{tab:piter_results}
\end{table}

\subsubsection{Portfolio optimization}

Here, we go beyond the case of solving a single linear system and focus on the setting of optimization methods, where one must solve sequences of multiple linear systems.

We consider a risk-adjusted return problem with a factor model~\cite{boyd2017multi,grinold2000active}, related to the seminal mean–variance model of portfolio selection~\cite{markowitz2008portfolio}. The problem can be formulated as the quadratic program
\begin{equation}\label{eq:portfolio}
\begin{split}
    \min_{\by \in \mathbb{R}^m, \bz \in \mathbb{R}^s} ~ & - \gamma^{-1} \bm ^\mathrm{T} \by + \by^\mathrm{T} \bD \by + \bz^\mathrm{T} \bz \\
    \text{s.t.}~ & \bF^\mathrm{T} \by = \bz \\
    &\bM \by \leq \bu \\
    &\mathbf{1}^\mathrm{T} \by  = 1, \by \geq \bnull.
\end{split}
\end{equation}
Here, $\bm \in \mathbb{R}^{m}$ is the vector of mean returns, $\bF \in \mathbb{R}^{m \times s}$ and $\bD \in \mathbb{R}^{m \times m}$ are the  loading and diagonal factor matrices that describe the covariance matrix of the assets, $\bM \in \mathbb{R}^{p \times m}$ and $\bu \in \mathbb{R}^{p}$ encode additional linear constraints on the holdings, and $\gamma > 0$ is the risk-aversion parameter.
The last set of constraints represents the fact the holdings $\by = [y_1, \hdots, y_m]$ form a probability distribution.

To solve~\cref{eq:portfolio} we implemented the inexact interior-point proximal method of multipliers (IP-PPM)~\cite{pougkakiotis2021interior,chu2024randomized}.
At each iteration, the predictor–corrector refinement requires solving two Newton systems of the form
\begin{equation}\label{eq:newton_system}
\left(\bZ^\mathrm{T} \widehat{\bD}_i \bZ + \mu_i \bI \right)\bx = \bb_i .
\end{equation}
Here, the diagonal matrix  $\widehat{\bD}_i$, the regularization parameter $\mu_i>0$, and right-hand-side $\bb_i$ depend on the iteration. The matrix $\bZ$, formed from the data matrices $\bF$ and $\bM$, is fixed across iterations.

Similar to~\cite{chu2024randomized}, we constructed a synthetic portfolio instance whose covariance spectrum has a few large exponentially decaying eigenvalues followed by a slowly decaying tail. In our experiment, we consider $m=10000$ assets, $p=6250$ inequality constraints, and $s=100$ factors. The resulting matrix $ \bZ$ has $n = p+s+1 = 6351$ columns and $m+n = 16351$ rows.
Further setup details are given in~\cref{appendix_c}.

The linear systems~\cref{eq:newton_system} arising from  each Newton step were solved with MINRES method to machine precision. Then we tested randomized preconditioners to solve these systems. The average number of (preconditioned) MINRES iterations to achieve accuracy tol $= 10^{-10}$ or $10^{-6}$, are shown in~\cref{tab:portfolio_dataset}.

\begin{table}[H]
   \centering
   \begin{tabular}{@{}lllrrrr@{}} \toprule
    Sketch size $l$ & tol & MINRES & Nystr\"{o}m & C-RAND* &  C-RAND & R-RAND \\ \midrule
    100 & 1e-6 & >900 & 214 & 216 &  216 & \textbf{173} \\ 
    100 & 1e-10 & >2500 & 619 & 616 &  612 & \textbf{514} \\ 
    400 & 1e-6 & >900 & 188 & 188 &  187 &  \textbf{155} \\ 
    400 & 1e-10 & >2500 & 542 & 547 &  547 & \textbf{464} \\   
    \bottomrule
   \end{tabular}
   \caption{Portfolio problem: average number of MINRES iterations to converge to relative residual tol $= 10^{-10}$ or $10^{-6}$. C-RAND* uses $\tau=\|\bPi\bA^{-1}\bPi\|$, whereas C-RAND uses $\tau$ chosen by minimizing~\cref{eq:cnd_bound2_ext} taking $F = \lambda_{\min}(\bA_\mu)$. }
   \label{tab:portfolio_dataset}
\end{table}

\begin{figure}[H]
   \begin{subfigure}{.253\textwidth}
     \centering
     \includegraphics[width=\linewidth]{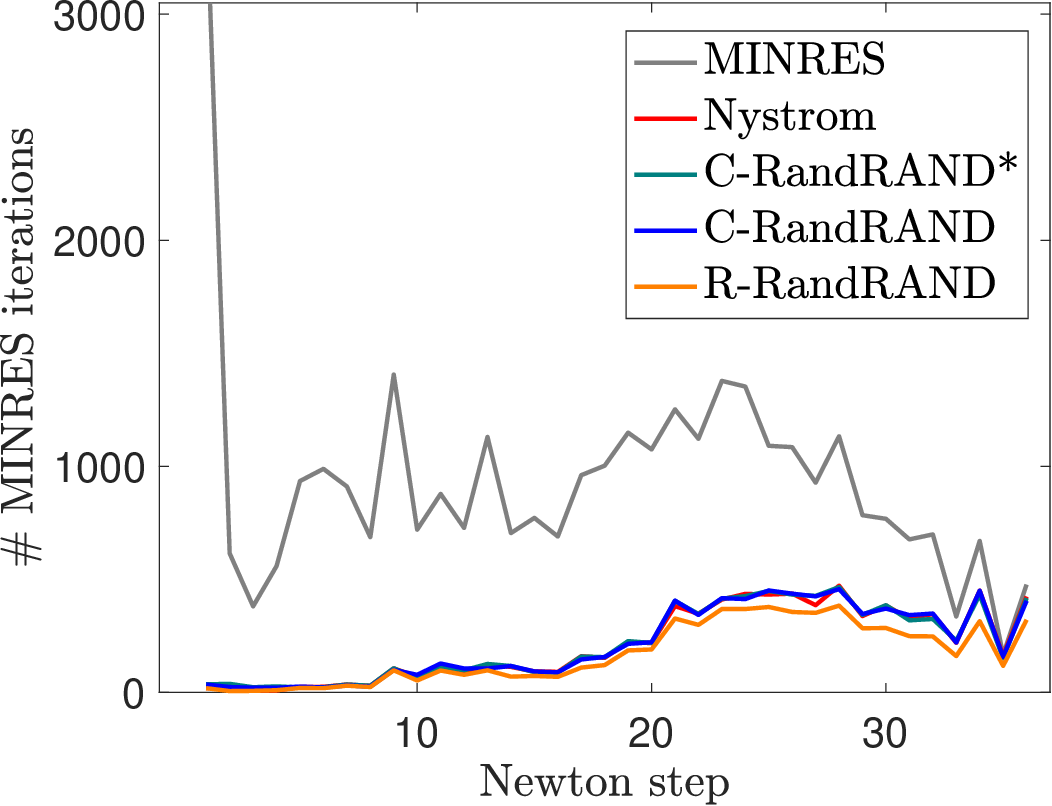}
     \caption{ $l=100$, tol = 1e-6}
   \end{subfigure}%
   \hspace{.0025\textwidth}
   \begin{subfigure}{.24\textwidth}
     \centering
     \includegraphics[width=\linewidth]{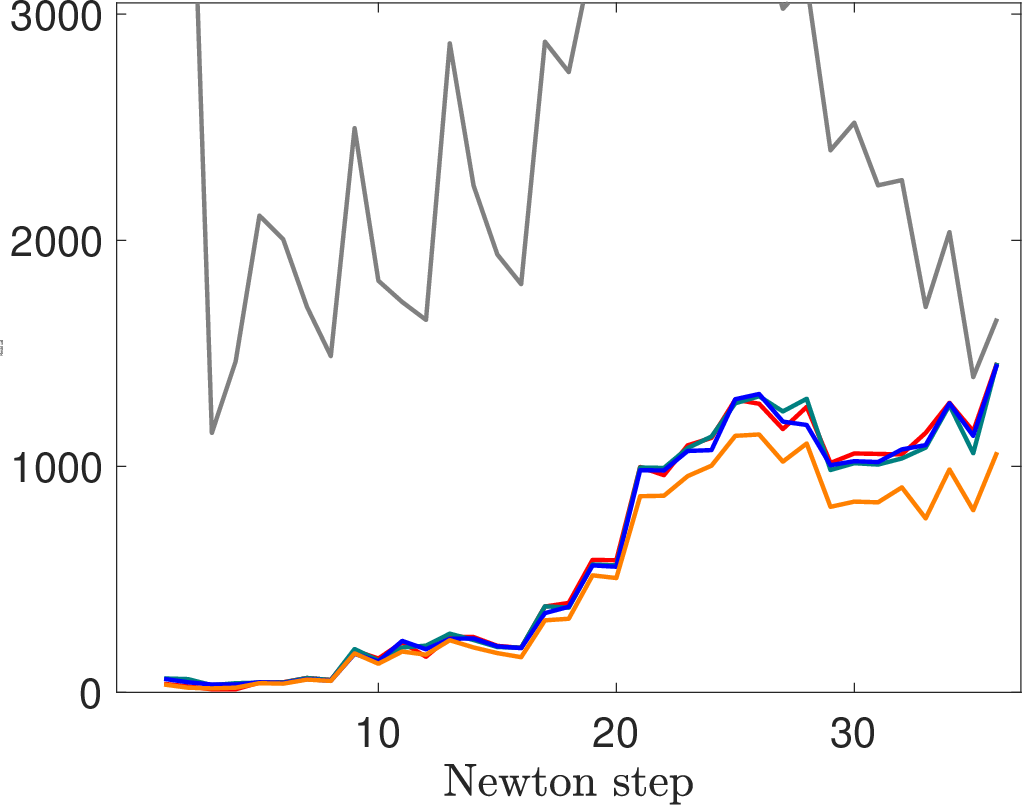}
     \caption{$l=100$, tol = 1e-10}
   \end{subfigure}
   \hspace{.0025\textwidth}
   \begin{subfigure}{.24\textwidth}
     \centering
     \includegraphics[width=\linewidth]{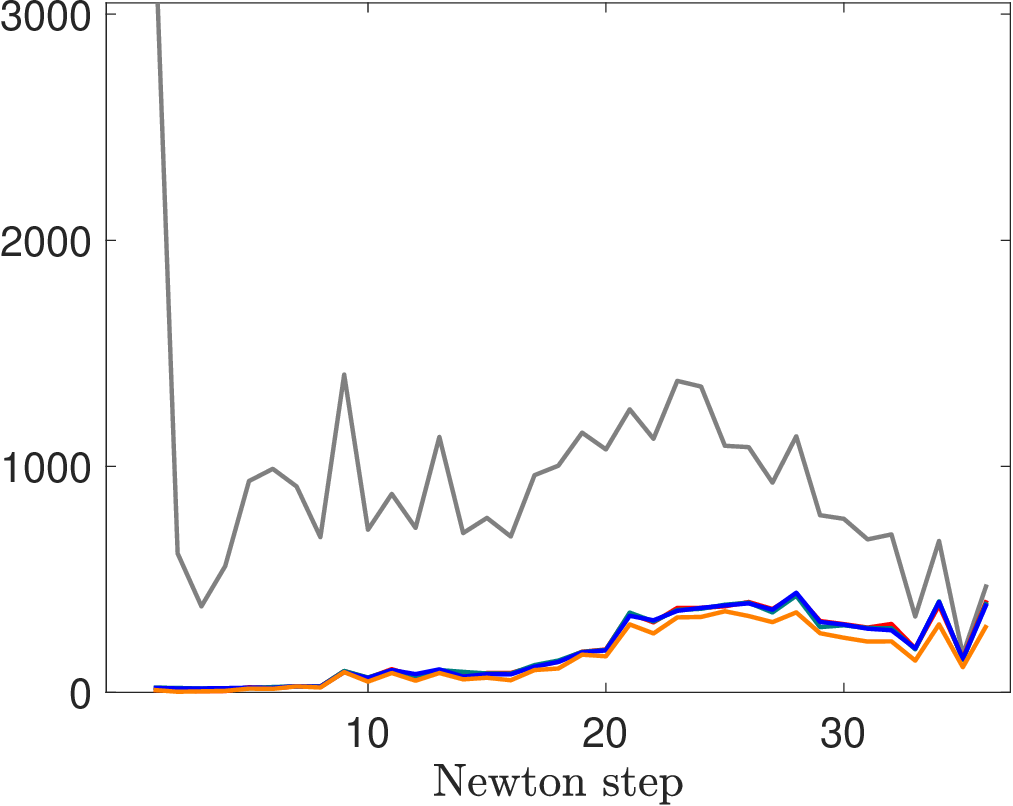}
     \caption{$l=400$, tol = 1e-6}
   \end{subfigure}%
   \hspace{.0025\textwidth}
   \begin{subfigure}{.24\textwidth}
     \centering
     \includegraphics[width=\linewidth]{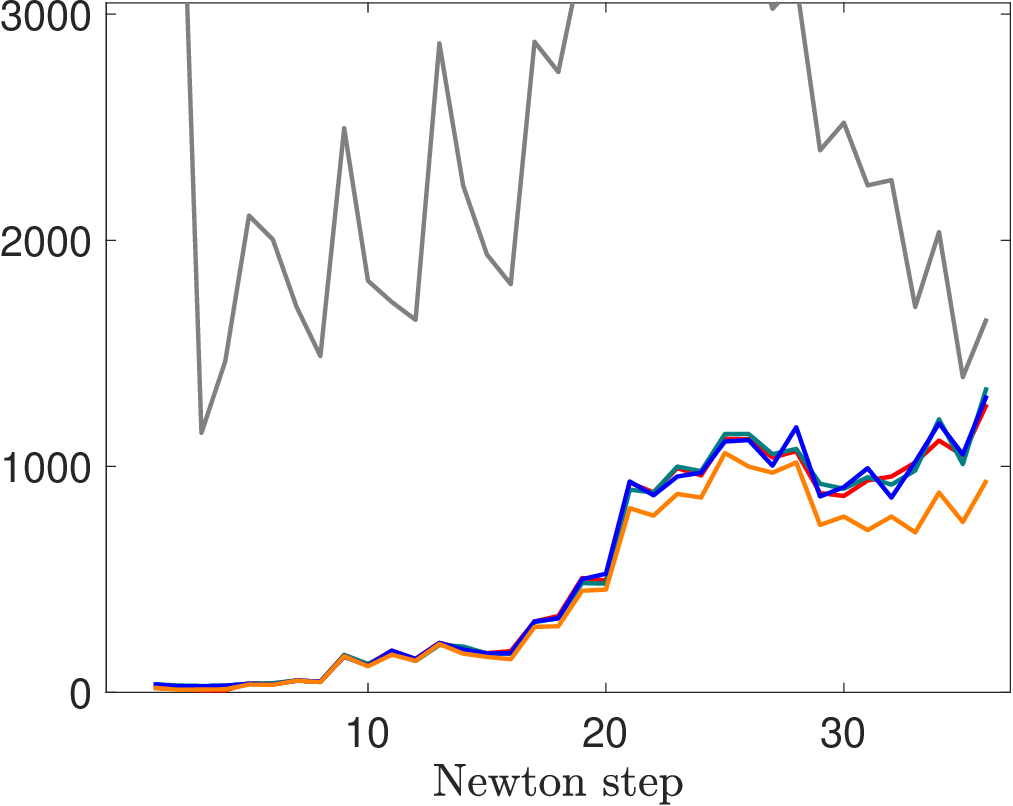}
     \caption{$l=400$, tol = 1e-10}
   \end{subfigure}%
   \caption{Number of MINRES iterations for solving portfolio optimization problems with different accuracy tolerances.}
   \label{fig:portfolio_converge}
\end{figure}

Like with previous experiments, the use of randomized preconditioning yields faster convergence. R-RandRAND has the lowest average number of iterations, requiring about 15\%-20\% fewer iterations than with C-RandRAND and Nystr\"{o}m preconditioner.
The number of MINRES iterations at each Newton step of the optimization algorithm is shown in~\cref{fig:portfolio_converge}. 
In the first half of Newton steps, all preconditioned methods appear similar and significantly outperform MINRES.
Later on, we start to see the increase in numbers of preconditioned iterations, and differences in performance between preconditioners.
%Indeed, the increased difficulty in solving systems in the latter Newton steps is well documented in the literature~\cite{gondzio2012interior}.
This can be explained by the fact that as the primal-dual solutions get closer to the boundary, the Newton systems start to have a spike of small eigenvalues. 

\begin{figure}[H]
   \begin{subfigure}{.253\textwidth}
     \centering
     \includegraphics[width=\linewidth]{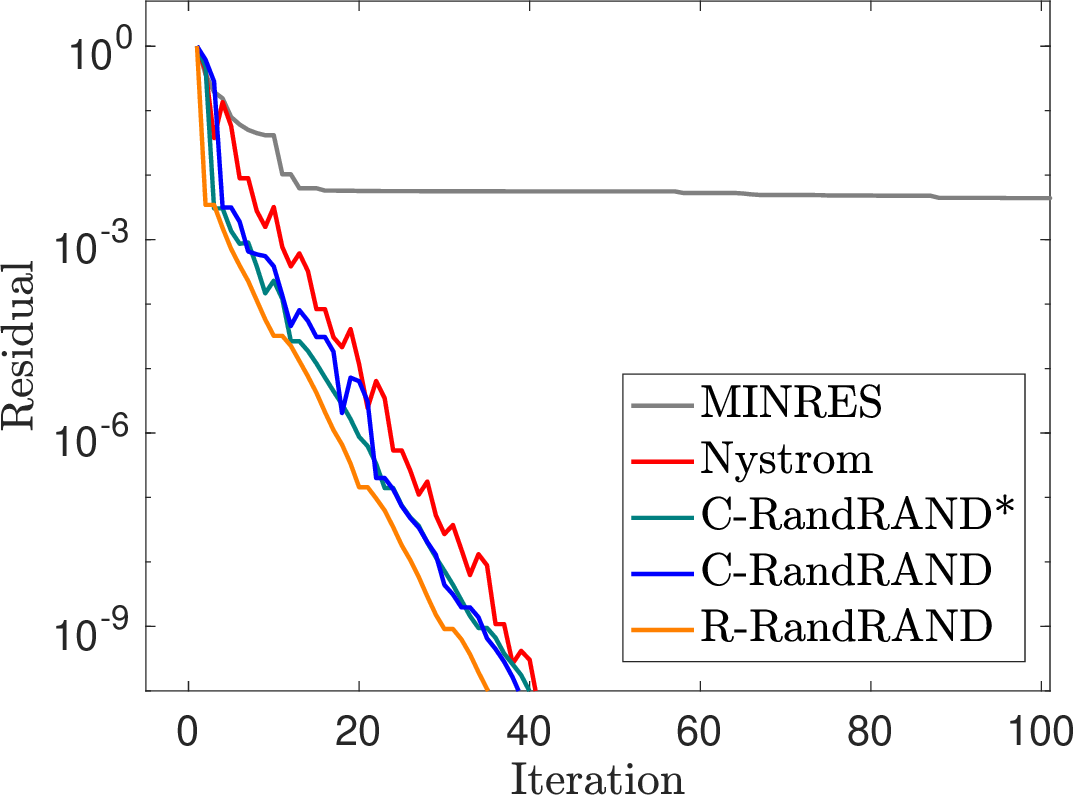}
     \caption{Newton step 5}
   \end{subfigure}%
   \hspace{.0025\textwidth}
   \begin{subfigure}{.24\textwidth}
     \centering
     \includegraphics[width=\linewidth]{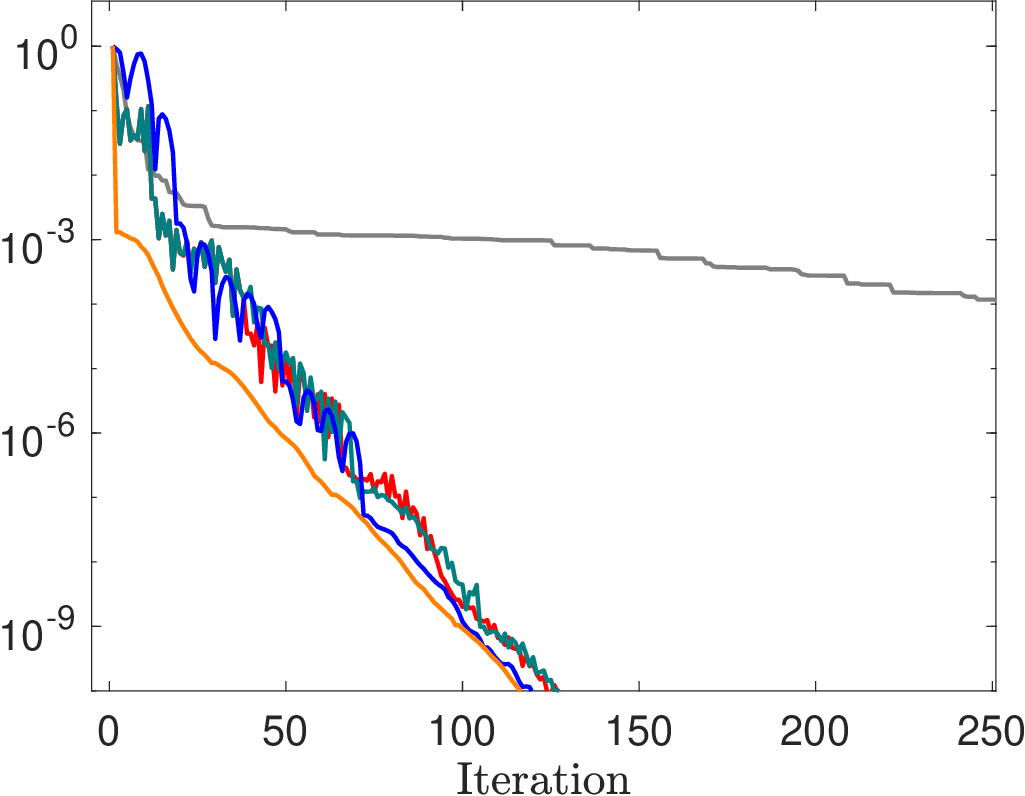}
     \caption{Newton step 10}
   \end{subfigure}
   \hspace{.0025\textwidth}
   \begin{subfigure}{.24\textwidth}
     \centering
     \includegraphics[width=\linewidth]{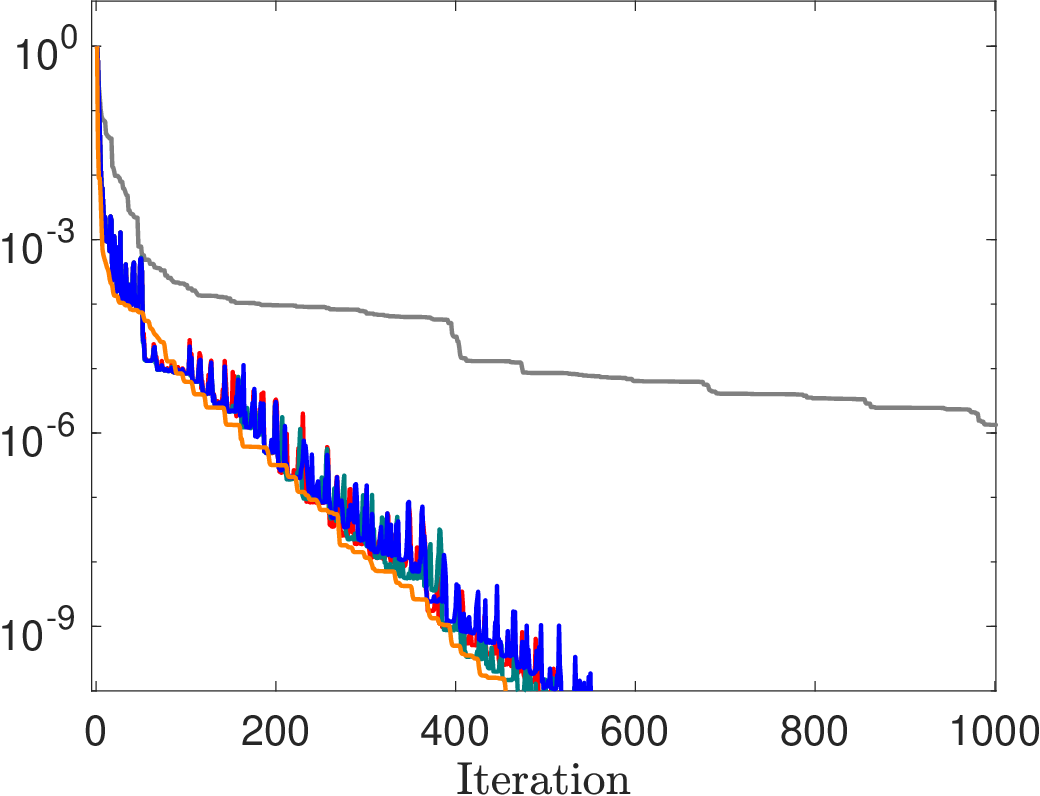}
     \caption{Newton step 20}
   \end{subfigure}%
   \hspace{.0025\textwidth}
   \begin{subfigure}{.24\textwidth}
     \centering
     \includegraphics[width=\linewidth]{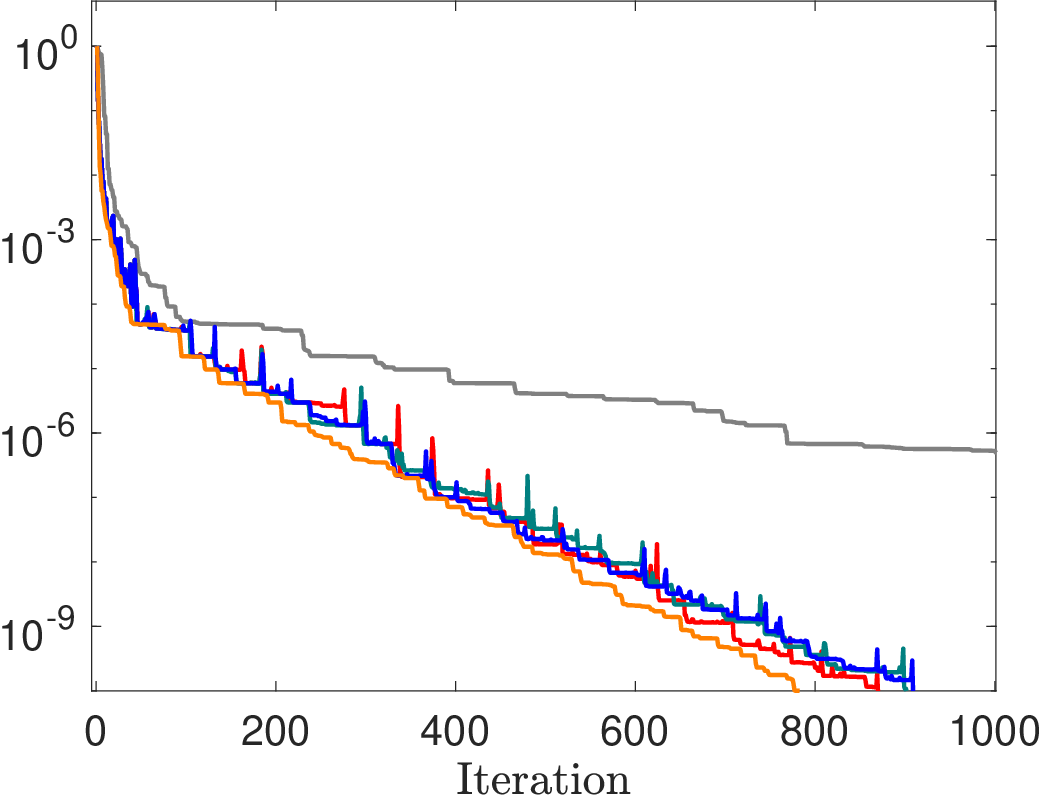}
     \caption{Newton step 30}
   \end{subfigure}%
    \caption{MINRES convergence for four different Newton systems.  Preconditioners are constructed with sketching size $l=400$.}
   \label{fig:portfolio_spectra}  
\end{figure}

Finally, to validate the improvement brought by preconditioning, we examine the convergences of several particular Newton systems in~\cref{fig:portfolio_spectra}.
Overall we observe similar trends as in the ridge regression experiments. C-RandRAND and Nystr\"{o}m exhibit similar convergence curves with occasional irregularities.  R-RandRAND exhibits a faster and smoother residual decay than other preconditioners, which suggests better robustness.

\subsection{Applications of Basis-less RandRAND}
\label{sec:basis-less-pde-appn}

 Finally, we showcase the performance of basis-less RandRAND preconditioners on limited-memory kernel ridge regression and an inverse PDE problem, where the basis-explicit preconditioning methods incur a dramatic computational overhead.

\subsubsection{PDE-constrained Inverse Problem}
We consider a PDE-constrained inverse problem, which represents the acoustic invisibility cloak benchmark from~\cite{balabanov2021randomized}. 
This problem describes a 2D acoustic wave scattering with a scatterer covered by a cloak composed of alternating layers of mercury and isotropic light liquids. 
The cloak aims to minimize the scatterer's visibility in a certain frequency band. 
Two unit sources are placed outside the cloak and need to be recovered from noisy pressure field measurements $u + \epsilon$. 
We set $\epsilon(x) = 0.25\,\omega(x)\,\lvert u(x)\rvert,$ where $\omega(x)$ is a standard normal random variable. 

After discretization, the problem becomes
\begin{equation}\label{eq:discrete_pde}
   \min_{\bff}\;\; \tfrac12 \bigl(\bK^{-1}\bM\,\bff - \mathbf{z}\bigr)^\mathrm{T}\,\bM\,\bigl(\bK^{-1}\bM\,\bff - \mathbf{z}\bigr) \;+\; \tfrac{\mu}{2} \,\bff^{\mathrm{T}} \bM\,\bff,
\end{equation}
where $\bM\in \mathbb{R}^{n \times n}$, with $n \approx 5 \cdot 10^5$, is the mass matrix and $\bK\in \mathbb{C}^{n \times n}$ is the stiffness matrix associated with the PDE.
The matrices $\bM$ and $\bK$ are sparse and permit precomputed efficient factorizations. 
The solution of~\cref{eq:discrete_pde} 
can be obtained by solving
$$
   (\bA + \mu \bI)\,\bx = \mathbf{b},
$$
with $\bA = \bC^{\mathrm{T}}\,\bK^{-1}\bM\,\bK^{-1}\bC$ and $\mathbf{b} = \bC^{\mathrm{T}}\,\bK^{-1}\bM\,\mathbf{z}$, where $\bC$ is (possibly permuted) Cholesky factor of $\bM$. We set $\mu = 10^{-6}\|\bA\|$. 
A detailed overview of this numerical example is provided in~\cref{appendix_c}.

\begin{figure}[ht!]
    \centering
    \begin{subfigure}{.25\textwidth}
        \includegraphics[width=\linewidth]{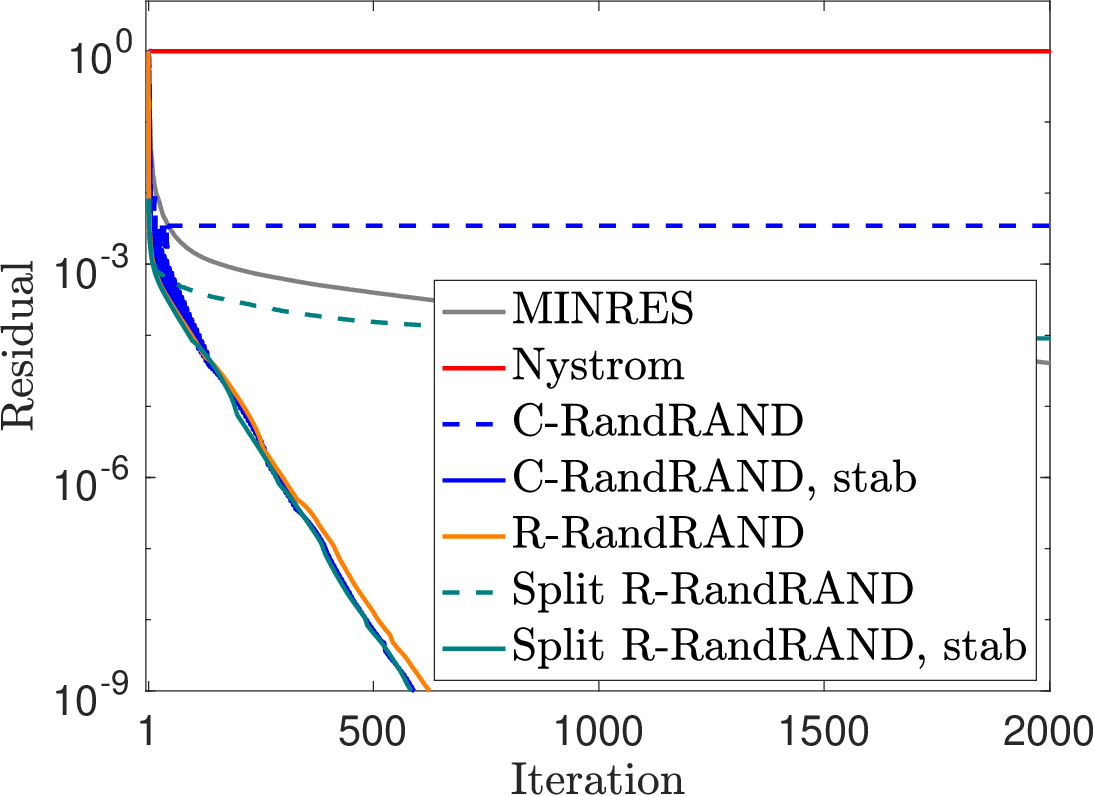}
        \caption{$l = 1000$}
    \end{subfigure}
    \hspace{.00125\textwidth}
    \begin{subfigure}{0.237\textwidth}
        \includegraphics[width=\linewidth]{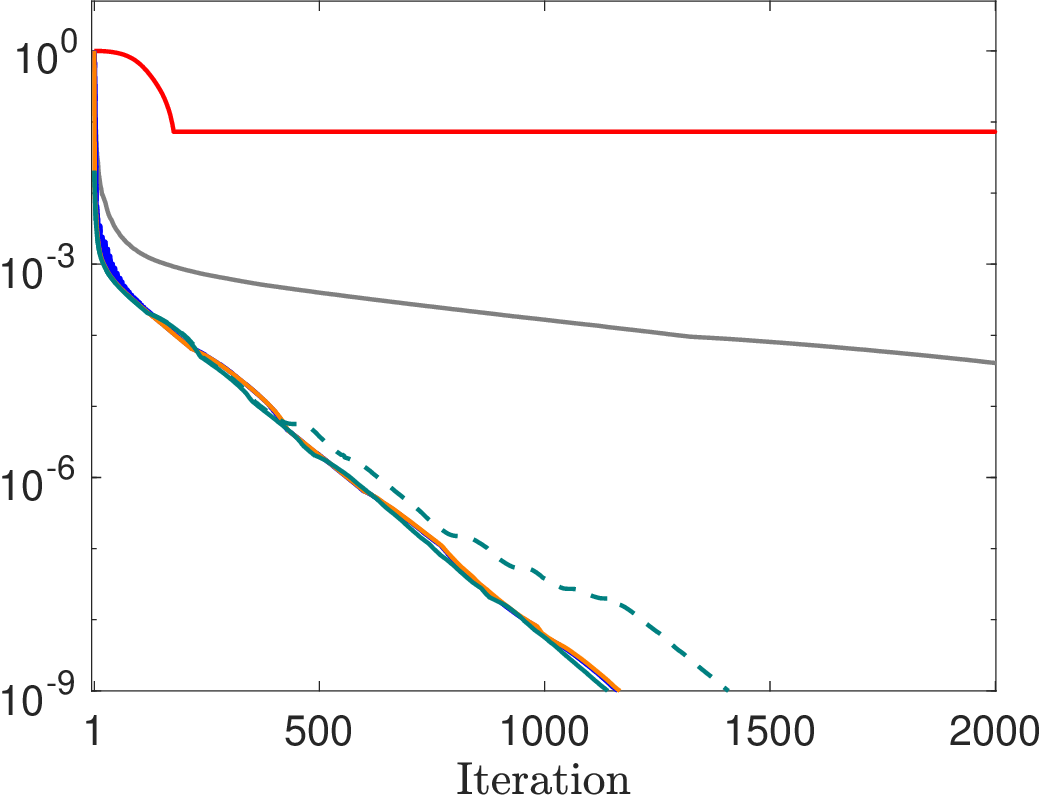}
        \caption{$l = 500$}
    \end{subfigure}
    \hspace{.00125\textwidth}
    \begin{subfigure}{0.237\textwidth}
        \includegraphics[width=\linewidth]{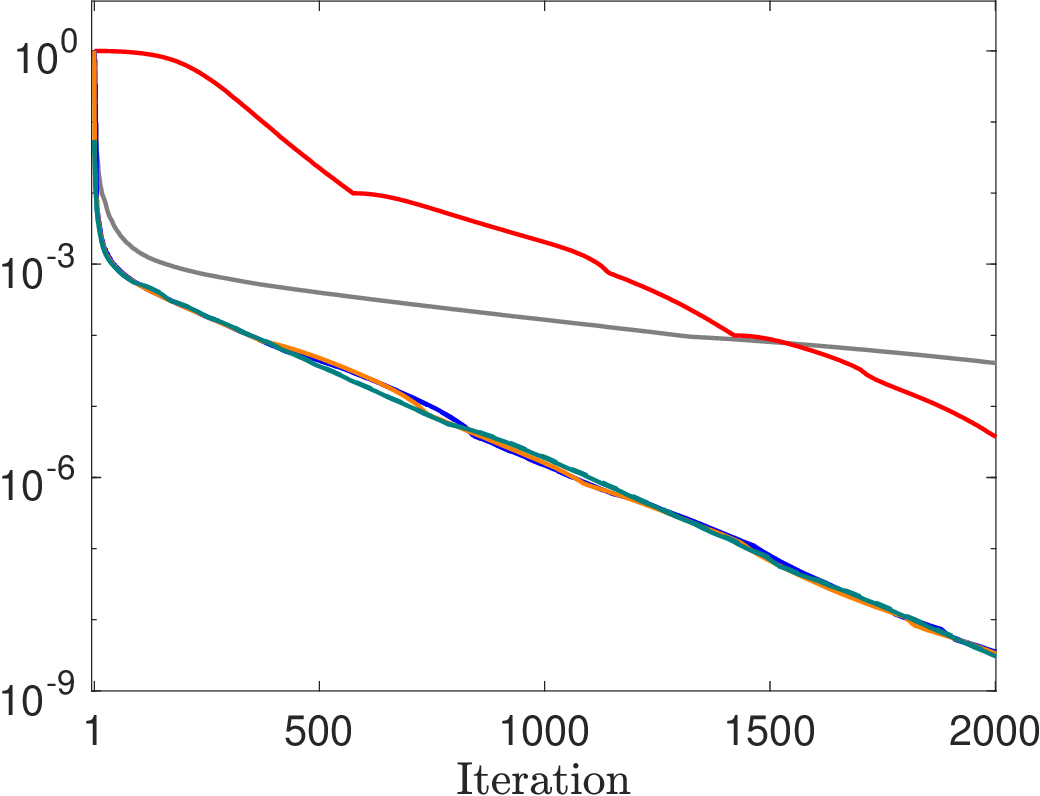}
        \caption{$l = 250$}
    \end{subfigure}
    \hspace{.00125\textwidth}
    \begin{subfigure}{0.237\textwidth}
        \includegraphics[width=\linewidth]{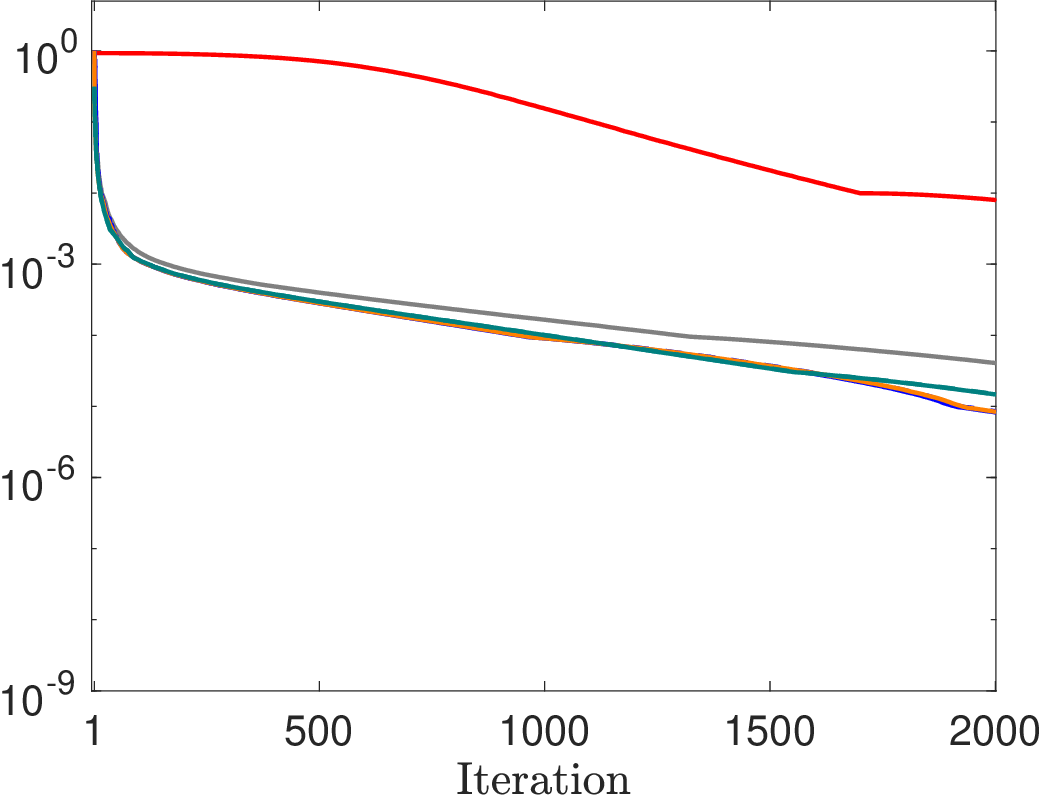}
        \caption{$l = 100$}
    \end{subfigure}
    \caption{RandRAND preconditioning for inverse PDE. C-RandRAND uses $\tau$ chosen by minimizing~\cref{eq:cnd_bound2_ext} taking $F = \mu^{-1}$.}
    \label{fig:invcloak}
\end{figure}

The MINRES here exhibits dramatically slow convergence, motivating the need for randomized preconditioning. 
Since in this numerical example, operations with dense sketching operators and explicit orthogonalization of the basis pose a bottleneck, we use basis-less approach with sparse+Gaussian embeddings. 
We construct $\bX$ as $\bX_1^\mathrm{T}\,\bX_2^\mathrm{T}$, where $\bX_1$ is an $l_1 \times n$ random sparse operator of density $\gamma = \log(n)$ and dimension $l_1 = \log(n)\,l$, and $\bX_2$ is an $l \times l_1$ Gaussian matrix. A matvec with such $\bX$ has negligible computational cost (from $5\%$ to $10\%$) compared to a matvec with $\bA$. 
The sketch sizes considered here are $l = 100$, $l = 250$, $l = 500$ or $l = 1000$, and the power parameter $q=1$. 

We tested basis-less implementations of the right R-RandRAND, split R-RandRAND, and C-RandRAND preconditioners. In our implementation, applying the right R-RandRAND requires four matvecs with $\bA_\mu\bOmega$, while the split R-RandRAND and C-RandRAND variants require only two matvecs\footnote{To attain this cost of C-RandRAND, in \cref{alg:BL_RandRAND_PCG} we used expression $\bP = \bI - \bPi +\tau \widetilde{\bPi} = \bI + \bQ(\tau\bG \bQ^\mathrm{T} - \bQ^\mathrm{T})$}. In turn, each matvec with $\bA_\mu\bOmega$ corresponds to $q+1 = 2$ matvecs with $\bA_\mu$. We considered two options for computing the R factor: the routine $\mathtt{fast{\_}q{\_}less{\_}chol{\_}qr}$ and its stabilized version, $\mathtt{fast{\_}q{\_}less{\_}precond{\_}chol{\_}qr}$ introduced in~\cref{basis_free_represent}.

For comparison, we also implemented preconditioning based on the basis-less Nystr{\"o}m approximation (described  e.g. in~\cite{derezinski2024faster}). See~\cref{appendix_c} for implementation details of this preconditioner. 

The resulting MINRES convergence curves are shown in~\cref{fig:invcloak}. 
It is revealed that i) the preconditioners drastically reduce the number of MINRES iterations when they are numerically stable and $l \geq 250$; ii) C-RandRAND and R-RandRAND perform similarly and yield a significant \blue{convergence} improvement over Nystr{\"o}m preconditioning; iii) the Nystr{\"o}m preconditioning exhibits dramatic numerical instability for $l\geq 500$ due to a high condition number of the basis matrix $\bA_\mu \bOmega$ while  R-RandRAND and C-RandRAND (when using $\mathtt{fast{\_}q{\_}less{\_}precond{\_}chol{\_}qr}$) are numerically robust.

\subsubsection{Limited-memory Kernel Ridge Regression} \label{lim_mem_krr}

Basis-less preconditioning can be particularly advantageous for large-scale kernel problems under extreme memory constraints. Consider kernel ridge regression~\cref{eq:kernel} in a scenario where storing several full vectors is computationally infeasible. The kernel matrix $\bA = \frac{1}{n} \bK$ can be accessed only in small blocks of columns, and requires matvecs to be computed on the fly. In this regime, one must rely on methods that construct and apply the preconditioner or the inverse low-rank approximation without explicit computation of the approximation basis.

Here, it is most appropriate to choose the basis $\bA_\mu \bOmega$ as uniformly sampled columns of $\bA_\mu$. This choice allows efficient on-the-fly application of $\bA_\mu \bOmega$ at a cost substantially lower than a full matvec with $\bA_\mu$. Furthermore, due to the low application cost of $\bA_\mu \bOmega$, the basis can be orthogonalized via the preconditioned basis-less QR orthogonalization ($\mathtt{fast{\_}q{\_}less{\_}precond{\_}chol{\_}qr}$) without a significant computational overhead.

We evaluate basis-less implementations of the right R-RandRAND and C-RandRAND preconditioners, and compare them with a basis-less Nystr{\"o}m  preconditioner from e.g.~\cite{derezinski2024faster} (see~\cref{appendix_c} for details).   
Furthermore, along with the basic representation of C-RandRAND given by~\cref{eq:g_randrand_explicit} we also test the re-orthogonalized version given by~\cref{eq:c_randrand_reorth} which is more numerically robust. Again, due to the relatively low application cost of $\bA_\mu \bOmega$, the re-orthogonalization here entails only a minor computational overhead. 

\Cref{tab:krr_dataset2} summarizes the datasets, rbf kernel hyperparameters, regularization parameters $\mu$, and residual tolerances. As in~\cref{krr}, labels are generated via one-vs-all encoding. The linear systems are here solved only for the first right-hand-side vector. For all randomized preconditioners, we use a sketching dimension $k=5000 \ll n$.

\begin{table}[H]
   \centering
   \begin{tabular}{@{}llrrrr@{}} \toprule
    Dataset &  $n$ & $\gamma$ & $\mu$ & tol  \\ \midrule
    EMNIST &  112800 &  1/64  & 1e-3/$n$,1e-4/$n$ & 5e-4 \\
    MNIST  & 60000  &  1e-2    & 1e-4/$n$, 1e-5/$n$ & 1e-5  \\
    ijcnn1 & 49990  &  1, 1/2    & 1e-5/$n$, 1e-6/$n$ & 1e-5  \\ \bottomrule
   \end{tabular}
   \caption{Description of kernel ridge regression datasets in the limited memory setting.}
   \label{tab:krr_dataset2}
\end{table}

\begin{figure}[H]
   \begin{subfigure}{.253\textwidth}
     \centering
     \includegraphics[width=\linewidth]{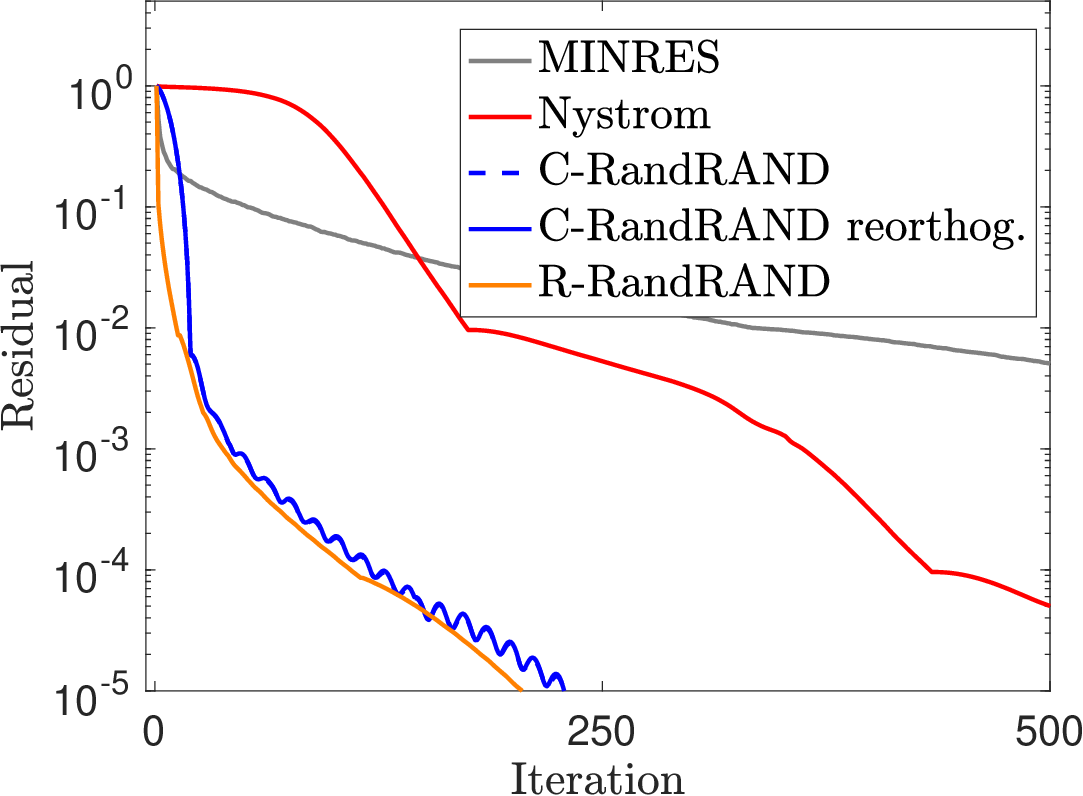}
     \caption{MNIST, $1\text{e-2}, 1\text{e-4}/n$ }
   \end{subfigure}%
    \hspace{.00125\textwidth}
 \begin{subfigure}{.24\textwidth}
     \centering
     \includegraphics[width=\linewidth]{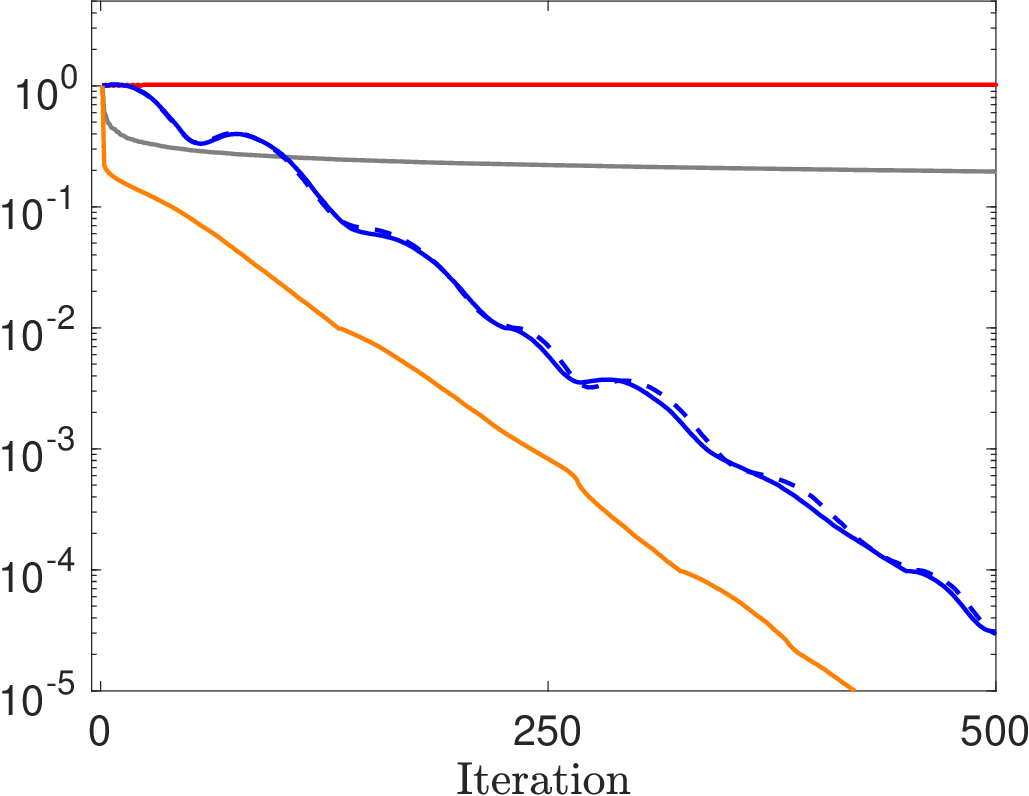}
     \caption{ijcnn1, $1, 1\text{e-5}/n$}
   \end{subfigure}
    \hspace{.00125\textwidth}
   \begin{subfigure}{.24\textwidth}
     \centering
     \includegraphics[width=\linewidth]{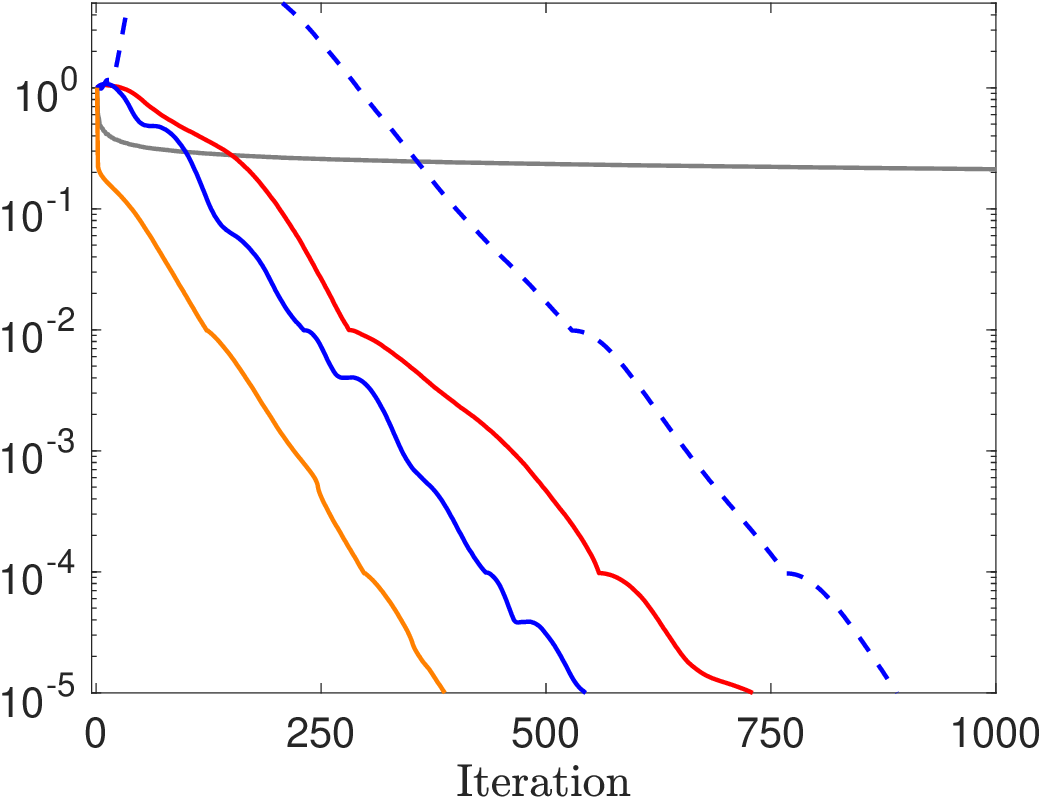}
     \caption{ijcnn1, $0.5, 1\text{e-6}/n$ }
   \end{subfigure}%
    \hspace{.00125\textwidth}
   \begin{subfigure}{.24\textwidth}
     \centering
     \includegraphics[width=\linewidth]{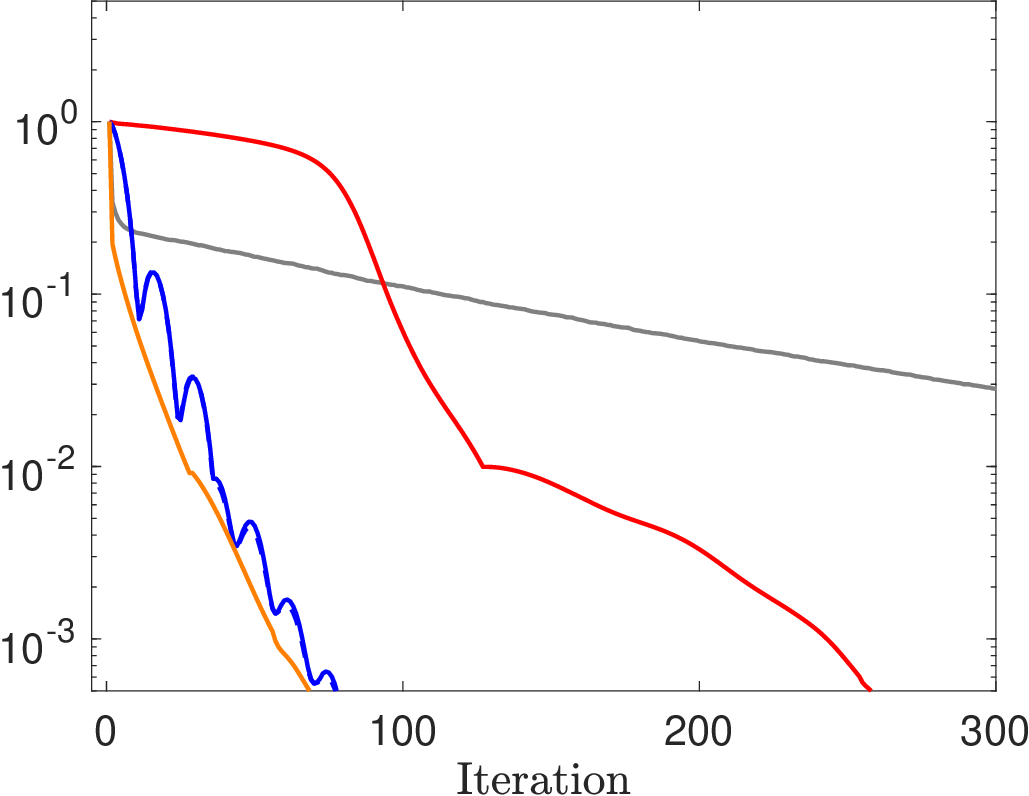}
     \caption{EMNIST, $1/64$, $1\text{e-3}/n$ }
   \end{subfigure}%
   \caption{MINRES convergence in limited memory kernel ridge regression test cases.}
   \label{fig:krr_convergence_progress2}
\end{figure}

\Cref{tab:krr_dataset3} reports the number of MINRES iterations required for convergence of each method. \Cref{fig:krr_convergence_progress2} shows the observed residual convergence in four test cases. RandRAND significantly reduces iteration counts compared to the unpreconditioned solver and consistently outperforms Nystr{\"o}m preconditioning in convergence rates and numerical robustness. For example, for the EMNIST dataset with $\gamma = 1/64$ and $\mu = 1\text{e-3}/n$, C-RandRAND and R-RandRAND require roughly $3.3\times $ and $3.8 \times$ fewer iterations, respectively, than Nystr{\"o}m method. On the MNIST test cases, the basis-less Nyström preconditioner suffers from severe numerical instabilities, leading to early termination of MINRES. In contrast, RandRAND preconditioners remain robust and  converge with similar rate as their basis-explicit counterparts. 

Notably, for the ijcnn1 dataset with $\gamma = 0.5$ and $\mu = 1\text{e-5}/n$, the basic C-RandRAND variant~\cref{eq:g_randrand_explicit} shows moderate instability and suboptimal convergence. This is fully remedied by re-orthogonalization~\cref{eq:c_randrand_reorth} of the C-RandRAND components. This example again demonstrates the benefit of representation of preconditioners through orthogonal projections in ensuring numerical robustness.

\begin{table}[H]
\centering
\begin{tabular}{@{}lcccccc@{}} \toprule
Dataset, $\gamma$, $\mu$ & \# MINRES & \# iter Nystr{\"o}m & \multicolumn{2}{c}{\# iter C-RAND} & \# iter R-RAND \\ 
\cmidrule(lr){4-5}
& & & basic & reorthog & \\ \midrule
EMNIST, $1/64$, $1\text{e-3}/n$ & >500 & 258   & 77  & 77  & \textbf{68}  \\
EMNIST, $1/64$, $1\text{e-4}/n$ & >500 & >500   & 83  & 83  & \textbf{76}  \\
MNIST, $1\text{e-2}, 1\text{e-4}/n$  & >1e3 & >1e3 & 228 & 228 & \textbf{205} \\
MNIST, $1\text{e-2}, 1\text{e-5}/n$  & >1e3 & >1e3 & 306 & 311 & \textbf{286} \\
ijcnn1, $1, 1\text{e-5}/n$   & >1e3 & >1e3 & 568 & 564 & \textbf{421} \\
ijcnn1, $1/2, 1\text{e-6}/n$ & >1e3 & 729   & 889 & 543 & \textbf{387} \\ \bottomrule
\end{tabular}
\caption{Numbers of (preconditioned) MINRES iterations needed to achieve the residual tolerance. C-RandRAND uses $\tau$ chosen by minimizing~\cref{eq:cnd_bound2_ext} taking $F = \mu^{-1}$.}
\label{tab:krr_dataset3}
\end{table}

    \section{Conclusion}
\label{sxn:conclusion}

We introduced RandRAND methods that accelerate linear solvers by deflating the spectrum of the operator via efficient orthogonal projections onto random subspaces. Unlike common deflation and randomized approaches, RandRAND avoids computing eigenpairs or low-rank approximations, leading to significant advantages in flexibility, performance, and applicability.

The orthogonal projectors in RandRAND provide flexible design and implementation choices, making it well-suited for diverse problem settings and computational architectures. We described: (i) \emph{basis-explicit variants}, effective for dense systems and for sparse systems with small deflation dimensions; (ii) \emph{basis-less variants} constructed from the random sketch(es) of the basis, especially well-suited for large-scale sparse problems and memory-constrained settings, particularly when combined with sparse sketching, fast transforms, or sampling; and (iii) extensions based on the \emph{least-squares representation}, allowing implementations with approximate orthogonalization and factorization-free approaches which can be appealing for some modern architectures and computation~recycling.

RandRAND comes with rigorous quasi-optimality guarantees that are stronger than those of the randomized preconditioners~\cite{frangella2023randomized,derezinski2024solving,derezinski2024faster} based on low-rank approximations. We established condition number bounds for RandRAND that are purely  \emph{projection-based} and integrated them with an extension of the randomized range approximation to shifted operators. In particular, we showed that the condition numbers of RandRAND's preconditioned operators are, with high probability, bounded by small constants when the sketch sizes are comparable to the effective spectral~dimension. 

Special attention was given to the numerical robustness of RandRAND. In particular, we described strategies for refining the stability of orthogonal projections in the basis-less setting and proposed novel randomized basis-less Cholesky QR algorithms.
These ideas can be used beyond RandRAND. For instance, they could improve  the numerical robustness of low-rank approximation algorithms.

RandRAND has been extensively tested on synthetic problems and real-world benchmarks, including kernel and random feature regression, portfolio optimization, and inverse PDE problems. It consistently delivers great  acceleration for CG and MINRES, and offers notable advantages over low-rank-based randomized preconditioners in terms of convergence rate, numerical robustness, and construction/application cost, especially in the basis-less~regime.

RandRAND is a general framework that applies to spd, indefinite and non-symmetric systems, and potentially can be extended to rectangular systems. Beyond linear systems, RandRAND also opens the door to efficient iterative methods for matrix function evaluations $f(\bA)\bx$, which will be addressed in future work. Another promising direction is to use RandRAND methods to construct second-level preconditioners for enhancing state-of-the-art deterministic~preconditioners.

We demonstrated how RandRAND can be effectively applied across a range of shift parameters $\mu$, a scenario common in kernel methods and numerical optimization. These ideas extend naturally to solving sequences of related linear systems, such as those arising in interior-point and other second-order optimization methods. The design of RandRAND makes it especially attractive for recycling computations in these settings, pointing to a promising research direction.

\paragraph{Acknowledgments.}
We thank Micha{\l} Derezi{\'n}ski, Petros Drineas and Raphael Meyer for insightful discussions.
We acknowledge DARPA, NSF, the DOE Competitive Portfolios grant, and the DOE SciGPT grant
for partial support of this~work. 
In addition, Caleb Ju acknowledges DOE CSGF grant under Award Number DE-SC0022158.

  \bibliographystyle{plain}
  \small 
  \bibliography{references}
  
  \normalsize
  \appendix
    \section{Proofs and Extensions of Results from~\cref{range_approx}}
\label{appendix_a}

This appendix establishes the guarantees for $\mathrm{range}(\bA_\mu \bOmega)$ in approximating the action of $\bA_\mu\in \mathbb{R}^{n \times n}$, with $\bOmega\in \mathbb{R}^{n \times l}$ constructed as in~\cref{eq:testmatrix}. These guarantees include and extend those in~\cref{range_approx}.  Some discussed scenarios, such as the range approximation without subspace iteration (see~\cref{thm:standard_gaussian,thm:standard_srht,thm:standard_sparse}), follow from the literature, while others are new.

\subsection{Structural Bounds} \label{appendix_general_bounds}

{We first establish deterministic structural bounds that apply to both randomized and non-randomized orthogonal projectors~$\bPi\in \mathbb{R}^{n \times n}$. These bounds will form the foundation for extending error guarantees of randomized range approximation from settings without subspace iteration to those incorporating it.}

\paragraph{Structural bounds for psd operators.}
We first consider settings with psd operators $\bA_\mu$. 
\begin{proposition}[Deterministic bound for psd operators] \label{thm:det_bound}
Let $\bA$ be psd and $\mu \geq 0$. For any orthogonal projector $\bPi $ and power $q \geq 1$, it holds that
$$
\|(\bI - \bPi) \bA\|^{q+1} \leq \|(\bI - \bPi) \bA_\mu \bA^q\|.
$$
\begin{proof}
Let $\bPi^* = \bI - \bPi$, and let $\bA = \bU \bD \bU^\mathrm{T}$ be the eigenvalue decomposition of $\bA$. Then,
\begin{align*}
\|\bPi^* \bA\|^{q+1} &= \|\bPi^* \bA^2 \bPi^*\|^{\frac{q+1}{2}} = \|(\bU^\mathrm{T} \bPi^* \bU)\, \bD^2\, (\bU^\mathrm{T} \bPi^* \bU)\|^{\frac{q+1}{2}} \\
&\leq \|(\bU^\mathrm{T} \bPi^* \bU)\, \bD^{q+1}\, (\bU^\mathrm{T} \bPi^* \bU)\| = \|\bPi^* \bA^{q+1} \bPi^*\| \\
&\leq \|\bPi^* (\bA^{q+1} + \mu \bA^q) \bPi^*\| = \|\bPi^* \bA_\mu \bA^q \bPi^*\| \leq \|\bPi^* \bA_\mu \bA^q\|.
\end{align*}
In the above, the first and fourth relations use the unitary invariance of the spectral norm. The second relation follows from inequality (8.5) in~\cite[Proposition 8.6]{halko2011finding} and the fact that $\bU^\mathrm{T} \bPi^* \bU$ is an orthogonal projector. The inequality (8.5) in~\cite[Proposition 8.6]{halko2011finding} itself follows from \cite[Theorem IX.2.10]{bhatia1997matrix}.
\end{proof}
\end{proposition}

This yields the following average bound for a randomized orthogonal projector.
\begin{proposition}[Average bound for psd operators] \label{thm:ave_bound}
Let $\bA$ be psd and $\mu \geq 0$, and let $\bPi$ be a random orthogonal projector. Then
$$
\mathbb{E}\left[\|(\bI - \bPi) \bA_\mu\|\right] \leq \mu + \mathbb{E}\left[\|(\bI - \bPi) \bA_\mu \bA^q\|\right]^{\frac{1}{q+1}}.
$$
\begin{proof}
The case $q = 0$ is immediate. For $q\geq 1$, 
    by the Holder's (or Jensen's) inequality, we have:
    \begin{equation*}
        \begin{split}
        \mathbb{E}\left [\|(\bI - \bPi) \bA_\mu\| \right] &\leq \mu + \mathbb{E}\left [\|(\bI - \bPi) \bA\|\right]  \leq \mu + \mathbb{E}\left [\|(\bI - \bPi) \bA\|^{q+1} \right]^{\frac{1}{q+1}}.
        \end{split}
    \end{equation*}
Combining this with~\cref{thm:det_bound} yields the result.
\end{proof}
\end{proposition}

\paragraph{Generalization to indefinite operators.}
The psd setting extends to general symmetric~$\bA_\mu$ as follows.

\begin{proposition}[Deterministic bound for symmetric operators] \label{thm:det_bound_indefinite}
Let $\bA$ be symmetric. For any orthogonal projector $\bPi$ and any $q \geq 1$,
$$
\|(\bI - \bPi) \bA\| \leq 2|\mu| + 2^{\frac{1}{q+1}} \|(\bI - \bPi) \bA_\mu \bA^q\|^{\frac{1}{q+1}}.
$$
\begin{proof}
As in~\cref{thm:det_bound}, we have
$$
\|(\bI - \bPi) \bA\|^{q+1} = \|(\bI - \bPi) \bA^2 (\bI - \bPi)\|^{\frac{q+1}{2}} \leq \|(\bI - \bPi) \bA^{q+1}\|,
$$
and similarly,
$$
\|(\bI - \bPi) \bA^q\|^{q+1} \leq \|(\bI - \bPi) \bA^{q+1}\|^q.
$$
Therefore,
$$
\|(\bI - \bPi) \bA^{q+1}\| \leq \|(\bI - \bPi) \bA_\mu \bA^q\| + |\mu| \|(\bI - \bPi) \bA^q\| \leq \|(\bI - \bPi) \bA_\mu \bA^q\| + |\mu| \|(\bI - \bPi) \bA^{q+1}\|^{\frac{q}{q+1}}.
$$
Set $x = \|(\bI - \bPi) \bA^{q+1}\|^{1/(q+1)}$. Then
$$
x^{q+1} \leq \|(\bI - \bPi) \bA_\mu \bA^q\| + \frac{|\mu|}{x} x^{q+1}.
$$
If $x < 2|\mu|$, the bound follows directly. Otherwise,
$$
x^{q+1} \leq 2 \|(\bI - \bPi) \bA_\mu \bA^q\|,
$$
and the claim follows.
\end{proof}
\end{proposition}

\begin{proposition}[Average bound for symmetric operators] \label{thm:ave_bound_ext}
Let $\bA$ be symmetric and $\bPi$ a random orthogonal projector. Then
$$
\mathbb{E}\left[\|(\bI - \bPi) \bA_\mu\|\right] \leq 3|\mu| + 2^{\frac{1}{q+1}} \mathbb{E}\left[\|(\bI - \bPi) \bA_\mu \bA^q\|\right]^{\frac{1}{q+1}}.
$$
\begin{proof}
The case $q = 0$ is trivial. For $q \geq 1$, we use~\cref{thm:det_bound_indefinite} and Holder's (or Jensen's) inequality:
\begin{align*}
\mathbb{E}[\|(\bI - \bPi) \bA_\mu\|] &\leq |\mu| + \mathbb{E}[\|(\bI - \bPi) \bA\|] \\
&\leq |\mu| + 2|\mu| + 2^{\frac{1}{q+1}} \mathbb{E}[\|(\bI - \bPi) \bA_\mu \bA^q\|^{1/(q+1)}] \\
&\leq 3|\mu| + 2^{\frac{1}{q+1}} \mathbb{E}[\|(\bI - \bPi) \bA_\mu \bA^q\|]^{1/(q+1)}.
\end{align*}
\end{proof}
\end{proposition}

\noindent

\subsection{Gaussian Matrices}

Here we derive extended results for range approximation with Gaussian sketching matrices~$\bX$, and use them to prove~\cref{thm:apriori_gen00,thm:apriori_gen2,thm:apriori2} from the main text.  
Our analysis uses the following theorem from~\cite[Section~10]{halko2011finding}.

\begin{theorem}[Section 10 of~\cite{halko2011finding}] \label{thm:standard_gaussian}
Let $\bC \in \mathbb{R}^{n \times n}$ be a matrix with singular values $\sigma_1 \geq \sigma_2 \geq \dots \geq \sigma_n$, and draw a Gaussian matrix $\bX \in \mathbb{R}^{l \times n}$ with sketching dimension $l \geq 4$. Let $\bPi$ be the orthogonal projector onto the column space of $\bC \bX^\mathrm{T}$. Then, for any rank parameter $3 \leq k \leq l-1$,
\begin{equation} \label{eq:apriori2}
\mathbb{E}[\|(\bI - \bPi) \bC\|] \leq \left(1 + \sqrt{\frac{k+1}{l-k}} \right)\sigma_k(\bC) + \frac{\mathrm{e} \sqrt{l}}{l-k-1} \sqrt{ \sum_{j \geq k} \sigma_j^2(\bC)}.
\end{equation}
Moreover, if $3 \leq k \leq l - 3$, then with probability at least $1 - 6^{-(l-k+1)}$,
\begin{equation} \label{eq:dev_bound}
\|(\bI - \bPi) \bC\| \leq \left(1 + 17 \sqrt{\frac{k+1}{l-k+1}} \right)\sigma_k(\bC) + \frac{8 \sqrt{l}}{l-k} \sqrt{ \sum_{j \geq k} \sigma_j^2(\bC)}.
\end{equation}
\end{theorem}

This result can be directly used to prove~\cref{thm:apriori_gen00}.
\begin{proof}[Proof of~\cref{thm:apriori_gen00}]
    Set $l = 2k+2$ and assume that $l\geq 6$, which implies $3 \leq k \leq l - 3$. According to~\cref{eq:apriori2} from~\cref{thm:standard_gaussian}, we have
    \begin{equation} 
        \mathbb{E}[\|(\bI - \bPi) \bA_\mu\|] \leq \left(1 + \sqrt{\frac{k+1}{k+2}} \right)\sigma_k(\bA_\mu) + \mathrm{e}\sqrt{\frac{{2k+2}}{k+1}} \sqrt{ \sum_{j \geq k} \sigma_j^2(\bA_\mu)} \leq 2 \sigma_k(\bA_\mu) + \sqrt{2}\mathrm{e} \sqrt{ \frac{1}{k+1}\sum_{j \geq k} \sigma_j^2(\bA_\mu)}.
    \end{equation}
    Furthermore, according to~\cref{eq:dev_bound}, it holds with probability at least $1 - 6^{-(2k+2-k+1)} = 1 - 6^{-(k+3)}$,
     \begin{equation} 
        \|(\bI - \bPi) \bA_\mu\| \leq \left(1 + 17 \sqrt{\frac{k+1}{k+3}} \right)\sigma_k(\bA_\mu) + \frac{8 \sqrt{2k+2}}{k+2} \sqrt{ \sum_{j \geq k} \sigma_j^2(\bA_\mu)} \leq  18 \sigma_k(\bA_\mu) + 8 \sqrt{2} \sqrt{ \frac{1}{k+2}\sum_{j \geq k} \sigma_j^2(\bA_\mu)}. 
    \end{equation}
    The lemma then follows immediately from the definition  $\mathrm{sr}_{k}(\bA_\mu^\mathrm{T}\bA_\mu):=    \sum_{j\geq k} \frac{\sigma^2_j(\bA_\mu)}{\sigma^2_k(\bA_\mu)}$.
\end{proof}

\paragraph{Randomized range approximation for psd operators.}
\cref{thm:apriori} extends the average bound in \cref{thm:standard_gaussian} to shifted subspace iteration for psd~$\bA$ via the structural bound~\cref{thm:ave_bound}.

\begin{theorem}[Average bound for Gaussian sketching] \label{thm:apriori}
Let $\bA$ be psd and $\mu \geq 0$. Draw a Gaussian matrix $\bX \in \mathbb{R}^{l \times n}$ with $l \geq 4$, and define the test matrix $\bOmega = \bA^q \bX^\mathrm{T}$. Then, the projection of $\bA_\mu$ onto the range of $\bA_\mu \bOmega$ satisfies
$$
\mathbb{E}[\|(\bI - \bPi) \bA_\mu\|] \leq \mu + \min_{3 \leq k \leq l-1} \left( \left(1 + \sqrt{\frac{k+1}{l-k}} \right) \lambda_k^q (\lambda_k + \mu) + \frac{\mathrm{e} \sqrt{l}}{l-k-1} \sqrt{ \sum_{j \geq k} \lambda_j^{2q} (\lambda_j + \mu)^2 } \right)^{\frac{1}{q+1}}.
$$
\begin{proof}
Apply~\cref{thm:ave_bound} with $\bC = \bA_\mu \bA^q$, and invoke~\cref{eq:apriori2} in~\cref{thm:standard_gaussian}. Note that the singular values of $\bC$ are $\sigma_j(\bC) = \lambda_j^q (\lambda_j + \mu)$.
\end{proof}
\end{theorem}

A deviation bound follows analogously.

\begin{theorem}[Deviation bound for Gaussian sketching] \label{thm:apriori20}
Adopt the setting of~\cref{thm:apriori}. Let  $p \geq 3$ be some integer. Then, with probability at least $1 - 6^{-p}$,
$$
\|(\bI - \bPi) \bA_\mu\| \leq \mu + \min_{3 \leq k \leq l-p} \left( \left(1 + 17 \sqrt{\frac{k+1}{l-k+1}} \right) \lambda_k^q (\lambda_k + \mu) + \frac{8 \sqrt{l}}{l-k} \sqrt{ \sum_{j \geq k} \lambda_j^{2q} (\lambda_j + \mu)^2 } \right)^{\frac{1}{q+1}}.
$$
\begin{proof}
 Choose a rank parameter $3 \leq k \leq l-p$ and an oversampling parameter $p' = l - k + 1$. Notice that $p+1 \leq p'$. 
    By taking $\bC = \bA_\mu$ in~\cref{thm:standard_gaussian} we get the statement of the proposition for $q=0$.   
    By taking $\bC = \bA_\mu \bA^q$ in~\cref{thm:standard_gaussian}, plugging the resulting relation into~\cref{thm:det_bound}, and using the union bound argument we deduce that  
    \begin{equation}\label{eq:dev_bound2}
    \|(\bI - \bPi) \bC \| \leq \mu + \|(\bI - \bPi) \bA \|  \leq \mu + \min_{3 \leq k \leq l-p} \left ( \left  (1+17 \sqrt{\frac{k+1}{l-k+1}} \right ) \sigma_{k}  + \frac{8\sqrt{l}}{l-k} \sqrt{\sum_{j\geq k} \sigma^{2}_j } \right )^\frac{1}{q+1}.  
    \end{equation} 
    holds with probability at least $\sum_{i\geq p+1} 6^{-i} \leq 6^{-p}$.
    The theorem then follows from the fact that the singular values of $\bC$ are $\sigma_j(\bC) =\lambda_{j}(\bA)^q(\lambda_{j}(\bA)+\mu)$, $1 \leq j \leq n$.
\end{proof}
\end{theorem}

Combining \cref{thm:apriori,thm:apriori20} with the definition of the stable rank yields \cref{thm:apriori2} in the main text.

\paragraph{Randomized range approximation for indefinite operators.}

Combining \cref{thm:standard_gaussian} with \cref{thm:det_bound_indefinite,thm:ave_bound_ext} yields a generalization to the shifted subspace iteration with general symmetric operators.

\begin{theorem}[Generalization of~\cref{thm:apriori} to symmetric indefinite $\bA$] \label{thm:apriori3_2}
Let $\bA$ be symmetric. Draw a Gaussian matrix $\bX \in \mathbb{R}^{l \times n}$ with $l \geq 4$, and set $\bOmega = \bA^q \bX^\mathrm{T}$. Then,
$$
\mathbb{E}[\|(\bI - \bPi) \bA_\mu\|] \leq 3|\mu| + 2^{\frac{1}{q+1}} \min_{3 \leq k \leq l-1} \left( \left(1 + \sqrt{\frac{k+1}{l-k}} \right) \sigma_k(\bA_\mu \bA^q) + \frac{\mathrm{e} \sqrt{l}}{l-k-1} \sqrt{ \sum_{j \geq k} \sigma_j^2(\bA_\mu \bA^q) } \right)^{\frac{1}{q+1}}.
$$
\begin{proof}
Combine~\cref{thm:ave_bound_ext} with~\cref{eq:apriori2} in~\cref{thm:standard_gaussian}.
\end{proof}
\end{theorem}

\begin{theorem}[Generalization of~\cref{thm:apriori20} to symmetric indefinite $\bA$] \label{thm:apriori20_ext}
Under the setting of~\cref{thm:apriori3_2}, let $p \geq 3$. Then, with probability at least $1 - 6^{-p}$,
$$
\|(\bI - \bPi) \bA_\mu\| \leq 3|\mu| + 2^{\frac{1}{q+1}} \min_{3 \leq k \leq l-p} \left( \left(1 + 17 \sqrt{\frac{k+1}{l-k+1}} \right) \sigma_k(\bA_\mu \bA^q) + \frac{8 \sqrt{l}}{l-k} \sqrt{ \sum_{j \geq k} \sigma_j^2(\bA_\mu \bA^q) } \right)^{\frac{1}{q+1}}.
$$
\begin{proof}
Apply~\cref{thm:det_bound_indefinite} and the relation~\cref{eq:dev_bound} from~\cref{thm:standard_gaussian}.
\end{proof}
\end{theorem}

\noindent
\cref{thm:apriori3_2,thm:apriori20_ext}, together with the bound
$$
\sigma_k(\bA_\mu \bA^q) \leq (\sigma_k(\bA) + |\mu|)\, \sigma_k(\bA)^q,
$$
directly yield \cref{thm:apriori_gen2}.

\subsection{SRHT Matrices}

In this subsection, we establish randomized range approximation guarantees for SRHT sketching that lead directly to \cref{thm:apriori3} in the main text.  
Our analysis begins with the following bound for the case without subspace iteration, adapted from~\cite[Theorem~2.1]{boutsidis2013improved} and~\cite{tropp2011improved}.

\begin{theorem}[Based on~\cite{boutsidis2013improved}~and~\cite{tropp2011improved}]\label{thm:standard_srht}
Let $\bC \in \mathbb{R}^{n \times n}$ be a matrix with singular values $\sigma_1 \geq \sigma_2 \geq \cdots \geq \sigma_n$. Draw an SRHT matrix $\bX \in \mathbb{R}^{l \times n}$, using zero padding if $n$ is not a power of two. Let $\bPi$ denote the orthogonal projector onto the column space of $\bC \bX^\mathrm{T}$.
Fix failure probability $\delta \in (0,1)$, and suppose the sketch size satisfies $3 \log^2(2n/\delta) \leq l \leq n$. Let $k \geq 3$ be the largest integer such that
$$
18 \left( \sqrt{k} + \sqrt{8 \log(2n/\delta)} \right)^2 \log(k/\delta) < l.
$$
Then, with probability at least $1 - \delta$, the following error bound holds:
\begin{equation} \label{eq:apriori_111}
\|(\bI - \bPi)\bC\| \leq 5 \sigma_{k+1}(\bC) + \sqrt{ \frac{3 \log(n/\delta)}{l} \sum_{j > k} \sigma_j^2(\bC) }.
\end{equation}
\begin{proof}
Let $\bC^* \in \mathbb{R}^{n \times s}$ be the matrix obtained by zero-padding $\bC$ so that the number of columns $s < 2n$ is a power of two. Let $\bX^* \in \mathbb{R}^{l \times s}$ be the corresponding SRHT matrix, and let $\bPi^*$ be the orthogonal projector onto the column space of $\bC^* (\bX^*)^\mathrm{T}$. Since $\bC$ is embedded into $\bC^*$, we have
$$
\|(\bI - \bPi)\bC\| = \|(\bI - \bPi^*)\bC^*\|.
$$

By~\cite[Theorem 2.1]{boutsidis2013improved}, with probability at least $1 - \delta$,
\begin{align*}
\|(\bI - \bPi^*) \bC^*\| &\leq \left(4 + \sqrt{ \frac{3}{l} \log(s/\delta) \log(n/\delta) } \right) \sigma_{k+1}(\bC) + \sqrt{ \frac{3 \log(n/\delta)}{l} \sum_{j > k} \sigma_j^2(\bC) }.
\end{align*}

Since $\mathrm{rank}(\bC) = n$ and $s < 2n$, we have $\log(s/\delta) \leq \log(2n/\delta)$, and thus
$
\sqrt{ \frac{3}{l} \log(s/\delta) \log(n/\delta) } \leq 1,
$
for the specified values of $l$. This yields the simplified bound~\eqref{eq:apriori_111}.
\end{proof}
\end{theorem}

\paragraph{Randomized range approximation for psd operators.}
The following result extends \cref{thm:standard_srht} to shifted subspace iteration for psd~$\bA$.

\begin{theorem}[Deviation bound for SRHT sketching] \label{thm:srht_dev}
Let $\bA$ be psd and $\mu \geq 0$. Draw an SRHT matrix $\bX \in \mathbb{R}^{l \times n}$ with $3 \log^2(2n/\delta) \leq l \leq n$, where $0 < \delta < 1$ is the failure probability. Let the test matrix be $\bOmega = \bA^q \bX^\mathrm{T}$. Let $k \geq 3$ be the largest integer such that
$$
18 \left( \sqrt{k} + \sqrt{8 \log(2n/\delta)} \right)^2 \log(k/\delta) < l.
$$
Then, with probability at least $1 - \delta$,
$$
\|(\bI - \bPi) \bA_\mu\| \leq \mu + \left( 5 \lambda_{k+1}^q (\lambda_{k+1} + \mu) + \sqrt{ \frac{3 \log(n/\delta)}{l} \sum_{j>k} \lambda_j^{2q} (\lambda_j + \mu)^2 } \right)^{\frac{1}{q+1}}.
$$
\begin{proof}
Apply~\cref{thm:det_bound} with $\bC = \bA_\mu \bA^q$, and use~\cref{thm:standard_srht}. Note that the singular values of $\bA_\mu \bA^q$ are $\lambda_j^q(\lambda_j + \mu)$.
\end{proof}
\end{theorem}

\paragraph{Randomized range approximation for indefinite operators.}
Finally, we combine \cref{thm:standard_srht} with \cref{thm:det_bound_indefinite} to obtain a generalization to the subspace iteration bounds for general shifted symmetric operators.

\begin{theorem}[Generalization of~\cref{thm:srht_dev} to symmetric indefinite $\bA$]\label{thm:srht_dev_ind}
Let $\bA$ be symmetric. Let $l$, $k$, and $\delta$ satisfy the assumptions of~\cref{thm:srht_dev}, and draw an SRHT matrix $\bX \in \mathbb{R}^{l \times n}$. Then, with probability at least $1 - \delta$,
$$
\|(\bI - \bPi) \bA_\mu\| \leq 3|\mu| + 2^{\frac{1}{q+1}} \left( 5 \sigma_{k+1}(\bA_\mu \bA^q) + \sqrt{ \frac{3 \log(n/\delta)}{l} \sum_{j>k} \sigma_j^2(\bA_\mu \bA^q) } \right)^{\frac{1}{q+1}}.
$$
\begin{proof}
Apply~\cref{thm:det_bound_indefinite} and use~\cref{thm:standard_srht} with $\bC = \bA_\mu \bA^q$.
\end{proof}
\end{theorem}

The results \cref{thm:standard_srht,thm:srht_dev,thm:srht_dev_ind} together imply \cref{thm:apriori3} in the main text.

\subsection{Sparse Embeddings}

We now establish range approximation guarantees for sparse embedding matrices, which lead directly to \cref{thm:apriori4} in the main text.  
The analysis parallels the Gaussian and SRHT cases and begins with the following result, adapted from~\cite[Theorem~8]{cohen2015optimal} and~\cite[Proof of Lemma~4.2]{derezinski2024faster}.

\begin{theorem}[Based on~\cite{cohen2015optimal}] \label{thm:standard_sparse}
Let $\bC \in \mathbb{R}^{n \times n}$ be a matrix with singular values $\sigma_1 \geq \sigma_2 \geq \cdots \geq \sigma_n$. Fix failure probability $0 < \delta < 1$ and target rank $k$.

Draw a sparse embedding matrix $\bX \in \mathbb{R}^{l \times n}$, with embedding dimension $l \geq c k \log(k/\delta)$ and $\gamma \geq c \log(k/\delta)$ nonzero entries per column, for some universal constant $c$. Let $\bPi$ be the orthogonal projector onto the column space of $\bC \bX^\mathrm{T}$. Then, with probability at least $1 - \delta$, the following bound holds:
\begin{equation} \label{eq:apriori_sparse}
\|(\bI - \bPi)\bC\| \leq \frac{3}{2} \sigma_{k+1}(\bC) + \sqrt{ \frac{1}{2k} \sum_{j > k} \sigma_j^2(\bC) }.
\end{equation}
\end{theorem}
\paragraph{Randomized range approximation for psd operators} 

\begin{theorem}(Deviation bound for sparse sketching)\label{thm:sparse_dev}
   Assume that $\bA$ is psd  and $\mu \geq 0$.  Set rank parameter $k$. 
    Draw a sparse randomized embedding $\bX \in \mathbb{R}^{l \times n}$ using sketching dimension $l \geq c k \log(k/\delta) $ and  $\gamma \geq c \log(k/\delta) $ non-zero entries per column, where $c$ is some universal constant.  
    Let the test matrix be $\bOmega =  \bA^q \bX^\mathrm{T}$.   Then the projection of $\bA_\mu$ onto the range of $\bA_\mu \bOmega$ satisfies the following error bound
    \begin{align} \label{eq:apriori}
    \| (\bI - \bPi) \bA_\mu\| \leq \mu + \left( \frac{3}{2} \lambda_{k+1}^q (\lambda_{k+1}+\mu) + \sqrt{\frac{1}{2k}  \sum_{j>k} \lambda_j^{2q} (\lambda_{j}+\mu)^2 } \right)^{\frac{1}{q+1}} 
    \end{align}
    with probability at least $1- \delta$. 
    \begin{proof}
        The theorem follows by combining~\cref{thm:det_bound} with~\cref{eq:apriori_sparse} in~\cref{thm:standard_sparse}, and using the fact that the singular values of $\bA_\mu \bA^q$ are $\lambda_{j}^q(\lambda_{j}+\mu)$, $1 \leq j \leq n$.
    \end{proof}
\end{theorem}

\paragraph{Randomized range approximation for indefinite operators}

\begin{theorem}(Deviation bound for sparse sketching)\label{thm:sparse_dev2}
   Assume that $\bA$ is symmetric.  Set rank parameter $k$. 
    Draw a sparse randomized embedding $\bX \in \mathbb{R}^{l \times n}$ using sketching dimension $l \geq c k \log(k/\delta) $ and  $\gamma \geq c \log(k/\delta) $ non-zero entries per column, where $c$ is some universal constant.  
    Let the test matrix be $\bOmega =  \bA^q \bX^\mathrm{T}$.   Then the projection of $\bA_\mu$ onto the range of $\bA_\mu \bOmega$ satisfies the following error bound
    \begin{align} 
    \| (\bI - \bPi) \bA_\mu\| \leq 3|\mu| + 2^{\frac{1}{q+1}}\left(  \frac{3}{2}  \sigma_{k+1}(\bA_\mu \bA^q) + \sqrt{\frac{1}{2k}  \sum_{j>k} \sigma_{j}(\bA_\mu \bA^q)^2 } \right)^{\frac{1}{q+1}} 
    \end{align}
    with probability at least $1- \delta$. 
    \begin{proof}
         The theorem follows by combining~\cref{thm:det_bound_indefinite} with~\cref{eq:apriori_sparse} in~\cref{thm:standard_sparse}.
    \end{proof}
\end{theorem}

\cref{thm:standard_sparse,thm:sparse_dev,thm:sparse_dev2} directly yield \cref{thm:apriori4} in the main text.

    \section{Proofs of Propositions from~\cref{sec:allprecond}}
\label{appendix_b}
This appendix presents the proofs of the condition number bounds from~\cref{sec:allprecond} for R-RandRAND, C-RandRAND and G-RandRAND  preconditioners. 
\paragraph{R-RandRAND.}
\begin{proof}[Proof of~\cref{thm:cnd_bound}]
Let $\bx \in \mathbb{R}^n$ be arbitrary vector. By definition of the preconditioner, we have:
\begin{equation} \label{eq:cnd_bound1}
\begin{split}
\langle \bx, \bA_\mu \bP \bx \rangle &= \langle \bx, (\bI - \bPi)\bA_\mu(\bI - \bPi)\bx \rangle + \tau \langle \bx, \bPi \bx \rangle \\
&= \| \bA_\mu^{1/2} (\bI - \bPi) \bx \|^2 + \tau \| \bPi \bx \|^2.
\end{split}
\end{equation}

We now derive an upper bound. Note that
\begin{equation}
\begin{split}
\langle \bx, \bA_\mu \bP \bx \rangle &\leq \| \bA_\mu^{1/2} (\bI - \bPi) \|^2 \| (\bI - \bPi) \bx \|^2 + \tau \| \bPi \bx \|^2 \\
&= \|\bE\| \| (\bI - \bPi)\bx \|^2 + \tau \| \bPi \bx \|^2 \\
&\leq \max(\tau, \|\bE\|) \left( \| (\bI - \bPi) \bx \|^2 + \| \bPi \bx \|^2 \right) = \max(\tau, \|\bE\|) \|\bx\|^2,
\end{split}
\end{equation}
where we used the fact that $\|\bE\| = \| \bA_\mu^{1/2} (\bI - \bPi) \|^2$.

Taking the maximum over unit vectors $\bx$, we conclude that
$$
\|\bA_\mu \bP\| = \max_{\bx \neq 0} \frac{\langle \bx, \bA_\mu \bP \bx \rangle}{\|\bx\|^2} \leq \max(\tau, \|\bE\|).
$$

For the lower bound, we return to~\eqref{eq:cnd_bound1} and use
$$
\| \bA_\mu^{1/2} (\bI - \bPi) \bx \|^2 \geq \lambda_{\min}(\bA_\mu) \| (\bI - \bPi) \bx \|^2.
$$
Hence,
\begin{equation}
\begin{split}
\langle \bx, \bA_\mu \bP \bx \rangle &\geq \lambda_{\min}(\bA_\mu) \| (\bI - \bPi) \bx \|^2 + \tau \| \bPi \bx \|^2 \\
&\geq \min\big( \tau, \lambda_{\min}(\bA_\mu) \big) \left( \| (\bI - \bPi) \bx \|^2 + \| \bPi \bx \|^2 \right) = \min(\tau, \lambda_{\min}(\bA_\mu)) \|\bx\|^2.
\end{split}
\end{equation}
It follows that
$$
\lambda_{\min}(\bA_\mu \bP) = \min_{\bx \neq 0} \frac{\langle \bx, \bA_\mu \bP \bx \rangle}{\|\bx\|^2} \geq \min\big( \tau, \lambda_{\min}(\bA_\mu) \big).
$$

Combining the bounds on the largest and smallest eigenvalues gives
$$
\cond(\bA_\mu \bP) \leq \frac{\max(\tau, \|\bE\|)}{\min(\tau, \lambda_{\min}(\bA_\mu))},
$$
which completes the proof.
\end{proof}

\begin{proof}[Proof of~\cref{thm:cnd_bound_spec}]
Let $\bC = \bA_\mu \bA^q$, and assume that $ck$ is an integer. By~\cref{thm:apriori2}, with probability at least $1 - 6^{-ck} \geq 1 - 6^{-k}$, it holds that
\begin{equation} \label{eq:cnd_bound_spec_1}
\| (\bI - \bPi) \bA_\mu \| \leq \mu + \left( \left( C_1 + C_2 \sqrt{\frac{1}{ck} \mathrm{sr}_{ck}(\bC^2)} \right) \lambda_{ck}(\bC) \right)^{\frac{1}{q+1}},
\end{equation}
where $C_1 = 18$, $C_2 = 8\sqrt{2}$ for the deviation bound, and $C_1 = 2$, $C_2 = \sqrt{2}\mathrm{e}$ in expectation.

We begin with a simple inequality:
$$
\frac{1}{(c - 1)k} \sum_{k \leq j < ck} \lambda_j(\bC)^2 \geq \lambda_{ck}(\bC)^2.
$$
Hence,
$$
\frac{1}{(c - 1)k} \sum_{j \geq k} \lambda_j(\bC)^2 \geq \frac{1}{ck} \sum_{j \geq ck} \lambda_j(\bC)^2 + \lambda_{ck}(\bC)^2.
$$
Taking square roots, we obtain
$$
\sqrt{ \frac{2}{(c - 1)k} \sum_{j \geq k} \lambda_j(\bC)^2 } \geq \lambda_{ck}(\bC) + \sqrt{ \frac{1}{ck} \sum_{j \geq ck} \lambda_j(\bC)^2 }.
$$

This yields the following inequality in terms of stable ranks:
\begin{equation} \label{eq:sr_lambda}
\sqrt{ \frac{2}{(c - 1)k} \mathrm{sr}_{k}(\bC^2) }  \lambda_k(\bC) \geq \left( 1 + \sqrt{ \frac{1}{ck} \mathrm{sr}_{ck}(\bC^2) } \right) \lambda_{ck}(\bC).
\end{equation}

Substituting~\eqref{eq:sr_lambda} into~\eqref{eq:cnd_bound_spec_1}, we obtain:
\begin{equation} \label{eq:cnd_bound_spec_2}
\begin{split}
\| (\bI - \bPi) \bA_\mu \| 
&\leq \mu + \left( \left( C_1 + C_2 \sqrt{ \frac{1}{ck} \mathrm{sr}_{ck}(\bC^2) } \right) \lambda_{ck}(\bC) \right)^{\frac{1}{q+1}} \\
&\leq \mu + \left( \max(C_1, C_2) \sqrt{ \frac{2}{(c - 1)k} \mathrm{sr}_{k}(\bC^2) }  \lambda_k(\bC) \right)^{\frac{1}{q+1}} \\
&\leq \mu + \left( \sqrt{ \frac{C_3}{(c - 1)k} \mathrm{sr}_{k}(\bC^2) } \lambda_k(\bC) \right)^{\frac{1}{q+1}},
\end{split}
\end{equation}
where $C_3 = 2 C_1^2 = 2 \cdot 18^2$ in the deviation case, and $C_3 = 4 \mathrm{e}^2$ in expectation.

We now combine~\eqref{eq:cnd_bound_spec_2} with the bound from~\cref{thm:cnd_bound}:
\begin{equation} \label{eq:cnd_bound_spec_3}
\mathrm{cond}(\bA_\mu \bP) 
\leq  \frac{\| \bE \|}{\lambda_{\min}(\bA_\mu)} 
\leq   \frac{ \| (\bI - \bPi) \bA_\mu \| }{ \lambda_{\min}(\bA_\mu) } 
\leq 1 +  \left( \frac{C_3}{(c - 1)k}  \mathrm{sr}_k(\bC^2)  \frac{ \lambda_k^2(\bC) }{ \lambda_{\min}(\bA_\mu)^{2q + 2} } \right)^{\frac{1}{2(q+1)}}.
\end{equation}

Now recall that
$$
\mathrm{sr}_k(\bC^2) = (n - k + 1) \frac{ \lambda_{\min}(\bC^2) + \alpha }{ \lambda_k(\bC^2) }  \mathrm{cond}_{k, \alpha}(\bC^2),
$$
which gives the result for $q = 0$ when we set $\alpha = 0$.

For $q \geq 1$, we set $\alpha = \mu^{2q+2}$ and use the inequality
$$
\lambda_{\min}(\bA_\mu)^{2q+2} = \lambda_{\min}\big( \bA_\mu^{2q} (\bA + \mu \bI)^2 \big) \geq \lambda_{\min}(\bA_\mu^{2q} \bA^2) + \mu^2 \lambda_{\min}(\bA_\mu^{2q}) \geq \lambda_{\min}(\bC^2) + \mu^{2q+2}.
$$
Substituting this into~\eqref{eq:cnd_bound_spec_3} completes the proof.
\end{proof}

\paragraph{C-RandRAND.}   
\begin{proof}[Proof of~\cref{thm:cnd_bound2}]
We begin by bounding $\|\bP^{1/2} \bA_\mu \bP^{1/2}\|$. Using the identity $\|\bP^{1/2} \bA_\mu \bP^{1/2}\| = \|\bA_\mu^{1/2} \bP \bA_\mu^{1/2}\|$, and applying the triangle inequality along with the fact that similarity transforms preserve eigenvalues, we obtain:
\begin{equation} \label{eq:crandrand1}
\begin{split}
    \|\bP^{1/2} \bA_\mu \bP^{1/2}\|
    &= \|\bA_\mu^{1/2}(\bI - \bPi + \tau \bPi \bA_\mu^{-1} \bPi)\bA_\mu^{1/2}\| \\
    &\leq \|\bA_\mu^{1/2} (\bI - \bPi) \bA_\mu^{1/2}\| + \tau \|\bA_\mu^{1/2} \bPi \bA_\mu^{-1} \bPi \bA_\mu^{1/2}\| \\
    &= \|\bE\| + \tau\, \lambda_{\max}(\bPi \bA_\mu \bPi \bA_\mu^{-1} \bPi).
\end{split}
\end{equation}

Next, we bound the norm of the inverse. Observe that
\begin{equation} \label{eq:crandrand2}
\begin{split}
    \lambda_{\min}(\bP^{1/2} \bA_\mu \bP^{1/2})^{-1}
    &= \|\bP^{-1/2} \bA_\mu^{-1} \bP^{-1/2}\| 
    = \|\bA_\mu^{-1/2} \bP^{-1} \bA_\mu^{-1/2}\| \\
    &\leq \|\bA_\mu^{-1/2} (\bI - \bPi) \bA_\mu^{-1/2}\| + \tau^{-1} \|\bA_\mu^{-1/2} (\bPi \bA_\mu^{-1} \bPi)^{-1} \bA_\mu^{-1/2}\| \\
    &= \|\bF\| + \tau^{-1} \| (\bPi \bA_\mu^{-1} \bPi)^{-1/2} (\bPi \bA_\mu^{-1} \bPi) (\bPi \bA_\mu^{-1} \bPi)^{-1/2} \| \\
    &= \|\bF\| + \tau^{-1}.
\end{split}
\end{equation}

The bounds~\cref{eq:crandrand1,eq:crandrand2} yield the bound~\cref{eq:c_randrand_eigenvalues} for the eigenvalues of $\bP^{1/2} \bA_\mu \bP^{1/2}$.

We now turn to obtaining the bound~\cref{eq:cnd_bound2} for the condition number:
\begin{align*}
    \mathrm{cond}(\bP^{1/2} \bA_\mu \bP^{1/2}) 
    &= \|\bP^{1/2} \bA_\mu \bP^{1/2}\|  \|\bP^{-1/2} \bA_\mu^{-1} \bP^{-1/2}\| \\
    &\leq \left( \|\bE\| + \tau\, \lambda_{\max}(\bPi \bA_\mu \bPi \bA_\mu^{-1} \bPi) \right)\left( \|\bF\| + \tau^{-1} \right) \\
    &= \|\bE\|\|\bF\| + \rho \|\bE\|\|\bF\| + \rho^{-1} \lambda_{\max}(\bPi \bA_\mu \bPi \bA_\mu^{-1} \bPi) + \lambda_{\max}(\bPi \bA_\mu \bPi \bA_\mu^{-1} \bPi) \\
    &\leq (1 + \rho)\, \frac{\|\bE\|}{\lambda_{\min}(\bA_\mu)} + (1 + \rho^{-1})\, \lambda_{\max}(\bPi \bA_\mu \bPi \bA_\mu^{-1} \bPi),
\end{align*}
where we used $\|\bF\| \leq \lambda_{\min}^{-1}(\bA_\mu)$. To complete the proof, we upper-bound the last term:
\begin{align*}
    \lambda_{\max}(\bPi \bA_\mu \bPi \bA_\mu^{-1} \bPi)
    &= \lambda_{\max}(\bA_\mu^{-1/2} \bPi \bA_\mu \bPi \bA_\mu^{-1/2}) 
    = \|\bA_\mu^{-1/2} \bPi \bA_\mu^{1/2}\|^2 \\
    &\leq \left(1 + \|\bA_\mu^{-1/2}(\bI - \bPi)\bA_\mu^{1/2}\|\right)^2
    \leq \left(1 + \|\bF\|^{1/2} \|\bE\|^{1/2} \right)^2.
\end{align*}
\end{proof}

\paragraph{G-RandRAND.}   
\begin{proof}[Proof of~\cref{thm:cnd_bound3}]
Let $\bM := \bA_\mu \bA_\mu^\mathrm{T}$. The result follows directly from~\cref{thm:cnd_bound2} by observing that $\bP^2$ corresponds to the C-RandRAND preconditioner~\cref{eq:c_randrand_def} for the operator $\bM$ (however the G-RandRAND preconditioner for $\bA_\mu$ may use a different approximation basis than C-RandRAND for $\bM$).

Specifically, note that
$$
\tau^2 = \rho \|(\bI - \bPi)\bM(\bI - \bPi)\| \|\bA_\mu^{-1} \bPi \bM \bPi \bA_\mu^{-\mathrm{T}}\|^{-1} = \rho \|(\bI - \bPi)\bM(\bI - \bPi)\| \lambda_{\mathrm{max}}(\bPi \bM \bPi \bM^{-1} \bPi).
$$
Therefore, by replacing $\bA_\mu$ with $\bM$ and $\tau$ with $\tau^2$ in the quasi-optimality guarantees from~\cref{thm:cnd_bound2}, the desired bounds follow. Here we used standard identities 
$$
\sigma_i(\bP \bA_\mu)^2 = \lambda_i(\bP \bM \bP) \quad \text{and} \quad \|\bA_\mu^{-1} \bPi \bA_\mu\|^2 = \lambda_{\mathrm{max}}(\bPi \bM \bPi \bM^{-1} \bPi).
$$
\end{proof}

\begin{proof}[Proof of~\cref{thm:cnd_bound_spec_ext}]
Define $\bC = \bA_\mu \bA^q$, and let $\mu_{\mathrm{rel}} := |\mu| / \sigma_{\min}(\bA_\mu)$. Assume that $ck$ is an integer. According to~\cref{thm:apriori3,thm:apriori20_ext}, it holds with probability at least $1 - 6^{-ck} \geq 1 - 6^{-k}$ that
\begin{equation} \label{eq:cnd_bound_spec_1_2}
\| (\bI - \bPi) \bA_\mu \| \leq 3|\mu| + \left( \left(2C_1 + 2C_2 \sqrt{\frac{1}{ck} \mathrm{sr}_{ck}(\bC^2)}\right) \sigma_{ck}(\bC) \right)^{\frac{1}{q+1}},
\end{equation}
where $C_1 = 18$, $C_2 = 8\sqrt{2}$; or, in expectation, $C_1 = 2$, $C_2 = \sqrt{2}e$.

Furthermore, it holds that
$$
\sqrt{\frac{2}{(c-1)k} \, \mathrm{sr}_{k}(\bC^2)} \sigma_{k}(\bC) \geq \left(1 + \sqrt{\frac{1}{ck} \, \mathrm{sr}_{ck}(\bC^2)}\right) \sigma_{ck}(\bC),
$$
which can be shown analogously to~\cref{eq:sr_lambda} from the proof of~\cref{thm:cnd_bound_spec}.

Substituting this into~\cref{eq:cnd_bound_spec_1_2}, we obtain with probability at least $1 - 6^{-k}$
\begin{equation} \label{eq:cnd_bound_spec_2_2}
\begin{split}
\|(\bI - \bPi)\bA_\mu\| 
&\leq 3|\mu| + \left(2(C_1 + C_2 \sqrt{\tfrac{1}{ck} \mathrm{sr}_{ck}(\bC^2)}) \, \sigma_{ck}(\bC)\right)^{\frac{1}{q+1}} \\
&\leq 3|\mu| + \left(2\max(C_1, C_2) \sqrt{\tfrac{2}{(c-1)k} \mathrm{sr}_{k}(\bC^2)} \, \sigma_{k}(\bC)\right)^{\frac{1}{q+1}} \\
&\leq 3|\mu| + \left(2 \sqrt{\tfrac{C_3}{(c-1)k} \mathrm{sr}_{k}(\bC^2)} \, \sigma_{k}(\bC)\right)^{\frac{1}{q+1}},
\end{split}
\end{equation}
where $C_3 = 2 \cdot 18^2$ in the deviation bound, and $C_3 = 4e^2$ in expectation.

According to~\cref{thm:cnd_bound3}, for $\rho = 1$, the condition number of $\bP\bA_\mu$ satisfies:
$$
\cond(\bP\bA_\mu)^2 \leq 2 \left(\frac{\|(\bI - \bPi)\bA_\mu\|}{\sigma_{\min}(\bA_\mu)}\right)^2 + 2 \left(1 + \frac{\|(\bI - \bPi)\bA_\mu\|}{\sigma_{\min}(\bA_\mu)}\right)^2 
\leq \left(2 + 2 \frac{\|(\bI - \bPi)\bA_\mu\|}{\sigma_{\min}(\bA_\mu)}\right)^2.
$$

We now consider two cases. For $q = 0$, we set  $\bC = \bA_\mu$ and $\alpha = 0$. Using the identity
\begin{equation}\label{eq:sr_cond_k}
\mathrm{sr}_k(\bC^2) = (n - k + 1) \frac{ \sigma_{\min}(\bC^2) + \alpha }{ \sigma_k(\bC^2) } \mathrm{cond}_{k,\alpha}(\bC^2),
\end{equation}
we obtain:
$$
\frac{\|(\bI - \bPi)\bA_\mu\|}{\sigma_{\min}(\bA_\mu)} 
\leq 3 \mu_{\mathrm{rel}} + 2 \sqrt{ \frac{C_3}{(c-1)k}  \mathrm{sr}_k(\bA_\mu^2) \frac{ \sigma_k(\bA_\mu^2) }{ \sigma_{\mathrm{min}}(\bA_\mu^2) } }
\leq 3 \mu_{\mathrm{rel}} + 2 \sqrt{ \frac{n C_3}{(c-1)k} \cond_k(\bA_\mu^2) }.
$$

For $q \geq 1$, we set $\bC = \bA_\mu \bA^q$ and  $\alpha = |\mu|^{2q+2}$ in~\cref{eq:sr_cond_k} and obtain
\begin{align*}
\frac{\|(\bI - \bPi)\bA_\mu\|}{\sigma_{\min}(\bA_\mu)} 
&\leq 3 \mu_{\mathrm{rel}} 
+ \frac{(\sigma_{\min}(\bC^{2}) + |\mu|^{2q+2})^{1/(2q+2)}}{\sigma_{\min}(\bA_\mu)} \left(2 \sqrt{ \frac{C_3}{(c-1)k} \mathrm{sr}_k(\bC^{2}) \frac{ \sigma_k(\bC^{2}) }{ \sigma_{\min}(\bC^{2}) + |\mu|^{2q+2} } } \right)^{1/(q+1)} \\
&\leq 3 \mu_{\mathrm{rel}} + (1 + \mu_{\mathrm{rel}})  \left( 2 \sqrt{ \frac{n C_3}{(c-1)k} \cond_{k, |\mu|^{2q+2}}(\bC^{2}) } \right)^{1/(q+1)},
\end{align*}
where the last step uses the bound:
\begin{align*}
\frac{(\sigma_{\min}(\bA_\mu^2 \bA^{2q}) + |\mu|^{2q+2})^{1/(2q+2)}}{\sigma_{\min}(\bA_\mu)} 
&\leq \max_j \left( \left( \frac{ |\lambda_j(\bA) + \mu|^2  |\lambda_j(\bA)|^{2q} + |\mu|^{2q+2} }{ |\lambda_j(\bA) + \mu|^{2q+2} } \right)^{1/(2q+2)} \right) \\
&= \max_j \left( \left( \left(\frac{|\lambda_j(\bA)|}{|\lambda_j(\bA) + \mu|}\right)^{2q} + \left(\frac{|\mu|}{|\lambda_j(\bA) + \mu|}\right)^{2q+2} \right)^{1/(2q+2)} \right) \\
&\leq \max_j \left(1 + \frac{|\mu|}{|\lambda_j(\bA) + \mu|}\right) = (1+\mu_{\mathrm{rel}}).
\end{align*}
\end{proof}

    \section{Extended Discussion of Numerical Experiments}
\label{appendix_c}
This appendix contains additional details on the experimental setup and implementation of the methods.

\subsection{Description of the Portfolio Optimization Problem}

We consider covariance matrix of the form $\bSigma = \bF\bF^\mathrm{T} + \bD$, where $\bF \in \mathbb{R}^{m \times s}$ is the loading factor and $\bD \in \mathbb{R}^{m \times m}$ is the diagonal factor. To construct the factor matrices, we first generate a synthetic ground-truth covariance $\bSigma_{\mathrm{true}} \in \mathbb{R}^{m \times m}$ whose eigenvectors are the columns of an orthogonalized Gaussian matrix, the first $20$ eigenvalues decay exponentially, and the remaining eigenvalues follow a power-law decay $i^{-1}$, where $i$ is the eigenvalue index. Then we approximate $\bSigma_{\mathrm{true}}$ with a rank-$s$ Nystr{\"o}m approximation $\bF\bF^\mathrm{T}$ (with $s = 100$), and define $\bD = \operatorname{diag}\big(\bSigma_{\mathrm{true}} - \bF\bF^\mathrm{T}\big)$.
The inequality constraint matrix $\bM \in \mathbb{R}^{p \times m}$ (with $p = 6250$) is given by $5 \cdot \bt \mathbf{1}^\mathrm{T}+ 10^{-3} \cdot \bT$, where $\bt$ is a standard Gaussian vector and $\bT$ a standard Gaussian matrix. The cost vector $\bm \in \mathbb{R}^m$ is also taken to be standard Gaussian, and we set $\gamma = 0.5$.
The matrix $\bZ$ in~\cref{eq:newton_system} is the transpose of the $(p + s + 1) \times (m + p + s + 1)$ equality-constraint matrix that arises after rewriting~\cref{eq:portfolio} in equality-and-box-constrained form. We set the tolerances for relative primal–dual feasibility and the optimality gap to $10^{-9}$.

\subsection{Description of the PDE-Constrained Inverse Problem}
We here provide a detailed overview of the PDE-constrained inverse problem used in~\cref{sec:basis-less-pde-appn}. It represents a source identification problem posed as the least-squares minimization
$$
   \min_{f,\,u} \, \frac{1}{2} \|Q u - z\|_{L^2(\Omega)}^2 \;+\; \frac{\mu}{2}\,\|f\|_{L^2(\Omega)}^2
   \quad
   \text{subject to}
   \quad \mathcal{L} u = f,
$$
where $\mu \ge 0$ is a regularization parameter, $\mathcal{L}$ is the PDE operator, $Q$ is the observation operator, and $z = Qu^* + \epsilon$ is the observation with noise $\epsilon$.  

\begin{figure}[ht!] 
    \centering
    \hspace{2em}
    \begin{subfigure}{0.25\textwidth}
        \includegraphics[page=2,width=\linewidth]{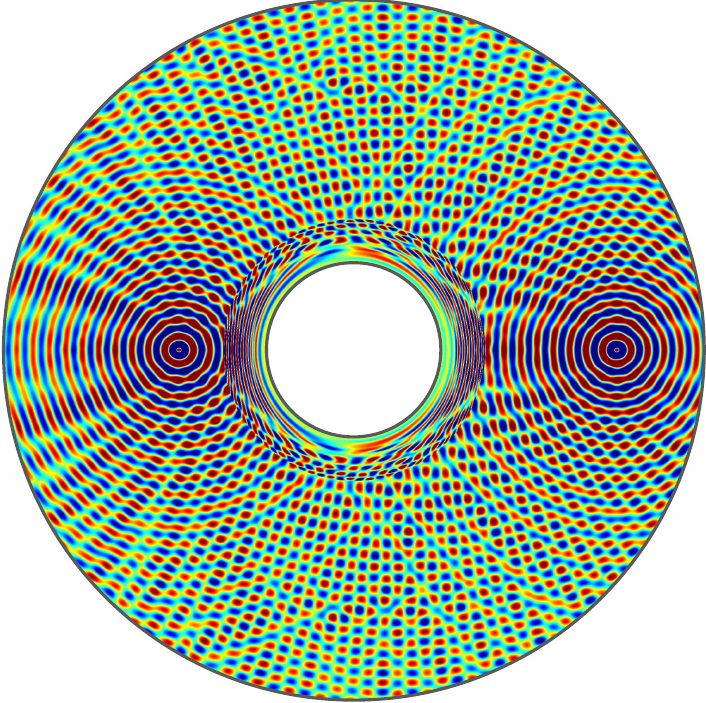}
        \caption{Exact field}
        \label{fig:invcloak1}
    \end{subfigure}
    \hspace{2em}
    \begin{subfigure}{0.25\textwidth}
        \includegraphics[page=3,width=\linewidth]{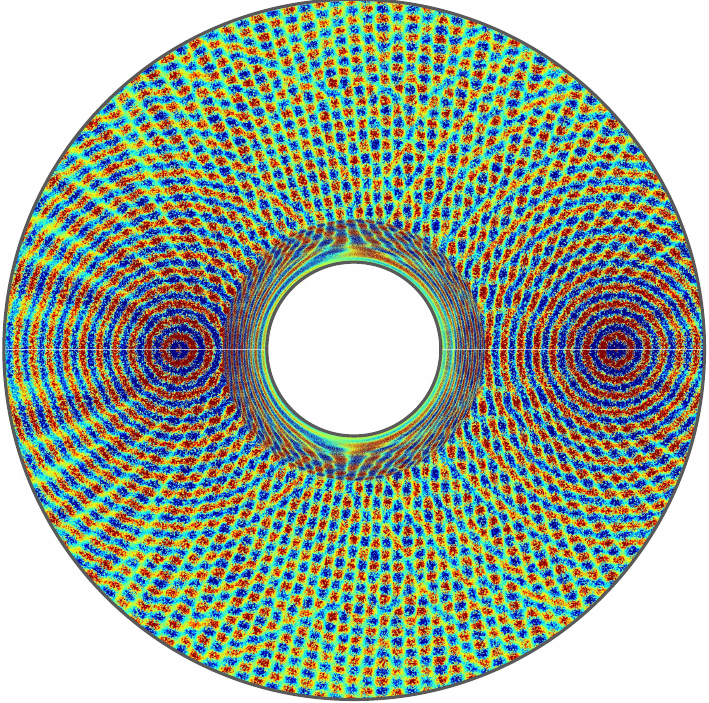}
        \caption{Observed noisy field}
        \label{fig:invcloak2}
    \end{subfigure}
   \hspace{2em}
    \begin{subfigure}{0.25\textwidth}
        \includegraphics[page=4,width=\linewidth]{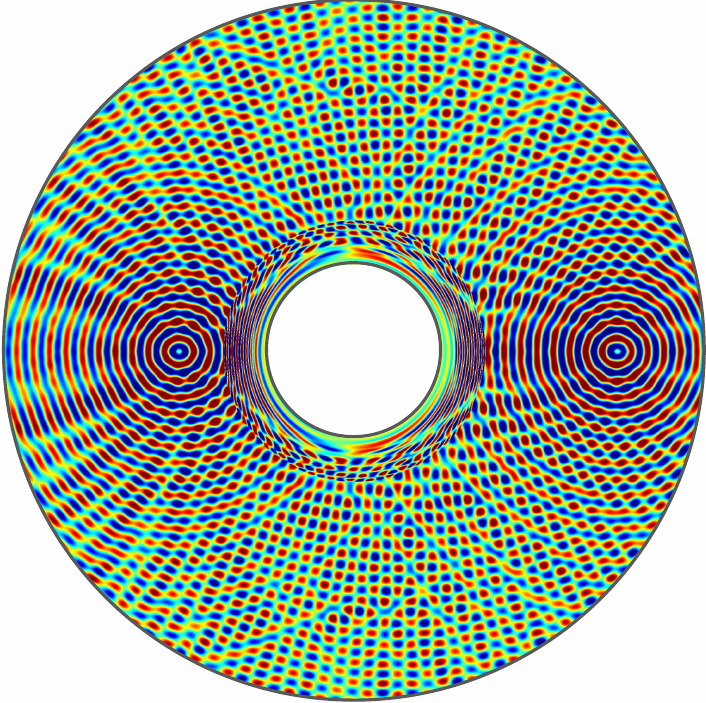}
        \caption{Recovered field}
        \label{fig:invcloak3}
    \end{subfigure}
    \hspace{2em}
    \caption{Invisibility Cloak Benchmark. The real component of pressure field $u$.}
    \label{fig:invcloak00}
\end{figure}

As a model problem, we consider the invisibility cloak benchmark from~\cite{balabanov2021randomized}, which describes a two-dimensional acoustic wave scattering, in which the scatterer is covered by a cloak composed of alternating layers of mercury and isotropic light liquids. The cloak is designed to reduce the scatterer’s visibility within some frequency band. The governing operator is the inhomogeneous Helmholtz operator (excluding absorbing boundary conditions), given by
\begin{equation}\label{eq:pde}
    \mathcal{L} u = \nabla \cdot(\rho^{-1} \nabla u) + \rho^{-1} \kappa^2 u,
\end{equation}
and we set $Q$ to be the identity. Two unit sources are placed outside the cloak and must be recovered from noisy measurements of the pressure field, $u + \epsilon$ (see~\cref{fig:invcloak1,fig:invcloak2}). We set $\epsilon(x) = 0.25\,\omega(x)\,\lvert u(x)\rvert,$ where $\omega(x)$ is a standard normal random variable.

The domain is discretized using approximately $2.5 \cdot 10^5$ second-order triangular finite elements, resulting in $n \approx 5 \cdot 10^5$ degrees of freedom. After eliminating the state variable $u$, the inverse problem reduces to the following regularized least-squares problem:
\begin{equation}\label{eq:discrete_pde_app}
    \min_{\bff}\;\; \tfrac12 \bigl(\bK^{-1}\bM\,\bff - \mathbf{z}\bigr)^\mathrm{T}\,\bM\,\bigl(\bK^{-1}\bM\,\bff - \mathbf{z}\bigr) \;+\; \tfrac{\mu}{2} \,\bff^{\mathrm{T}} \bM\,\bff,
\end{equation}
where $\bM \in \mathbb{R}^{n \times n}$ is the spd mass matrix, and $\bK \in \mathbb{C}^{n \times n}$ is the Hermitian stiffness matrix associated with~\cref{eq:pde}. The matrices $\bM$ and $\bK$ are sparse and admit efficient precomputed factorizations.

The solution to~\cref{eq:discrete_pde_app} can be computed by solving the linear system
$$
    (\bA + \mu \bI)\,\bx = \mathbf{b},
$$
where $\bA = \bC^{\mathrm{T}}\,\bK^{-1}\bM\,\bK^{-1}\bC$, $\mathbf{b} = \bC^{\mathrm{T}}\,\bK^{-1}\bM\,\mathbf{z}$, $\bC$ is a (possibly permuted) Cholesky factor of $\bM$. We set the regularization parameter to $\mu = 10^{-6}\|\bA\|$. The recovered solution corresponding to this value of $\mu$ is shown in~\cref{fig:invcloak3}.

Furthermore, to avoid numerical singularity of $\bA$, we apply a shift to the operator $\bA \mapsfrom \bA + 0.01 \mu \bI$ and update $\mu \mapsfrom 0.99 \mu$, as discussed in~\cref{rmk:operator_shift}.

\subsection{Basis-less Nystr{\"o}m Preconditioning Implementation.}
We implemented Nystr{\"o}m preconditioning in a basis-less form as follows (see e.g.~\cite{derezinski2024faster}):
\begin{equation}\label{eq:nys_precond}
\bP = \mu^{-1} \left( \bI - \bA \bOmega \bC^{-1} \bC^{-\mathrm{T}} \bOmega^{\mathrm{T}}\bA \right),
\end{equation}
where $\bC$ is the Cholesky factor of the matrix
$\bOmega^\mathrm{T} \bA^2 \bOmega + \mu \bOmega^\mathrm{T} \bA \bOmega$. This preconditioner is mathematically equal to the inverse of the regularized Nystr{\"o}m approximation of $\bA$: $$\bP = \left (\bA \bOmega (\bOmega^\mathrm{T} \bA \bOmega)^{-1} \bOmega^{\mathrm{T}}\bA + \mu \bI \right )^{-1}.$$

\subsection{Limited-memory Kernel Ridge Regression} 
First, to avoid numerical singularity of $\bA$, we shift operators $\bA \mapsfrom \bA + 0.01 \mu \bI$ and update the parameters $\mu \mapsfrom 0.99 \mu$ (see~\cref{rmk:operator_shift}).

Second, in some cases the matrix $\bW = \bOmega^\mathrm{T} \bA^2 \bOmega + \mu \bOmega^\mathrm{T} \bA \bOmega$ needed for obtaining a Nystr{\"o}m preconditioner~\cref{eq:nys_precond} is not numerically spd, which prevents a valid Cholesky factorization. In such cases, we define $\bC$ as the Cholesky factor of $\bW + \alpha \bI$, where $\alpha$ is the smallest value in the set ${10^{-14}|\bW|,10^{-13}|\bW|,\hdots,10^{-1}|\bW|}$ that provides a successful factorization.

\section{Additional Aspects of C-RandRAND} \label{appendix_d_0}
\subsection{Eigenvalue Surrogate for C-RandRAND} \label{surrogate}
For the cases when the estimate $F$ of $\lambda_{\min}(\bA_\mu)$ is not available, one may choose the coefficient $\rho$ in~\cref{eq:c_randrand_tau} apriori as a small factor like $0.5$ and achieve quasi-optimal range deflation. Alternatively, we propose the following computable surrogate for $F$
\begin{equation}\label{eq:optimal_rho2}
F = \max \left (\|\bA_\mu^{-1} \bPi\|^2\|\bPi \bA_\mu^{-1} \bPi\|^{-1}, \|\bE\|^{-1}, \|\bE\| \|\bA_\mu^{-1} \bPi\|^2 \lambda_{\max}^{-1}(\bPi \bA_\mu \bPi \bA_\mu^{-1} \bPi) \right ),
\end{equation}
which potentially can slightly reduce the quasi-optimality constant compared to selecting $\rho$ as $0.5$. The first two terms in this expression provide a lower bound on $\lambda_{\mathrm{min}}^{-1}(\bA_\mu)$, while the third $\|\bE\|\|\bA_\mu^{-1} \bPi\|^2$ $\lambda_{\mathrm{max}}^{-1}(\bPi \bA_\mu \bPi \bA_\mu^{-1} \bPi)$ ensures that $\rho$ satisfies 
$$\rho \|\bE\|  \leq \|\bA^{-1}_\mu \bPi\|^{-1}  \leq \lambda_{l}(\bA_\mu) \leq \|(\bI - \bPi)\bA_\mu \| \approx \|\bE\|.    $$

\subsection{Improved Condition Number Bound for C-RandRAND} \label{appendix_d}

The condition number bound from~\cref{thm:cnd_bound2} can be improved in scenarios where $\|\bPi \bA_\mu(\bI - \bPi)\|/\|\bE\|$ and $\|\bPi \bA_\mu^{-1} (\bI - \bPi)\|/\|\bF\|$ are small. Such scenarios are common in practice, as $\mathrm{range}(\bA_\mu \bOmega)$ typically is roughly an invariant subspace of $\bA_\mu$.

By combining~\cref{thm:cnd_bound2_ext} with~\cref{thm:cnd_bound2}, we obtain the following approximate condition number bound for the preconditioned operator $\bB = \bP^{\frac{1}{2}} \bA_\mu \bP^{\frac{1}{2}}$ obtained via C-RandRAND:
\begin{equation} \label{eq:cnd_bound2_ext}
    \cond(\bB) \lesssim  \min \left (\max(F,  \tau^{-1}) +  \tau^{-1/2} c_1, F +  \tau^{-1}  \right ) \cdot \min \left ( \|\bE\| \max(1, \rho ) +  \tau^{1/2} c_2,  \|\bE\| (1 + \rho ) \right ),
\end{equation}
where $F$ is an estimate of $\lambda_{\mathrm{min}}^{-1}(\bA_\mu) = \|\bA_\mu^{-1}\|$ (e.g., $F = \mu^{-1}$), $\rho = \tau \lambda_{\mathrm{max}}(\bPi \bA_\mu \bPi \bA_\mu^{-1} \bPi)/\|\bE\|$, and $c_1 = \|\bC_1 + \bC_1^\mathrm{T}\|$, $c_2 = \|\bC_2 + \bC_2^\mathrm{T}\|$.

In practice, the quantity $\|\bE\|$ can be estimated using a few power iterations. The terms $\lambda_{\mathrm{max}}(\bPi \bA_\mu \bPi \bA_\mu^{-1} \bPi)$, $\|\bC_1 + \bC_1^\mathrm{T}\|$, and $\|\bC_2 + \bC_2^\mathrm{T}\|$ can either be computed explicitly from a QR factorization of $\bA_\mu \bOmega$ or estimated via power iterations, using the identity $\bA_\mu^{-1} \bPi = \bOmega (\bA_\mu \bOmega)^\dagger$. Furthermore, we found empirically that $c_1 = \|\bC_1 + \bC_1^\mathrm{T}\|$ and $c_2 = \|\bC_2 + \bC_2^\mathrm{T}\|$ can be sufficiently well estimated by the following easily computable upper bounds
\begin{equation}\label{eq:C_est}
\|\bC_1 + \bC_1^\mathrm{T}\| \leq 2\|(\bI - \bPi)\bA_\mu^{-1} \bPi \| \lambda_{\min}(\bPi \bA^{-1}_{\mu} \bPi)^{1/2},~~~\|\bC_2 + \bC_2^\mathrm{T}\| \leq 2\|(\bI - \bPi)\bA_\mu \bPi \| \|\bPi \bA^{-1}_{\mu} \bPi\|^{1/2}.
\end{equation}
The coefficient $\tau$ in C-RandRAND can then be selected to minimize the bound that results by combining~\cref{eq:C_est} and~\cref{eq:cnd_bound2_ext}. We use this way of selecting $\tau$ in the empirical evaluation of C-RandRAND throughout~\cref{emp_eval}. Since $\bPi \bA^{-1}_{\mu} \bPi$ represents the pseudo-inverse of the Nystr{\"o}m approximation of $\bA_{\mu}$ (see~\cref{Nystrom_vs_CRandRAND}) the quantity $\lambda_{\min}(\bPi \bA^{-1}_{\mu} \bPi)$ can be computed as the reciprocal of the norm of Nystr{\"o}m approximation of $\bA_{\mu}$. 

If computing~\cref{eq:C_est} is undesirable we found that a similar performance can be achieved by selecting $c_1$ and $c_2$ in~\cref{eq:cnd_bound2_ext} as
$$c_1= \max(F^{1/2},  \tau^{-1/2}),~~~\text{and } c_2= \|\bE\|^{1/2} \max(1, \rho^{1/2} ).$$

\begin{proposition} 
\label{thm:cnd_bound2_ext} 
    Let $\bP$ be a C-RandRAND preconditioner
    $$\bP = (\bI - \bPi) + \tau \bPi \bA_\mu^{-1} \bPi$$
    for spd matrix $\bA_\mu$. Denote $\bE = (\bI-\bPi)\bA_\mu(\bI-\bPi)$, $\bF = (\bI-\bPi)\bA^{-1}_\mu(\bI-\bPi)$, $ \bC_1 = (\bI-\bPi) (\bA_\mu^{-1} \bPi) (\bPi \bA_\mu^{-1} \bPi)^{-1/2}$,$\bC_2 = (\bI-\bPi) (\bA_\mu \bPi) (\bPi \bA_\mu^{-1} \bPi)^{1/2}$.  The eigenvalues of the preconditioned operator $\bB = \bP^{1/2} \bA_\mu \bP^{1/2}$ are bounded by
        \begin{equation} \label{eq:c_randrand_eigenvalues2_ext}
            \left (\max(\lambda^{-1}_{\mathrm{min}}(\bA_\mu),  \tau^{-1}) +  \tau^{-1/2} \| \bC_1+\bC_1^\mathrm{T}\| \right  )^{-1} \leq \lambda_i(\bB) \leq \left ( \max(\|\bE\|,  \lambda_{\mathrm{max}}(\bPi \bA_\mu \bPi \bA_\mu^{-1} \bPi)) +  \tau^{1/2} \| \bC_2+\bC_2^\mathrm{T}\| \right ).
        \end{equation}
        \begin{proof}
            We have, 
            \begin{align*}
            \|\bP^{1/2} \bA_\mu \bP^{1/2}\| &= \|\bE + \tau (\bPi \bA_\mu^{-1} \bPi)^{1/2} \bA_\mu (\bPi \bA_\mu^{-1} \bPi)^{1/2} + \tau^{1/2}(\bC_2 + \bC_2^\mathrm{T})\|   \\
            &\leq \|\bE + \tau(\bPi \bA_\mu^{-1} \bPi)^{1/2} \bA_\mu (\bPi \bA_\mu^{-1} \bPi)^{1/2}\| + \tau^{1/2} \|\bC_2 + \bC_2^\mathrm{T}\| \\
            &\leq \max(\|\bE\|, \tau \|(\bPi \bA_\mu^{-1} \bPi)^{1/2} \bA_\mu (\bPi \bA_\mu^{-1} \bPi)^{1/2}\|) + \tau^{1/2}\|\bC_2 + \bC_2^\mathrm{T}\| \\
            & = \max(\|\bE\|, \lambda_{\mathrm{max}}(\bPi \bA_\mu \bPi \bA_\mu^{-1} \bPi)) + \tau^{1/2} \|\bC_2 + \bC_2^\mathrm{T}\|,
            \end{align*}
            where we used the fact that for any symmetric operators $\bM$, $\bN$ satisfying $\bM\bN = \bnull$, it holds  $\|\bM+\bN\| = \max(\|\bM\|,\|\bN\|)$. 

            Similarly, 
            \begin{align*}
            \|\bP^{-1/2} \bA_\mu^{-1} \bP^{-1/2}\| &= \|\bF + \tau^{-1} (\bPi \bA_\mu^{-1} \bPi)^{-1/2} \bA_\mu^{-1} (\bPi \bA_\mu^{-1} \bPi)^{-1/2} + \tau^{-1/2}(\bC_1 + \bC_1^\mathrm{T})\|   \\
            &\leq \|\bF + \tau^{-1} (\bPi \bA_\mu^{-1} \bPi)^{-1/2} \bA_\mu^{-1} (\bPi \bA_\mu^{-1} \bPi)^{-1/2}\| + \tau^{-1/2} \|\bC_1 + \bC_1^\mathrm{T}\| \\
            &\leq \max(\|\bF\|, \tau^{-1} \|(\bPi \bA_\mu^{-1} \bPi)^{-1/2} \bA_\mu^{-1} (\bPi \bA_\mu^{-1} \bPi)^{-1/2}\|) + \tau^{-1/2}\|\bC_1 + \bC_1^\mathrm{T}\| \\
            & = \max(\|\bF\|, \tau^{-1}) + \tau^{-1/2}\|\bC_1 + \bC_1^\mathrm{T}\|. 
            \end{align*}
            
        \end{proof}
\end{proposition}

\end{document}